\pdfoutput=1
\RequirePackage{ifpdf}
\ifpdf 
\documentclass[pdftex]{sigma}
\else
\documentclass{sigma}
\fi

\numberwithin{equation}{section}

\newtheorem{Theorem}{Theorem}[section]
\newtheorem{Corollary}[Theorem]{Corollary}
\newtheorem{Lemma}[Theorem]{Lemma}
\newtheorem{Proposition}[Theorem]{Proposition}
 { \theoremstyle{definition}
\newtheorem{Definition}[Theorem]{Definition}
\newtheorem{Example}[Theorem]{Example}
\newtheorem{Remark}[Theorem]{Remark} }

\usepackage[shortlabels]{enumitem}
\usepackage{mathdots}
\usepackage[normalem]{ulem}
\usepackage{chngcntr}

\usepackage{mathrsfs}
\usepackage{euscript}
\usepackage[all,cmtip]{xy}
\usepackage{xstring}
\usepackage{slashed}
\usepackage{tensor}
\usepackage{mathtools}

\usepackage{tikz-cd}

\def\mathscr#1{\EuScript{#1}}

\entrymodifiers={+!!<0pt,\fontdimen22\textfont2>}

\def\N {\mathbb{N}}
\def\Z {\mathbb{Z}}

\def\C {\mathbb{C}}

\def\hc#1{\mathrm{h}_{#1}}
\def\h {\mathrm{H}}

\def\subset{\subseteq}

\newcommand{\incl}{\hookrightarrow}

\newlength{\myeqm} 
\newlength{\myeqn} 
\def\brackets#1{\IfStrEq{#1}{-}{}{(#1)}}
\def\subindex#1{\IfStrEq{#1}{-}{}{_{#1}}}

\newcommand{\alxydim}[2]{\begin{aligned}\xymatrix#1{#2}\end{aligned}}



\newlength{\myl}
\newcommand\sheaf[1]{\unitlength 0.1mm
 \settowidth{\myl}{$#1$}
 \addtolength{\myl}{-0.8mm}
 \begin{picture}(0,0)(0,0)
 \put(2,0){\text{\underline{\hspace{\myl}}}}
 \end{picture}#1\hspace{-0.15mm}}

\def\ddt#1#2#3{\left.\frac{\mathrm{d}^{\IfStrEq{#1}{1}{}{#1}}}{\mathrm{d}#2}\IfStrEq{#2}{#3}{\right.}{\right|_{#3}}}

\def\liegrpd{\mathscr{L}\text{ie}\mathscr{G}\text{rpd}}

\def\hom{\mathcal{H}\!{\rm om}}

\def\act#1#2{#1/\!\!/#2}

\newlength{\widthtmp}
\def\length#1{\settowidth{\widthtmp}{#1}\the\widthtmp}

\def\ttimes#1#2{\hspace{-0.15em}\tensor[_{#1}]{\times}{_{#2}}}

\newcommand{\ie}{i.e., }

\newcommand{\U}{\operatorname{U}}

\newcommand{\pr}{\operatorname{pr}}
\newcommand{\lact}{\triangleright}

\newcommand{\id}{\operatorname{id}}

\newcommand{\Hom}{\mathscr{H}\mathrm{om}}

\newcommand{\Aut}{\mathrm{Aut}}

\newcommand{\VectBdl}{\mathscr{V}\mathscr{B}\mathrm{dl}}
\newcommand{\LineBdl}{\mathscr{L}\mathscr{B}\mathrm{dl}}
\newcommand{\sVectBdl}{\mathrm{s}\VectBdl}

\newcommand{\twoLineBdl}{2\LineBdl}
\newcommand{\twoLineBdltriv}{2\LineBdl_{\mathrm{triv}}}
\newcommand{\stwoLineBdl}{\mathrm{s}2\LineBdl}

\newcommand{\stwoLineBdltriv}{\mathrm{s}2\LineBdl_{\mathrm{triv}}}
\newcommand{\stwoLineBdlref}{\mathrm{s}2\LineBdl^{\mathrm{ref}}}

\newcommand{\Grb}{\mathscr{G}\mathrm{rb}}
\newcommand{\Grbtriv}{\mathscr{G}\mathrm{rb}_{\mathrm{triv}}}
\newcommand{\sGrb}{\mathrm{s}\Grb}
\newcommand{\sGrbtriv}{\mathrm{s}\Grbtriv}
\newcommand{\hGrb}{\mathrm{h}\Grb}

\newcommand{\sGrbref}{\mathrm{s}\Grb^{\mathrm{ref}}}

\newcommand{\Alg}{\Incl\mathrm{lg}}

\newcommand{\sAlg}{\mathrm{s}\Alg}

\newcommand{\sAlgbi}{\sAlg^{\mathrm{bi}}}

\newcommand{\csAlg}{\mathrm{cs}\Incl\mathrm{lg}}
\newcommand{\cssAlg}{\mathrm{cs\text{-}s}\Alg}

\newcommand{\sAlgBdlbi}{\sAlg\mathscr{B}\mathrm{dl}^{\mathrm{bi}}}

\newcommand{\csAlgBdl}{\csAlg\mathscr{B}\mathrm{dl}}
\newcommand{\cssAlgBdl}{\cssAlg\mathscr{B}\mathrm{dl}}
\newcommand{\cssAlgBdlbi}{\cssAlg\mathscr{B}\mathrm{dl}^{\mathrm{bi}}}
\newcommand{\csAlgBdlbi}{\csAlg\mathscr{B}\mathrm{dl}^{\mathrm{bi}}}

\newcommand{\MultGrb}{\mathscr{M}\mathrm{ult}\Grb}

\def\sExt{\normalfont\text{s}\mathscr{E}\text{xt}}
\def\RsExt{\normalfont\mathscr{R}\text{s}\mathscr{E}\text{xt}}
\def\Mou{\operatorname{Mou}}

\def\Gr{\mathrm{Gr}}

\def\Rliegrpd{\normalfont\mathscr{R}\liegrpd}
\def\Grliegrpd{\normalfont\mathscr{G}\mathrm{r}\liegrpd}

\def\RVectBdl{\mathscr{R}\VectBdl}

\newcommand{\Incl}{\mathscr{A}}

\def\quand{\qquad\text{and}\qquad}

\usetikzlibrary{patterns}
\usetikzlibrary{decorations.pathreplacing,angles,quotes}
\usetikzlibrary{decorations.text}

\begin{document}

\allowdisplaybreaks

\newcommand{\arXivNumber}{2502.18102}

\renewcommand{\PaperNumber}{059}

\FirstPageHeading

\ShortArticleName{Real Twistings are 2-Line Bundles}

\ArticleName{Real Twistings are 2-Line Bundles}

\Author{Tim L\"UDERS~$^{\rm a}$, Lynn OTTO~$^{\rm b}$ and Konrad WALDORF~$^{\rm b}$}

\AuthorNameForHeading{T.~L\"uders, L.~Otto and K.~Waldorf}

\Address{$^{\rm a)}$~Universit\"at Wien, Fakult\"at f\"ur Physik, Boltzmanngasse 5, 1090 Wien, Austria}
\EmailD{\mail{tim.lueders@univie.ac.at}}

\Address{$^{\rm b)}$~Universit\"at Greifswald, Institut f\"ur Mathematik und Informatik,\\
\hphantom{$^{\rm b)}$}~Walther-Rathenau-Str.~47, 17487 Greifswald, Germany}
\EmailD{\mail{lynn.otto@uni-greifswald.de}, \mail{konrad.waldorf@uni-greifswald.de}}

\ArticleDates{Received May 13, 2025, in final form June 04, 2026; Published online June 18, 2026}

\Abstract{We construct and study a~bicategory of super 2-line bundles over graded Lie groupoids, providing a~unified framework for geometric models of twistings of (real) K-theory. The core of our work is to exhibit a~wide range of models from the literature as special cases, among them several variants of bundle gerbes (real/equivariant/Jandl), Freed--Moore's twisted groupoid extensions, Freed--Hopkins--Teleman's K-theory twistings, Moutuou's real twistings, Freed's invertible algebra bundles, and Distler--Freed--Moore's orientifold twistings.}


\Classification{19L50; 22A22; 53C08; 18N10}

\section{Introduction}

Over the last few decades, twisted K-theory has attracted considerable interest, much of it driven by applications in theoretical physics, e.g., \cite{Bouwknegt2000,Diaconescu2000ADO,Freeda,kapustin1,Maldacena2001DBraneIA,witten4}. Twisted cohomology theories in general can be approached via homotopy-theoretic methods, and twisted K-theory has additionally elegant operator-algebraic descriptions, e.g., \cite{Kasparov1979TheKI,Rosenberg1989}.
However, when one wishes to make close contact with physical models, differential-geometric realizations of the twistings and the twisted K-theory classes become very useful.

It is well known that twistings are often easier to describe geometrically than the corresponding twisted K-theory classes, which can be difficult to treat in a~strictly finite-dimensional, smooth framework. The focus of this article lies in differential-geometric models for the twistings, and, more specifically, on their typical low-degree part located in the \v{C}ech cohomology groups
\begin{equation}
\label{classification}
\check\h^0(X,\Z_2) \times\check \h^1(X,\Z_2) \times \check\h^2(X,\sheaf{\C^{\times}}).
\end{equation}
 Here, the first factor represents the $\Z_2$-grading of twisted K-theory (for $\mathrm{KO}$-theory it would have values in $\Z_8$). We remark at this point that we work without hermitian structures, which makes no difference on the cohomological level since $ \check\h^2(X,\sheaf{\C^{\times}})\cong \check\h^2(X,\sheaf{\U(1)})\cong \h^3(X,\Z)$.

A range of geometric models for these twists already exists. Classical work by Donavan and Karoubi
\cite{DK70} characterizes the torsion components of~\eqref{classification} via central simple super algebra bundles. The non-torsion part -- i.e., general classes in
$\check \h^2(X,\sheaf{\C^{\times}})$
-- is elegantly realized by bundle gerbes~\cite{bouwknegt1,carey2,murray}, objects that naturally appear in conformal field theory and higher gauge theory~\cite{gawedzki1,Mickelsson2004}. More recently, nice geometric objects with precisely the classification~\eqref{classification} have been introduced: \emph{super $2$-line bundles}, a~categorified version of super line bundles
\cite{Kristel2020,Mertsch2020}. Super 2-line bundles can be viewed as a~combination of bundle gerbes and central simple algebra bundles, and contain both structures as special cases.
The basic idea is that a~super 2-line consists of locally defined central simple super algebra bundles, which are glued together along Morita equivalences in a~consistent way. Super 2-line bundles over $X$ form a~bicategory~$\stwoLineBdl(X)$, and satisfy a~gluing law (they form a~2-stack).
Moreover, they have a~symmetric monoidal structure
capturing precisely the ring structure of twistings in K-theory. Thus, super 2-line bundles are a~complete and nice differential-geometric model for twistings on manifolds~$X$.

Twisted K-theory arises naturally in string theory when D-branes wrap submanifolds of the target space and B-fields necessitate a~``twist'' of Chan--Paton gauge bundles~\cite{gawedzki1,kapustin1,witten4}. Further generalizations of twisted K-theory become relevant when the target space admits additional ``orbifold'' symmetries. These may be implemented by group actions, or more generally, by replacing target manifolds $X$ by Lie groupoids $\Gamma$. For example, equivariant versions of bundle gerbes -- pioneered by Gaw\c{e}dzki--Reis and Meinrenken~\cite{gawedzki2, meinrenken1} -- play a~key role in describing WZW-models on non-simply connected groups.
Likewise, the Verlinde ring of loop group representations can be identified with a~twisted K-theory associated to the action Lie groupoid $\act GG$, as shown in the seminal work of Freed, Hopkins, and Teleman~\cite{Freed2011a,Freed2011,Freed2012LoopGA}. They introduced a~new model for twistings over groupoids, based on central extensions of groupoids.

The full set of twistings (again, below degree three) on a~Lie groupoid $\Gamma$ is classified by the Lie groupoid \v Cech cohomology groups
\begin{equation}
\label{classification-2}
\check\h^0(\Gamma,\Z_2) \times\check \h^1(\Gamma,\Z_2) \times \check\h^2(\Gamma,\sheaf{\C^{\times}}),
\end{equation}
of which equivariant bundle gerbes realize the last factor, while
 the model of Freed--Hopkins--Teleman realizes the last two factors (the first factor, the grading, is added manually).

A further layer of complexity arises in ``orientifold'' setups, where certain symmetries must act by complex conjugation rather than in the usual $\C$-linear way. This leads to so-called \emph{graded Lie groupoids} $(\Gamma,\phi)$, in which each morphism $\gamma$ carries a~$\Z_2$-label $\phi(\gamma)$ indicating if $\gamma$ acts linearly or anti-linearly. Graded Lie groupoids contain the previously considered cases of manifolds and Lie groupoids, as well as further configurations such as manifolds with an involution (``real spaces'') and Lie groupoids with an involution functor (``real groupoids''). One can extend the notion of Lie groupoid cohomology to such graded situations, thereby obtaining classification data of the form
\begin{equation}
\label{classification-3} \check \h^0(\Gamma,\Z_2)\times\check \h^1(\Gamma,\Z_2)\times\check \h^2((\Gamma,\phi),\sheaf{\C^{\times}}). \end{equation}
We remark that the grading $\phi$ has necessarily no effect on the first two factors.
Variations of bundle gerbes (so-called Jandl gerbes) \cite{gawedzki8,schreiber1}, real~\cite{mohamedmoutuou2012} or ``$\phi$-twisted'' groupoid extensions~\cite{Freeda}, and further constructions~\cite{Distler2011,distler2010spin,gomi2021freedmoorektheory} have been developed to handle twistings in these graded contexts. They appear, for instance, also in the classification of topological insulators, where Bloch bundles over the Brillouin torus become twisted depending on symmetry types~\cite{Freeda}.

In this paper, we propose a~new, comprehensive differential-geometric model for twistings over graded Lie groupoids, one which unifies \emph{all} of the aforementioned versions. Specifically, we extend super 2-line bundles -- which already capture all twistings on ordinary manifolds, \eqref{classification} -- first to ordinary Lie groupoids and then to graded Lie groupoids. The passage from manifolds to Lie groupoids follows a~standard homotopy-limit construction applicable to any presheaf (of sets, categories, or higher categories) on the site of manifolds. The extension to graded Lie groupoids uses a~further ``involution'' ingredient, but remains formally similar.

We thus arrive in Definition~\ref{The bicategory of super 2-line bundles definition} at a~bicategory $\stwoLineBdl(\Gamma,\phi)$ of super 2-line bundles over a~graded Lie groupoid $(\Gamma,\phi)$ that reduces -- for a~trivial grading $\phi=1$ -- to the canonical homotopy limit construction of super 2-line bundles over Lie groupoids, and moreover -- for a~Lie groupoid with only identity morphisms -- to the ordinary bicategory of super 2-line bundles over manifolds. Up to isomorphism, super 2-line bundles over graded Lie groupoids, over Lie groupoids, and over manifolds are classified by the cohomologies~\eqref{classification-3}, \eqref{classification-2}, and~\eqref{classification}, respectively (see Theorem~\ref{classification-of-graded-equivariant-two-line-bundles}). A crucial feature of our bicategory $\stwoLineBdl(\Gamma,\phi)$ is that it is invariant under Morita equivalences of graded Lie groupoids (see Theorem~\ref{descent-for-graded-equivariant-2-line-bundles}).

After setting up our formalism, we show how super 2-line bundles recover or generalize several well-known types of twistings from the literature; the most important ones are
\begin{itemize}\itemsep=0pt
\item
Jandl and other real or equivariant/anti-equivariant bundle gerbes~\cite{gawedzki8,nikolaus2,schreiber1}, see Theorem~\ref{comparison-equivariant-jandl-gerbes}. This result comes essentially from the fact that already over manifolds bundle gerbes are a~subclass of 2-line bundles.
\item
Freed--Moore's $\phi$-twisted groupoid extensions~\cite{Freeda}, see Theorem~\ref{The main theorem of the Freed-Moore chapter}. This result is obtained by a~detailed comparison of the relevant definitions.

\item
 Freed--Hopkins--Teleman's twistings by central extensions of groupoids~\cite{Freed2011a}, see Theorem~\ref{The main theorem of the Freed Hopkins Teleman chapter}. Our work here is to compare spans and refinements of weak equivalences of Lie groupoids (used by Freed--Hopkins--Teleman) with spans and refinements of surjective submersions (used in our definition of super 2-line bundles).

\item
Moutuou's real groupoid picture of twistings~\cite{mohamedmoutuou2012}, extending the concept of real gerbes and real central extensions, see Theorem~\ref{main-theorem-moutuou}.
\end{itemize}
Hence, super 2-line bundles provide a~single, coherent framework that not only reproduces these and many further special cases but clarifies their underlying geometry and consistency conditions.
We hope that this unification will be useful for researchers working on twisted K-theory, higher geometry, and related applications in physics.

In order to keep the present paper more transparent, we do not consider additional structures (e.g., hermitian metrics or connections) on our super 2-line bundles. We plan to include such structures in forthcoming work. Likewise, we devote little attention to the symmetric monoidal structure on super 2-line bundles over graded Lie groupoids. Finally, although we
carry out some comparisons with the twisted vector bundles belonging to each of the models, we do not delve too far into the discussion of the corresponding twisted K-theory groups.

This article is organized in the following way. Section~\ref{main-super-2-line-bundles} contains a~comprehensive introduction to super 2-line bundles. In Section~\ref{extension-to-graded-lie-groupoids}, we set up the formalism of super 2-line bundles over graded Lie groupoids, and provide related results. We present in Section~\ref{twistings} our list of twistings from the literature that we recover. Section~\ref{section-comparison-work} consists of separate sections in which detailed comparison work from the list of Section~\ref{twistings} is carried out. Finally, we include an appendix about Lie groupoids and weak equivalences, which appear at several places.

\section{Super 2-line bundles}
\label{main-super-2-line-bundles}

In this section, we set the stage for the fundamental objects studied in this paper, \emph{$2$-line bundles}, i.e., categorified line bundles. In the first two subsections, we recall material derived from a~broader discussion of 2-vector bundles in~\cite{Kristel2020}. In the third subsection, we introduce the notion of complex conjugate 2-line bundles.
\subsection{Definitions}
\label{section-super-2-line-bundles}

We recall from~\cite{Kristel2020,Mertsch2020} the bicategory of super 2-line bundles $\stwoLineBdl(X)$ over a~smooth mani\-fold~$X$. Throughout, we use the field of complex numbers $\C$ as our ground field.
\begin{Definition}
\label{super-2-line-bundle}
A \emph{super $2$-line bundle} over $X$ is a~quadruple $\mathscr{L}=(\pi,\mathcal{A},\mathcal{M},\mu)$ consisting of a~surjective submersion $\pi\colon Y \to X$, a~central simple super algebra bundle $\mathcal{A}$ over $Y$,
an invertible bimodule bundle $\mathcal{M}$ over $Y^{[2]}$ whose fibre $\mathcal{M}_{y_1,y_2}$ over a~point $(y_1,y_2)\in Y^{[2]}$ is an $\mathcal{A}_{y_2}$-$\mathcal{A}_{y_1}$-bimodule, and
an invertible even intertwiner $\mu$ of bimodule bundles over $Y^{[3]}$, which restricts over each point $(y_1,y_2,y_3)\in Y^{[3]}$ to an $\mathcal{A}_{y_3}$-$\mathcal{A}_{y_1}$-bimodule intertwiner
\begin{equation*}
\mu_{y_1,y_2,y_3}\colon\ \mathcal{M}_{y_2,y_3} \otimes_{\mathcal{A}_{y_2}} \mathcal{M}_{y_1,y_2} \to \mathcal{M}_{y_1,y_3}.
\end{equation*}
This structure is subject to the condition that
 $\mu$ is
 associative, \ie the diagram
\begin{align*}
\xymatrix@C=5em{\mathcal{M}_{y_3,y_4} \otimes_{\mathcal{A}_{y_3}} \mathcal{M}_{y_2,y_3} \otimes_{\mathcal{A}_{y_2}} \mathcal{M}_{y_1,y_2} \ar[d]_{\id\otimes \mu_{y_1,y_2,y_3}} \ar[r]^-{\mu_{y_2,y_3,y_4} \otimes \id} & \mathcal{M}_{y_2,y_4} \otimes_{\mathcal{A}_{y_2}} \mathcal{M}_{y_1,y_2} \ar[d]^{\mu_{y_1,y_2,y_4}} \\ \mathcal{M}_{y_3,y_4} \otimes_{\mathcal{A}_{y_3}} \mathcal{M}_{y_1,y_3} \ar[r]_-{\mu_{y_1,y_3,y_4}} & \mathcal{M}_{y_1,y_4} }
\end{align*}
is commutative for all $(y_1,y_2,y_3,y_4)\in Y^{[4]}$.
\end{Definition}

\begin{Remark}\quad
\begin{itemize}\itemsep=0pt
\item
Definition~\ref{super-2-line-bundle} is almost the definition of a~super 2-\emph{vector} bundle given in~\cite{Kristel2020}, just that we require here that the algebra bundle has central simple fibres, corresponding to the condition of having rank one for ordinary line bundles.
\item
For precise definitions of algebra bundles and bimodule bundles we refer to~\cite{Kristel2022}. In particular, all algebras are associative and unital. Due to the restriction to central simple algebra bundles, the ``insidious'' problems in the construction of the relative tensor product of bimodule bundles that have been observed in~\cite{Kristel2022} do not appear here.

\item
If $\mathscr{L}=(\pi,\mathcal{A},\mathcal{M},\mu)$ is a~super 2-line bundle, and if $\Delta\colon Y \to Y^{[2]}$ denotes the diagonal map, then there exists a~canonical invertible intertwiner $\Delta^{*}\mathcal{M} \cong \mathcal{A}$ of $\mathcal{A}$-$\mathcal{A}$-bimodule bundles over $Y$. Second, if $s\colon Y^{[2]} \to Y^{[2]}$ swaps the factors, then $s^{*}\mathcal{M}$ is inverse to $\mathcal{M}$.

\end{itemize}
\end{Remark}
We suppose that $\mathscr{L}_1=(\pi_1,\mathcal{A}_1,\mathcal{M}_1,\mu_1)$ and $\mathscr{L}_2=(\pi_2,\mathcal{A}_2,\mathcal{M}_2,\mu_2)$ are super 2-line bundles over $X$.
\begin{Definition}
\label{1-morphism-between-2-line-bundles}
A \emph{$1$-morphism}
$\mathscr{P}\colon \mathscr{L}_1 \to \mathscr{L}_2$
is a~triple $\mathscr{P}=(\zeta,\mathcal{P},\phi)$ consisting of a~surjective submersion $\zeta\colon Z \to Y_1 \times_X Y_2$, a~bimodule bundle $\mathcal{P}$ over $Z$, whose fibre $\mathcal{P}_z$ over a~point~${z\in Z}$ with $\zeta(z)=:(y_1,y_2)$ is an $(\mathcal{A}_2)_{y_2}$-$(\mathcal{A}_1)_{y_1}$-bimodule, and
an invertible even intertwiner $\phi$ of bimodule bundles over $Z^{[2]}$, which restricts over a~point $(z,z')\in Z^{[2]}$ with $\zeta(z')=:(y_1',y_2')$ to an $(\mathcal{A}_2)_{y_2'}$-$(\mathcal{A}_1)_{y_1}$-intertwiner
\begin{equation*}
\phi_{z,z'}\colon\ \mathcal{P}_{z'} \otimes_{(\mathcal{A}_1)_{y_1'}} (\mathcal{M}_1)_{y_1,y_1'} \to (\mathcal{M}_2)_{y_2,y_2'} \otimes_{(\mathcal{A}_2)_{y_2}} \mathcal{P}_{z} .
\end{equation*}
This intertwiner is a~``homomorphism'' with respect to the intertwiners $\mu_1$ and $\mu_2$, in the sense that the diagram
\begin{equation}\label{diag:1Morphisms}
\begin{gathered}
\xymatrix@C=6em{
\mathcal{P}_{z^{\prime\prime}} \otimes_{(\mathcal{A}_1)_{y_1^{\prime\prime}}}(\mathcal{M}_1)_{y_1',y_1''} \otimes_{(\mathcal{A}_1)_{y_1'}} (\mathcal{M}_1)_{y_1,y_1'} \ar[r]^-{ \id \otimes (\mu_1)_{y_1,y_1',y_1''}} \ar[d]_{\phi_{z^\prime,z^{\prime\prime}}\otimes \id} &
\mathcal{P}_{z^{\prime\prime}} \otimes_{(\mathcal{A}_1)_{y_1''}} (\mathcal{M}_1)_{y_1,y_1''} \ar[dd]^{\phi_{z,z''}} \\
 (\mathcal{M}_2)_{y_2',y_2''} \otimes_{(\mathcal{A}_2)_{y_2'}} \mathcal{P}_{z'} \otimes_{(\mathcal{A}_1)_{y_1'}} (\mathcal{M}_1)_{y_2,y_2'} \ar[d]_{\id \otimes \phi_{z,z'}} & \\
 (\mathcal{M}_2)_{y_2',y_2''} \otimes_{(\mathcal{A}_2)_{y_2'}} (\mathcal{M}_2)_{y_2,y_2'} \otimes_{(\mathcal{A}_2)_{y_2}} \mathcal{P}_{z} \ar[r]_-{(\mu_2)_{y_2,y_2',y_2''} \otimes \id} &
 (\mathcal{M}_2)_{y_2,y_2''} \otimes_{(\mathcal{A}_2)_{y_2}} \mathcal{P}_{z}
 }
\end{gathered}
\end{equation}
is commutative for all $(z,z',z'')\in Z^{[3]}$, where $\zeta(z'')=: (y_1'',y_2'')$.
\end{Definition}

Finally, we recall how the 2-morphisms are defined. Consider two 1-morphisms between the same super 2-line bundles,
\begin{equation*}
\xymatrix{\mathscr{L}_1 \ar@/^1pc/[r]^{\mathscr{P}}\ar@/_1pc/[r]_{\mathscr{P}'} & \mathscr{L}_2,}
\end{equation*}
with all structure denoted and labelled as above.
For abbreviation, we set $Y_{12} := Y_1 \times_X Y_2$.

\begin{Definition}
A \emph{$2$-morphism} $\mathscr{P} \Rightarrow \mathscr{P}'$ is an equivalence class of pairs $(\rho,\varphi)$ consisting of a~surjective submersion $\rho\colon W \to Z \times_{Y_{12}} Z'$, and
an intertwiner $\varphi$ of bimodule bundles over $W$ that restricts over a~point $w\in W$ with $\rho(w)=:(z,z')$ to an intertwiner
$
\varphi_w\colon \mathcal{P}_{z} \to \mathcal{P}'_{z'}
$
of $(\mathcal{A}_1)_{y_1}$-$(\mathcal{A}_2)_{y_2}$-bimodules, where $\zeta(z)=\zeta'(z')=:(y_1,y_2)$. This intertwiner is subject to the condition that it commutes with the intertwiners $\phi$ and $\phi'$, in the sense that the diagram
\begin{equation*}
\xymatrix@C=4em{
\mathcal{P}_{\tilde z} \otimes_{(\mathcal{A}_1)_{\tilde y_1}} (\mathcal{M}_1)_{y_1,\tilde y_1} \ar[d]_{\varphi_{\tilde w} \otimes \id} \ar[r]^-{\phi_{z,\tilde z}} & (\mathcal{M}_2)_{y_2,\tilde y_2} \otimes_{(\mathcal{A}_2)_{y_2}} \mathcal{P}_{z} \ar[d]^{\id \otimes \varphi_w}
 \\
 \mathcal{P}'_{\tilde z'} \otimes_{(\mathcal{A}_1)_{\tilde y_1'}} (\mathcal{M}_1)_{y_1',\tilde y_1'} \ar[r]_-{\phi'_{z',\tilde z'}} &(\mathcal{M}_2)_{y_2',\tilde y_2'} \otimes_{(\mathcal{A}_2)_{y'_2}} \mathcal{P}'_{z'}
 }
\end{equation*}
is commutative for all $(w,\tilde w) \in W \times_X W$, where $\rho(w)=(z,z')$, $\rho(\tilde w)=(\tilde z, \tilde z')$.
Two pairs $(\rho,\varphi)$ and $(\rho',\varphi')$ are equivalent if the pullbacks of $\varphi$ and $\varphi'$ coincide over \smash{$W \times_{Z \times_{Y_{12}} Z'} W'$}.
\end{Definition}

The structure collected above forms the bicategory $\stwoLineBdl(X)$ of super 2-line bundles over~$X$.
We recall the following result about the invertibility of 1-morphisms and 2-morphisms \cite[Lemma~2.3.4]{Kristel2020}.

\begin{Proposition}\label{LemmaInvertibility}\quad
\begin{enumerate}[$(a)$]\itemsep=0pt
\item \label{LemmaInvertibilityA}
A $1$-morphism $\mathscr{P}=(\zeta,\mathcal{P},\phi)$ is invertible if and only if its bimodule bundle $\mathcal{P}$ is invertible in the bicategory of central simple super algebra bundles.

\item \label{LemmaInvertibilityA2}
A $1$-morphism $\mathscr{P}=(\zeta,\mathcal{P},\phi)$ has a~right $($left$)$ adjoint if and only if its bimodule bundle $\mathcal{P}$ has a~right $($left$)$ adjoint in the bicategory of central simple super algebra bundles.

\item \label{LemmaInvertibilityB}
A $2$-morphism $(\rho,\varphi)$ is invertible if and only if its intertwiner $\varphi$ is invertible.

\end{enumerate}
\end{Proposition}

We also recall the following classification result \cite[Theorem~4.4.4]{Kristel2020}.

\begin{Theorem}
\label{super-2-line-bundles}
The set of isomorphism classes of super $2$-line bundles over $X$ is in bijection with the set
\begin{equation*}
\h^0(X,\mathrm{BW}_\C) \times \h^1(X,\Z_2) \times \h^3(X,\Z),
\end{equation*}
where $\mathrm{BW}_\C\cong \Z_2$ is the Brauer--Wall group of Morita equivalence classes of central simple super algebras over $\C$.
\end{Theorem}

\begin{Remark}
\label{extracting-cocycle}
We recall how the classification of Theorem~\ref{super-2-line-bundles} comes about by extracting cocycles from a~super 2-line bundle $\mathscr{L}=(\pi,\mathcal{A},\mathcal{M},\mu)$; see \cite[Section~4]{Kristel2020} and \cite[Section~2.2]{Mertsch2020}. We choose an open cover $\mathcal{U}=(U_i)_{i\in I}$ of $X$ with simply-connected open sets $U_i$ that allow sections~${s_i\colon U_i \to Y}$.
The pullback $s_i^{*}\mathcal{A}$ is then trivializable, $s_i^{*}\mathcal{A} \cong U_i \times A_i'$, where $A_i'$ is a~central simple super algebra. We fix for each element in~$\mathrm{BW}_{\C}\cong \Z_2$ a~representative; for instance, the neutral element may be represented by $A=\C$ (purely even), and the non-trivial element may be represented by a~super algebra $A=\C \oplus \C u$, where $u$ is an odd element and~${u^2=1}$.
 Let $A_i$ denote the representative of the class of $A_i'$. Then, with $\mathcal{A}_i:=U_i \times A_i$, we obtain an invertible $s_i^{*}\mathcal{A}$-$\mathcal{A}_i$-bimodule bundle $\mathcal{N}_i$ over $U_i$. We may also view $\mathcal{N}_i$ as an isomorphism $\mathscr{N}_i\colon \mathcal{A}_i \to \mathscr{L}|_{U_i}$ in~$\stwoLineBdl(U_i)$, a~local trivialization. We obtain a~(constant) map~${\alpha_i\colon U_i \to \mathrm{BW}_{\C}}$ with value $A_i$.
On a~double overlap $U_{ij} :=U_i \cap U_j$, we consider the corresponding section $s_{ij} := (s_i,s_j)\colon U_{ij} \to Y^{[2]}$, and the $s_j^{*}\mathcal{A}$-$s_i^{*}\mathcal{A}$-bimodule bundle $s_{ij}^{*}\mathcal{M}$.
 We obtain the invertible $\mathcal{A}_j$-$\mathcal{A}_i$-bimodule bundle \smash{$\mathcal{M}_{ij} := \mathcal{N}^{-1}_{j} \otimes_{s_j^{*}\mathcal{A}} s_{ij}^{*}\mathcal{M} \otimes_{s_i^{*}\mathcal{A}} \mathcal{N}_i$} over $U_{ij}$.
 We may also view \smash{$\mathcal{M}_{ij} = \mathscr{N}_j^{-1} \circ \mathscr{N}_i\colon \mathcal{A}_i \to \mathcal{A}_j$} as the ``transition function'' for the local trivializations $\mathscr{N}_i$. We also remark that $\mathcal{M}_{ij} \cong U_{ij}\times M_{ij}$, with $M_{ij}$ an invertible $A_j$-$A_i$-bimodule. Since~$A_i$ and~$A_j$ are -- by construction -- both the fixed representative of the same Morita equivalence class, we have $A:=A_i=A_j$ and $\alpha_i = \alpha_j$. Now, $M_{ij}$ is an invertible $A$-$A$-bimodule of the central simple super algebra $A$, whose well-known classification implies that either $M_{ij}\cong A$ or~${M_{ij} \cong \Pi A}$, the grading reversal of $A$.
Writing $\varepsilon A$ to mean either $A$ (when $\varepsilon=1$) or $\Pi A$ (when $\varepsilon=-1$), we obtain a~map $\epsilon_{ij}\colon U_{ij} \to \Z_2$ such that $M_{ij}\cong \varepsilon_{ij}A$. On a~triple overlap~${U_{ijk} :=U_i \cap U_j \cap U_k}$, the isomorphism $\mu$ becomes, under the above isomorphisms, an invertible intertwiner $\mu_{ijk}'\colon U_{ijk} \times (\varepsilon_{jk}A \otimes_A \varepsilon_{ij}A) \to U_{ijk} \times \varepsilon_{ik}A$, where $A:=A_i=A_j=A_k$. Thus, we have $\varepsilon_{jk}A \otimes_A \varepsilon_{ij}A \cong (\varepsilon_{jk}\varepsilon_{ij})A$; this implies $\varepsilon_{jk}\varepsilon_{ij}=\varepsilon_{ik}$.
 Automorphisms of the $A$-$A$-bimodule~$A$ are given by multiplication with elements in $Z(A^{\times})$, and since $A$ is central, \smash{$\mu'_{ijk}$} becomes a~smooth map \smash{$\mu_{ijk}\colon U_{ijk} \to \C^{\times}$}. The associativity condition in Definition~\ref{super-2-line-bundle} implies the cocycle condition for $\mu_{ijk}$. The triple $(\alpha_i,\varepsilon_{ij},\mu_{ijk})$ is a~\v Cech cocycle representative of the class of $\mathscr{L}$ under the classification of Theorem~\ref{super-2-line-bundles} (under the usual identification $\check \h^2(X,\sheaf {\C^{\times}}) \cong \h^3(X,\Z)$).
\end{Remark}

\begin{Remark}
A non-super version of super 2-line bundles is included into our discussion as those super 2-line bundles where the odd components of all involved super algebra bundles and bimodule bundles are zero. This sets up the bicategory $\twoLineBdl(X)$ of 2-line bundles, and the classification of Theorem~\ref{super-2-line-bundles} reduces to the set $ \h^3(X,\Z)$, as the Brauer group of $\C$ is trivial.
\end{Remark}

If one restricts Definition~\ref{super-2-line-bundle} to the trivial algebra bundle with fibre $\C$, $\mathcal{A}=Y \times \C$, the bimodule bundles in Definitions~\ref{super-2-line-bundle} and~\ref{1-morphism-between-2-line-bundles} become just super vector bundles. Moreover, in Definition~\ref{super-2-line-bundle}, the requirement that the bimodule bundle $\mathcal{M}$ is invertible implies that $\mathcal{M}$ has rank one, i.e., it is a~super line bundle.
Thus, above definitions reduce to the ones of (super) bundle gerbes~\cite{murray}, and we obtain the following result, see \cite[Proposition~3.2.2]{Kristel2020}.

\begin{Proposition}
\label{embedding-of-bundle-gerbes}
There are, respectively, a~fully faithful functor and an equivalence of bicategories,\begin{equation*}
\sGrb(X) \to \stwoLineBdl(X)
\qquad
\Grb(X) \stackrel{\sim}\to \twoLineBdl(X)
\end{equation*}
that embed $($super$)$ bundle gerbes into $($super$)$ {\rm2}-line bundles. Under the classification of Theorem~{\rm\ref{super-2-line-bundles}}, their images correspond to the subsets $\h^1(X,\Z_2) \times \h^3(X,\Z)$ and $\h^3(X,\Z)$, respectively.
\end{Proposition}

We denote by $\cssAlgBdlbi(X)$ the bicategory of central simple super algebra bundles over~$X$, with bimodules as the 1-morphisms.
The definition of the bicategory $\stwoLineBdl(X)$ recalled above was obtained in~\cite{Kristel2020} by applying the plus construction of~\cite{nikolaus2} to the presheaf $\cssAlgBdlbi$. Recall that the classical plus construction of Grothendieck assigns to a~presheaf $\mathcal{F}$ a~sheaf $\mathcal{F}^+$, its sheafification. The sheafification functor $\mathcal{F} \mapsto \mathcal{F}^+$ can be characterized abstractly as the left adjoint to the inclusion functor of sheaves into presheaves, or can be described by a~concrete model. In the concrete model, $\mathcal{F}^+(X)$ consists of matching families of local sections of $\mathcal{F}$ modulo a~suitable equivalence relation: two families are identified if they agree on a~common refinement.

Nikolaus and Schweigert~\cite{nikolaus2} have generalized this concrete description to presheaves of bicategories. Intuitively, $\mathcal{F}^+$ assembles all compatible collections of locally defined objects together with the gluing data relating them.
Applying their construction here yields the definition
\begin{equation*}
 \stwoLineBdl(X) := \bigl(\cssAlgBdlbi\bigr)^{+}(X).
\end{equation*}
In other words, super 2-line bundles over $X$ are objects obtained by gluing central simple super algebra bundles along invertible bimodule bundles.
As a~sheafification, the plus construction comes with an embedding of the underlying presheaf; see \cite[Section~3.3]{Kristel2020}.

\begin{Proposition}
\label{embedding-of-algebra-bundles}
There are fully faithful functors
\begin{equation*}
\cssAlgBdl^{bi}(X)\to \stwoLineBdl(X)
\qquad
\csAlgBdl^{bi}(X) \to \twoLineBdl(X)
\end{equation*}
that embed central simple $($super$)$ algebra bundles into $($super{\rm) 2}-line bundles. Under the classification of Theorem~{\rm\ref{super-2-line-bundles}}, their images corresponds to the subsets
\begin{equation*}
\h^0(X,\mathrm{BW}_{\C})\times \h^1(X,\Z_2) \times \operatorname{Tor}\bigl(\h^3(X,\Z)\bigr) \quand \operatorname{Tor}\bigl(\h^3(X,\Z)\bigr),
\end{equation*}
respectively.
\end{Proposition}

The functors in Proposition~\ref{embedding-of-algebra-bundles} have the simple description to use the identity surjective submersion $\pi=\id_X$, the given algebra bundles, and identity bimodules and intertwiners.
We remark that the cohomological classification of the bicategories $\cssAlgBdl^{bi}(X)$ and $\csAlgBdl^{bi}(X)$ reproduces a~classical result of Donavan--Karoubi~\cite{DK70}.

\begin{Remark}
Theorem~\ref{super-2-line-bundles} can be enhanced to a~group isomorphism by equipping the bicategory $\stwoLineBdl(X)$ with a~symmetric monoidal structure. Then, by \cite[Corollary~4.4.10]{Kristel2020}, every super 2-line bundle is the tensor product of a~central simple algebra bundle (under the functor of Proposition~\ref{embedding-of-algebra-bundles}) and of a~super bundle gerbe (under the functor of Proposition~\ref{embedding-of-bundle-gerbes}).
\end{Remark}

Another consequence of the definition of (super) 2-line bundles as the sheafification is the following result \cite[Theorem~2.3.3]{Kristel2020}, which we state for completeness.

\begin{Theorem}\label{2-line-bundles-form-stack}
$($Super$)$ $2$-line bundles form a~$2$-stack on the site of smooth manifolds.
\end{Theorem}

\subsection{Framing}
\label{framing}

Many bicategories provide a~way to induce 1-morphisms from a~simpler class of morphisms that compose with strict associativity. We handle this situation within the following framework.

\begin{Definition} \label{Definition: framed bicategory}
 A \emph{framed bicategory} consists of a~bicategory $\mathscr{B}$, a~category $\mathscr{C}$ with the same objects,
 and a~functor $F \colon \mathscr{C} \to \mathscr{B}$ such that
 \begin{itemize}\itemsep=0pt
 \item $F$ is the identity on the level of objects and
 \item the image of every morphism of $\mathscr{C}$ admits a~right adjoint in $\mathscr{B}$.
 \end{itemize}
\end{Definition}

In this definition, by functor we mean a~\emph{strong pseudofunctor}, i.e., when $f$ and $g$ are composable morphisms in $\mathscr{C}$, it comes equipped with a~2-isomorphism $\eta_{f,}\colon F(g) \circ F(f) \Rightarrow F(g\circ f)$, and these are coherent on triples of composable morphisms. There is a~relation between framed bicategories and double categories~\cite{shulman1}; see \cite[Remark~2.1.7]{Kristel2022}.
Synonymously to Definition~\ref{Definition: framed bicategory}, we also use the terminology that ``$\mathscr{B}$ is framed under $\mathscr{C}$'', and that ``$F$ is a~framing of $\mathscr{B}$''.
\begin{Example}
The bicategory $\sAlgbi$ of super algebras is framed under the category $\sAlg$, by assigning to a~homomorphism $\varphi\colon A \to B$ between super algebras the $B$-$A$-bimodule $B_{\varphi}$, whose underlying vector space is $B$, the left $B$-action is multiplication, and the right $A$-action is multiplication along $\varphi$. The ``insidious'' nature of the bicategory $\sAlgBdlbi(X)$ prevents this framing to carry over to algebra bundles, see \cite[Section~3.3]{Kristel2022}. However, the bicategory~$\cssAlgBdlbi(X)$ of \emph{central simple} super algebra bundles is framed under the category of central simple super algebra bundles.
\end{Example}

\begin{Example}\label{refinements-of-bundle-gerbes}
The bicategory $\sGrb(X)$ of super bundle gerbes over $X$ is framed under a~category $\sGrbref(X)$.
Recall from Section~\ref{section-super-2-line-bundles} that a~super bundle gerbe is a~super 2-line whose super algebra bundle is the trivial one with fibre $\C$, in the following denoted by just $\mathscr{G}=(\pi,\mathcal{M},\mu)$. Objects in $\sGrbref(X)$ are super bundle gerbes $\mathscr{G}=(\pi,\mathcal{M},\mu)$ over $X$, and morphisms $\mathscr{G}_1 \to \mathscr{G}_2$ in $\sGrbref(X)$ are pairs $(\rho,u)$
consisting of a~smooth map $\rho\colon Y_1 \to Y_2$ and an invertible super line bundle morphism $u \colon \mathcal{M}_1 \to \rho^*\mathcal{M}_2$ over \smash{$Y_1^{[2]}$} such that $\pi_2 \circ \rho = \pi_1$
and the following diagram commutes:
\begin{equation}
\label{condition-for-refinements}
\begin{gathered}
 \xymatrix@C=4em{(\mathcal{M}_1)_{y',y''} \otimes (\mathcal{M}_1)_{y,y'}\ar[d]_{u_{y',y''} \otimes u_{y,y'}} \ar[r]^-{ \mu_1} & (\mathcal{M}_1)_{y,y''} \ar[d]^{u_{y,y''}}\\
 (\mathcal{M}_2)_{\rho(y'),\rho(y'')} \otimes (\mathcal{M}_2)_{\rho(y),\rho(y')} \ar[r]_-{\mu_2} & (\mathcal{M}_2)_{\rho(y),\rho(y'')}.}
\end{gathered}
\end{equation}
The morphisms in $\sGrbref(X)$ are called \emph{refinements}. We will not describe the framing functor $\sGrbref(X) \to \sGrb(X)$ here; instead,
we provide in Definition~\ref{refinements-of-super-2-line-bundles} below a~more general construction that can easily be reduced to the case of super bundle gerbes. \end{Example}

In the following we describe a~canonical framing for the bicategory of super 2-line bundles over $X$, deduced in \cite[Section~3.5]{Kristel2020} as a~feature of the plus construction.

\begin{Definition}
\label{refinements-of-super-2-line-bundles}
Let $\mathscr{L}_1=(\pi_1,\mathcal{A}_1,\mathcal{M}_1,\mu_1)$ and $\mathscr{L}_2=(\pi_2,\mathcal{A}_2,\mathcal{M}_2,\mu_2)$ be super 2-line bundles over $X$.
A \emph{refinement} $\mathscr{R}\colon \mathscr{L}_1\to \mathscr{L}_2$ is a~triple $\mathscr{R}=(\rho,\varphi,u)$ consisting of a~smooth map $\rho\colon Y_1 \to Y_2$ such that $\pi_2\circ \rho=\pi_1$, of a~homomorphism $\varphi\colon \mathcal{A}_1\to \rho^{*}\mathcal{A}_2$ of super algebra bundles over $Y_1$ and
of an invertible bundle morphism $u\colon \mathcal{M}_1 \to \rho^* \mathcal{M}_2$ over \smash{$Y_1^{[2]}$}
that over a~point \smash{$(y,y')\in Y_1^{[2]}$} restricts to an intertwiner
$
u_{y,y'}\colon (\mathcal{M}_1)_{y,y'} \to (\mathcal{M}_2)_{\rho(y),\rho(y')}
$
along the algebra homomorphisms $\varphi_{y'}\colon (\mathcal{A}_1)_{y'} \to (\mathcal{A}_2)_{\rho(y')}$ and $\varphi_{y}\colon (\mathcal{A}_1)_y \to (\mathcal{A}_2)_{\rho(y)}$, and renders the diagram
\begin{equation*}
\begin{aligned}
\xymatrix@C=4em{(\mathcal{M}_1)_{y',y''} \otimes_{(\mathcal{A}_1)_{y'}} (\mathcal{M}_1)_{y,y'}\ar[d]_{u_{y',y''} \otimes u_{y,y'}} \ar[r]^-{ \mu_1} & (\mathcal{M}_1)_{y,y''} \ar[d]^{u_{y,y''}}\\ (\mathcal{M}_2)_{\rho(y'),\rho(y'')} \otimes_{(\mathcal{A}_2)_{\rho(y')}} (\mathcal{M}_2)_{\rho(y),\rho(y')} \ar[r]_-{\mu_2} & (\mathcal{M}_2)_{\rho(y),\rho(y'')}}
\end{aligned}
\end{equation*}
commutative for all \smash{$(y,y',y'')\in Y_1^{[3]}$}.
\end{Definition}

Given two refinements
$
\mathscr{R}_{12}=(\rho_{12},\varphi_{12},u_{12})\colon \mathscr{L}_1 \to \mathscr{L}_2$,
$
\mathscr{R}_{23}=(\rho_{23},\varphi_{23},u_{23})\colon \mathscr{L}_2 \to \mathscr{L}_3$,
their composition is given by
\begin{equation*}
\mathscr{R}_{23} \circ \mathscr{R}_{12} := \bigl(\rho_{23} \circ \rho_{12},\rho_{12}^{*}\varphi_{23} \circ \varphi_{12},\bigl(\rho_{12}^{[2]}\bigr)^{*}u_{23} \circ u_{12}\bigr),
\end{equation*}
and the identity morphism of $\mathscr{L}$ is $(\id_Y,\id_{\mathcal{A}},\id_{\mathcal{M}})$.
This defines the category $\stwoLineBdlref_k(X)$.

The framing
\begin{equation}
\label{framing-of-2-line-bundles}
\stwoLineBdlref(X) \to \stwoLineBdl(X)
\end{equation}
is defined as follows. On the level of objects, it is the identity. On the level of morphisms, they associate to a~refinement $\mathscr{R}=(\rho,\varphi,u)\colon\mathscr{L}_1 \to \mathscr{L}_2$ the following 1-morphism $(\zeta,\mathcal{P},\phi)$. We define~${Z:= Y_1 \times_X Y_2}$ and $\zeta=\id_Z$. Consider the smooth map \smash{$\tilde \rho\colon Z \to Y_2^{[2]}$} with $\tilde\rho(y_1,y_2):=(\rho(y_1),y_2)$.
We define the bimodule bundle $\mathcal{P} := (\tilde\rho^{*}\mathcal{M}_2)_{\pr_1^{*}\varphi}$ over $Z$, in which the right $\pr_1^{*}\mathcal{A}_1$-action is along the algebra homomorphism $\pr_1^{*}\varphi\colon \pr_1^{*}\mathcal{A}_1\to \pr_1^{*}\rho^{*}\mathcal{A}_2$.
Over a~point $(y_1,y_2)\in Z$, its fibre is $\mathcal{P}_{y_1,y_2}=((\mathcal{M}_2)_{\rho(y_1),y_2})_{\varphi_{y_1}}$, which is an $(\mathcal{A}_2)_{y_2}$-$(\mathcal{A}_1)_{y_1}$-bimodule.
Finally, we define the intertwiner $\phi$ fibre-wise over a~point \smash{$((y_1,y_2),(y_1',y_2'))\in Z^{[2]}$} as the following composite:
\begin{equation*}
\xymatrix{
\hspace{-2em} \mathcal{P}_{y_1',y_2'} \otimes_{(\mathcal{A}_1)_{y_1'}} (\mathcal{M}_1)_{y_1,y_1'}=((\mathcal{M}_2)_{\rho(y_1'),y_2'})_{\varphi_{y_1'}} \otimes_{(\mathcal{A}_1)_{y_1'}} (\mathcal{M}_1)_{y_1,y_1'}
 \ar[d]^{\id \otimes u_{y_1,y_1'}}
\\
(\mathcal{M}_2)_{\rho(y_1^\prime), y_2^\prime} \otimes_{(\mathcal{A}_2)_{\rho(y_1^\prime)}} ((\mathcal{M}_2)_{\rho(y_1),\rho(y_1^\prime)})_{\varphi_{y_1}}
 \ar[d]^{(\mu_2)_{\rho(y_1),\rho(y_1'),y_2'}}
\\
((\mathcal{M}_2)_{\rho(y_1),y_2'})_{\varphi_{y_1}} \ar[d]^{(\mu_2)_{\rho(y_1),y_2,y_2'}^{-1}}
\\
(\mathcal{M}_2)_{y_2,y_2'} \otimes_{(\mathcal{A}_2)_{y_2}} ((\mathcal{M}_2)_{\rho(y_1,)y_2})_{\varphi_{y_1}}
= (\mathcal{M}_2)_{y_2,y_2'} \otimes_{(\mathcal{A}_2)_{y_2}} \mathcal{P}_{y_1,y_2}. \hspace{-20em} }
\end{equation*}
It is shown in \cite[Proposition~3.5.2]{Kristel2020} that this is indeed a~framing.

For later use, we state the following result.

\begin{Lemma}
\label{refinements-induce-isos}
A refinement $\mathscr{R}=(\rho,\varphi,u)$ induces an invertible $1$-morphism if the algebra bundle morphism $\varphi$ is invertible. \end{Lemma}

\begin{proof}
The induced 1-morphism is by Proposition~\ref{LemmaInvertibility}\,\ref{LemmaInvertibilityA} invertible if and only if its bimodule bundle $\mathcal{P} = (\tilde\rho^{*}\mathcal{M}_2)_{\pr_1^{*}\varphi}$ over $Z$ is invertible. In turn, a~bimodule bundle is invertible if and only if it is fibre-wise invertible \cite[Lemma~4.2.8]{Kristel2022}.
This is the case when $\mathcal{M}_2$ is invertible and~$\varphi$ is an isomorphism.
\end{proof}

\begin{Example}
\label{refinement-maps}
A standard source for refinements are situations in which a~super 2-line bundle $\mathscr{L}=(\pi,\mathcal{A},\mathcal{M},\mu)$ is given together with a~refinement of its surjective submersion ${\pi\colon Y \to X}$, i.e., given is another surjective submersion $\pi'\colon Y' \to X$ together with a~smooth map $\rho\colon Y' \to Y$ such that $\pi \circ \rho=\pi'$. Then, one may define a~second super 2-line bundle
$
\mathscr{L}^{\rho} := \bigl(\pi',\rho^{*}\mathcal{A},\bigl(\rho^{[2]}\bigr)^{*}\mathcal{M}, \allowbreak\bigl(\rho^{[3]}\bigr)^{*}\mu\bigr)
$
and observe that $\mathscr{R}^{\rho}:=(\rho,\id,\id)$ is a~refinement $\mathscr{L}^{\rho} \to \mathscr{L}$ that induces by Lemma~\ref{refinements-induce-isos} an isomorphism $\mathscr{L}^{\rho} \cong \mathscr{L}$ in $\stwoLineBdl(X)$.
\end{Example}

For later use, we record the appropriate notion of morphisms between framed bicategories.
\begin{Definition}
 Let $F \colon \mathscr{C}\to \mathscr{B}$ and $F' \colon \mathscr{C}' \to \mathscr{B}'$ be framed bicategories. A \emph{framed functor} consists of a~pair of functors $G_1\colon \mathscr{C} \to \mathscr{C}'$ and $G_2\colon \mathscr{B} \to \mathscr{B}'$ and a~pseudonatural equivalence
\[
 \alxydim{}{
 \mathscr{C}\ar[rr]^-{F} \ar[d]_-{G_1} & &\mathscr{B} \ar@{=>}[dll]|*+{\rho} \ar[d]^-{G_2} \\
 \mathscr{C}'\ar[rr]_-{F'} & & \mathscr{B}'.
 }
 \]
A framed functor is called an \emph{equivalence of framed bicategories} if $G_1$ and $G_2$ are equivalences.
\end{Definition}

\begin{Example}
There are natural and evident functors
\[\sGrbref(X)\to \stwoLineBdlref(X)\qquad \text{and} \qquad\cssAlgBdl(X) \to \stwoLineBdlref(X),
\]
 fitting into strictly commutative diagrams
\begin{equation*}
\alxydim{}{\sGrbref(X) \ar[r]\ar[d] & \sGrb(X) \ar[d] \\ \stwoLineBdlref(X) \ar[r] & \stwoLineBdl(X)}
\quand
\alxydim{}{ \cssAlgBdl(X) \ar[r] \ar[d]_{} & \cssAlgBdlbi(X) \ar[d] \\ \stwoLineBdlref(X) \ar[r] & \stwoLineBdl(X), }
\end{equation*}
which are framed functors, see \cite[Section~3.5]{Kristel2020}. We remark that all vertical functors in these diagrams are in fact fully faithful.
\end{Example}

\subsection{Complex conjugation}
\label{complex-conjugation}

Let $V$ be a~finite dimensional complex vector space. We recall that the \emph{complex conjugate} vector space $\overline{V}$ consists of the same additive group
$V$ and the scalar multiplication given by $\lambda \cdot v := \bar\lambda v$. Further, we recall that any linear map $f\colon V \to W$ between complex
vector spaces is a~linear map $\overline{V} \to \overline{W}$, for which we will use the notation $\overline{f}$ in order to emphasize the changed domain and codomain.

If $A$ is a~complex algebra, then we define the complex conjugate algebra $\overline{A}$ to have the complex conjugate of the underlying vector space, and the same multiplication (which is, as a~map $\overline{A} \times \overline{A} \to \overline{A}$, bilinear). Similarly, in the situation of bimodules, we complex conjugate the underlying vector space and leave the actions as they are. This way, complex conjugation becomes a~functor
\begin{equation*}
\overline{(..)}\colon\ \sAlgbi \to \sAlgbi,
\end{equation*}
which is strictly involutive.
It is clear that complex conjugation extends to algebra bundles and bimodule bundles, hence to the bicategory of algebra bundles, and moreover even to a~morphism of presheaves of bicategories.

It is now straightforward
to define the \emph{complex conjugate} of a~super 2-line bundle $\mathscr{L}=(\pi,\mathcal{A},\mathcal{M},\mu)$: indeed, we use the same surjective submersion
$\pi$, the complex conjugate algebra bundle $\overline{\mathcal{A}}$, the complex conjugate bimodule bundle $\overline{\mathcal{M}}$, and the intertwiner $\overline{\mu}$.
The complex conjugate super 2-line bundle defined in this way is denoted $\overline{\mathscr{L}}$. Similarly, we define complex conjugate 1-morphisms and 2-morphisms of super 2-line bundles.
It is clear that this construction extends to a~2-functor
\begin{equation}
\label{complex-conjugation-2-functor}
\overline{(..)} \colon\ \stwoLineBdl(X) \to \stwoLineBdl(X),
\end{equation}
which is strictly involutive, and again a~morphism of presheaves of bicategories.

\section{Super 2-line bundles over graded Lie groupoids}
\label{extension-to-graded-lie-groupoids}

The purpose of this section is to introduce a~very general theory
 of super 2-line bundles over \emph{graded Lie groupoids}, a~notion that includes various structures that have appeared before in the literature. We provide a~descent result and a~cohomological classification.

\subsection{Graded Lie groupoids}

 Lie groupoids encode all sorts of equivariance; they include orbifolds
and even more general structures. Traditional equivariance, where a~Lie group $G$ acts on a~smooth manifold $X$, is included by considering the action groupoid $\act XG$. The case of no equivariance is included by associating to a~smooth manifold $X$ the Lie groupoid $X_{\rm dis}$ with only identity morphisms.

If $\Gamma$ is a~Lie groupoid, we will denote by $s,t\colon \Gamma_1\to \Gamma_0$ the source and target maps, by $\Gamma_k := \Gamma_1 \ttimes st \Gamma_1 \ttimes st\cdots \ttimes st \Gamma_1$ the smooth manifold of $k$-tuples of composable morphisms, and by~${c\colon \Gamma_2 \to \Gamma_1}$ the composition map.

Every presheaf $\mathscr{F}$ on the category of smooth manifolds extends canonically to a~presheaf on the category of Lie groupoids (whose morphisms are smooth functors), via the homotopy limit over the associated simplicial manifold,
\begin{equation}
\label{homotopy-limit-equivariance}
\mathscr{F}(\Gamma)
:= \mathrm{holim} (
\alxydim{}{\mathscr{F}(\Gamma_0) \ar@<0.7ex>[r]\ar@<-0.7ex>[r] & \mathscr{F}(\Gamma_1) \ar@<1.4ex>[r]\ar[r]\ar@<-1.4ex>[r] & \mathscr{F}(\Gamma_2) \hdots}
)
,
\end{equation}
see, e.g.,~\cite{nikolaus2}. In particular, as super 2-line bundles form a~sheaf of bicategories on the site of smooth manifolds (see Theorem~\ref{2-line-bundles-form-stack}), we obtain a~bicategory $\stwoLineBdl(\Gamma)$ of super 2-line bundles over a~Lie groupoid $\Gamma$. If $G$ is a~Lie group acting on $X$, then, by definition, a~$G$-equivariant super 2-line bundle over $X$ is a~super 2-line bundle over the action groupoid $\Gamma=\act XG$.

In this paper, we need a~yet more general version of equivariance, namely, equivariance with respect to \emph{graded} Lie groupoids. These appeared in \cite[Definition~2.4]{Freed2011a} (in a~topological setting).

\begin{Definition}
\label{graded-groupoid}
 A \emph{graded Lie groupoid} is a~Lie groupoid $\Gamma$ equipped with a~smooth functor
 $\phi \colon \Gamma \to \mathrm{B}\mathbb{Z}_2$, where the latter denotes the groupoid with a~single object and two morphisms.
\end{Definition}

We use the convention that $\Z_2=\{1,-1\}$. Saying that $\phi$ is a~smooth functor is equivalent to saying that $\phi\colon \Gamma_1 \to \Z_2$ is a~locally constant map such that $\phi(\gamma_2\circ\gamma_1)=\phi(\gamma_2)\phi(\gamma_1)$ holds for all composable morphisms $\gamma_2$ and $\gamma_1$, and $\phi(\id_x)=1$ holds for all objects $x$. If $(\Gamma,\phi)$ is a~graded Lie groupoid, we denote by $(\Gamma,\phi)^{\rm even}$ the Lie groupoid with the same objects and only those morphisms $\gamma$ with $\phi(\gamma)=1$; this is a~sub-Lie groupoid of $\Gamma$.
The general idea is that the sub-groupoid $(\Gamma,\phi)^{\rm even}$ captures all sort of ``orbifold''-equivariance, for instance, with respect to group actions, while the additional morphisms with $\phi=-1$ parameterize possible ``orientifold''-equivariance (a.k.a.\ anti-equivariance, real-equivariance, etc.).

\begin{Example}\label{ex1}
We have the following main examples of graded Lie groupoids in mind:
\begin{enumerate}[(a)]\itemsep=0pt
\item
\label{ex-manifolds}
Smooth manifolds $X$, considered as groupoids with only identity morphisms, $\Gamma=X_{\rm dis}$. Note that necessarily $\phi=1$ in this case.

\item
\label{graded-Lie-group}
Graded Lie groups $(G,\varepsilon)$, i.e., Lie groups $G$ with a~smooth group homomorphism $\varepsilon\colon G \to \Z_2$. We consider $G$ as a~groupoid with a~single object, $\Gamma=\mathrm{B}G$, and equip it with the grading $\phi := \varepsilon$.

\item
\label{graded-action-groupoid}
Graded action groupoids, formed by a~smooth manifold $X$, a~graded Lie group $(G,\varepsilon)$, and a~smooth (right) action of $G$ on $X$. Then, $\Gamma := \act XG$ is the usual action groupoid, and the grading $\phi\colon X \times G \to \Z_2$ is $\phi := \varepsilon \circ \pr_G$. We also write $\act X{(G,\varepsilon)}$ for the graded action groupoid. It reduces to~\ref{ex-manifolds} when $G=\{e\}$ and to~\ref{graded-Lie-group} when $X=\{\ast\}$.

\item
\label{real-manifold}
Real manifolds $(X,\tau)$, i.e., smooth manifolds $X$ with a~smooth involution $\tau\colon X \to X$. We consider $\tau$ as an action of the group $\Z_2$ on $X$, and we consider $\Z_2$ graded under the identity map $\varepsilon:=\id_{\Z_2}$. Then, we form the graded action groupoid $\mathrm{Gr}(X,\tau) :=\act X{(\Z_2,\id)}$ of~\ref{graded-action-groupoid}.
\end{enumerate}
\end{Example}

Analogous to the extension of presheaves on manifolds to (ordinary, ungraded) Lie groupoids, one can extend them to \emph{graded} Lie groupoids whenever the presheaf is equipped with a~strictly involutive presheaf automorphism. In the following subsection, we illustrate this process for the presheaf of super 2-line bundles, whose automorphism is given by complex conjugation (see Section~\ref{complex-conjugation}).

\subsection{Main definitions}
\label{graded-equivariant-definitions}

We fix the convention that an upper index $(..)^{\phi}$ on super 2-line bundles, 1-morphisms, or 2-morphisms, where $\phi\colon X \to \Z_2$ is a~smooth map, means the application of the complex conjugation functor~\eqref{complex-conjugation-2-functor} on connected components with $\phi=-1$, and has no effect on components with $\phi=1$.

\begin{Definition} \label{equivariant structure}
 Let $(\Gamma,\phi)$ be a~graded Lie groupoid. A \emph{super $2$-line bundle} over $(\Gamma,\phi)$ is a~triple~$(\mathscr{L},\mathscr{P},\psi)$ consisting of
 \begin{itemize}\itemsep=0pt
 \item a~super 2-line bundle $\mathscr{L}$ over $\Gamma_0$,
 \item a~1-isomorphism of super 2-line bundles
 \begin{equation} \label{equivariant structure: 1-morphism}
 \mathscr{P}\colon \ s^*\mathscr{L} \to (t^*\mathscr{L})^\phi
 \end{equation}
 over $\Gamma_1$, and

 \item a~2-isomorphism of super 2-line bundles
 \begin{equation} \label{equivariant structure: 2-morphism}
 \psi\colon\ (\operatorname{pr}_1^*\mathscr{P})^{\operatorname{pr}_2^*\phi} \circ \operatorname{pr}_2^*\mathscr{P}
 \Rightarrow c^*\mathscr{P}
 \end{equation}
 over $\Gamma_2$,
 \end{itemize}
such that the following diagram over $\Gamma_3$ is commutative
\begin{equation} \label{equivariant structure: diagram}
\alxydim{@C=1em}{\pr_{12}^{*}\bigl(\operatorname{pr}_1^*\mathscr{P} ^{\operatorname{pr}_2^*\phi}\circ \operatorname{pr}_2^*\mathscr{P}\bigr)^{\pr_3^{*}\phi} \circ \pr_3^{*}\mathscr{P} \ar@{=>}[rr] \ar@{=>}[d]_{\pr_{12}^{*}\psi^{\pr_3^{*}\phi} \circ \id} && \pr_1^{*}\mathscr{P}^{\pr_2^{*}\phi\cdot \pr_3^{*}\phi} \circ \pr_{23}^{*}\bigl(\pr_1^{*}\mathscr{P}^{\pr_2^{*}\phi} \circ \pr_2^{*}\mathscr{P}\bigr) \ar@{=>}[d]^{\id \circ \pr_{23}^{*}\psi} \\
\pr_{12}^{*}c^{*}\mathscr{P}^{\pr_3^{*}\phi} \circ \pr_3^{*}\mathscr{P}\ar@{=>}[d] && \pr_1^{*}\mathscr{P}^{\pr_{23}^{*}c^{*}\phi} \circ \pr_{23}^{*}c^{*}\mathscr{P} \ar@{=>}[d] \\
(c \times \id)^{*}\bigl(\pr_1^{*}\mathscr{P}^{\pr_2^{*}\phi} \circ \pr_2^{*}\mathscr{P}\bigr) \ar@{=>}[d]_{(c\times \id)^{*}\psi} && (\id \times c)^{*}\bigl(\pr_1^{*}\mathscr{P}^{\pr_2^{*}\phi} \circ \pr_2^{*}\mathscr{P}\bigr) \ar@{=>}[d]^{(\id \times c)^{*}\psi}\\ (c \times \id)^{*}c^{*}\mathscr{P} \ar@{=>}[rr] && (\id \times c)^{*}c^{*}\mathscr{P}. }
\end{equation}
\end{Definition}

Here, and in what follows, unlabelled arrows represent the canonical 2-isomorphisms arising from either the bicategorical structure or the presheaf structure of super 2-line bundles.
Whenever appropriate, we will also call $(\mathscr{P},\psi)$ a~\emph{graded-equivariant structure} on the super 2-line bundle $\mathscr{L}$ over $\Gamma_0$. In fibre-wise notation, we have
\begin{equation*}
 \mathscr{P}_\gamma\colon\ \mathscr{L}_{s(\gamma)} \to \mathscr{L}_{t(\gamma)}^{\phi(\gamma)}
\quand
 \psi_{\gamma_1,\gamma_2} \colon \ \mathscr{P}_{\gamma_1}^{\phi(\gamma_2)} \circ \mathscr{P}_{\gamma_2} \Rightarrow \mathscr{P}_{\gamma_1 \circ \gamma_2}.
\end{equation*}

\begin{Remark}
\label{normalized graded-equivariant structures}
The pullback of $\mathscr{P}$ along $\id\colon \Gamma_0 \to \Gamma_1$ is a~1-isomorphism $\id^{*}\mathscr{P}\colon \mathscr{L} \to \mathscr{L}$, and the pullback of $\psi$ to $\Gamma_0$ gives a~2-isomorphism $(\id^{*}\mathscr{P})^2 \cong \id^{*}\mathscr{P}$. This shows that $\id^{*}\mathscr{P}\cong \id_{\mathscr{L}}$; moreover, under this identification the pullback of $\psi$ along either of the maps $\gamma \mapsto (\gamma,\id_{s(\gamma)})$ and $\gamma \mapsto (\id_{t(\gamma)},\gamma)$ is $\id_{\mathscr{P}}$.
\end{Remark}

\begin{Definition} \label{equivariant 1-morphisms}
 Let $(\mathscr{L}_1,\mathscr{P}_1,\psi_1)$ and $(\mathscr{L}_2,\mathscr{P}_2,\psi_2)$ be
 super 2-line bundles over $(\Gamma,\phi)$. A \emph{$1$-morphism} $(\mathscr{L}_1,\mathscr{P}_1,\psi_1)\to (\mathscr{L}_2,\mathscr{P}_2,\psi_2)$ is a~pair
 $(\mathscr{B},\eta)$ consisting of
 \begin{itemize}\itemsep=0pt
 \item a~1-morphism of super 2-line bundles
 \begin{equation} \label{equivariant 1-morphism: 1-morphism}
 \mathscr{B} \colon\ \mathscr{L}_1 \to \mathscr{L}_2
 \end{equation}
 over $\Gamma_0$, and
 \item a~2-isomorphism of super 2-line bundles
 \begin{equation} \label{equivariant 1-morphism: 2-morphism}
 \eta \colon\ \mathscr{P}_2 \circ s^*\mathscr{B} \Rightarrow (t^*\mathscr{B})^\phi \circ \mathscr{P}_1
 \end{equation}
 over $\Gamma_1$,
 \end{itemize}
such that the following diagram over $\Gamma_2$ is commutative
\begin{equation}
\label{equivariant 1-morphism: diagram}
\strut\hspace{-2.6em}\alxydim{@C=0.9em}{\operatorname{pr}_1^*\mathscr{P}_2^{\operatorname{pr}_2^*\phi} \circ \operatorname{pr}_2^*(\mathscr{P}_2 \circ s^{*}\mathscr{B}) \ar@{=>}[r] \ar@{=>}[d]_{\id \circ \pr_2^{*}\eta} & \bigl(\operatorname{pr}_1^*\mathscr{P}_2^{\operatorname{pr}_2^*\phi} \circ \operatorname{pr}_2^*\mathscr{P}_2 \bigr)\circ \pr_2^{*}s^{*}\mathscr{B} \ar@{=>}[rr]^-{\psi_2 \circ \id} && c^{*}\mathscr{P}_2 \circ \pr_2^{*}s^{*}\mathscr{B} \ar@{=>}[d] \\
\pr_1^{*}\mathscr{P}_2^{\pr_2^{*}\phi} \circ \pr_2^{*}\bigl(t^{*}\mathscr{B}^{\phi} \circ \mathscr{P}_1\bigr) \ar@{=>}[d] &&& c^{*}(\mathscr{P}_2 \circ s^{*}\mathscr{B}) \ar@{=>}[d]_{c^{*}\eta} \\ \pr_1^{*}(\mathscr{P}_2\circ s^{*}\mathscr{B})^{\pr_2^{*}\phi} \circ \pr_2^{*}\mathscr{P}_1 \ar@{=>}[d]_{\pr_1^{*}\eta^{\pr_2^{*}\phi} \circ \id} &&& c^{*}\bigl(t^{*}\mathscr{B}^{\phi} \circ \mathscr{P}_1\bigr) \ar@{=>}[d] \\
\pr_1^{*}\bigl(t^{*}\mathscr{B}^{\phi} \circ \mathscr{P}_1\bigr)^{\pr_2^{*}\phi} \circ \pr_2^*\mathscr{P}_1 \ar@{=>}[r] & \pr_1^{*}t^{*}\mathscr{B}^{\pr_1^{*}\phi \cdot \pr_2^{*}\phi} \circ \bigl(\pr_1^{*}\mathscr{P}_1^{\pr_2^{*}\phi} \circ \pr_2^{*}\mathscr{P}_1\bigr) \ar@{=>}[rr]_-{\id \circ \psi_1} && c^{*}t^{*}\mathscr{B}^{c^{*}\phi} \circ c^{*}\mathscr{P}_1. }\hspace{-4em}
\end{equation}
\end{Definition}

The identity 1-morphism of a~super 2-line bundle $(\mathscr{L},\mathscr{P},\psi)$ over $(\Gamma,\phi)$ is $(\id_{\mathscr{L}},\id_{\mathscr{P}})$.
In the situation of three super 2-line bundles $(\mathscr{L}_1,\mathscr{P}_1,\psi_1)$, $(\mathscr{L}_2,\mathscr{P}_2,\psi_2)$, and $(\mathscr{L}_3,\mathscr{P}_3,\psi_3)$ over $(\Gamma,\phi)$, the composition of two 1-morphisms
\begin{equation*}
\alxydim{@C=4em}{(\mathscr{L}_1,\mathscr{P}_1,\psi_1) \ar[r]^-{(\mathscr{B},\eta)} & (\mathscr{L}_2,\mathscr{P}_2,\psi_2) \ar[r]^-{(\mathscr{B}',\eta')} & (\mathscr{L}_3,\mathscr{P}_3,\psi_3)}
\end{equation*}
consists of the 1-morphism $\mathscr{B}' \circ \mathscr{B}$ and the 2-isomorphism
\begin{equation*}
\alxydim{@C=1.5em}{\mathscr{P}_3 \circ s^*(\mathscr{B}' \circ \mathscr{B}) \ar@{=>}[rr]^-{\eta' \circ \id} &&
 t^*\mathscr{B}'^\phi \circ \mathscr{P}_2 \circ s^*\mathcal{B} \ar@{=>}[rr]^-{\id \circ \eta } &&
 t^*(\mathscr{B}'\circ \mathscr{B})^\phi \circ \mathscr{P}_1.}
\end{equation*}

\begin{Definition}
\label{equivariant 2-morphisms}
Let $(\mathscr{L}_1,\mathscr{P}_1,\psi_1)$ and $(\mathscr{L}_2,\mathscr{P}_2,\psi_2)$ be
super 2-line bundles over $(\Gamma,\phi)$, and let~$(\mathscr{B},\eta)$ and $(\mathscr{B}',\eta')$
be 1-morphisms $(\mathscr{L}_1,\mathscr{P}_1,\psi_1)\to (\mathscr{L}_2,\mathscr{P}_2,\psi_2)$.
A 2-morphism
\begin{equation*}
\xi\colon\ (\mathscr{B},\eta)\Rightarrow (\mathscr{B}',\eta')
\end{equation*}
is a~2-morphism
$\xi\colon \mathscr{B} \Rightarrow \mathscr{B}'$
over $\Gamma_0$ such that the following diagram over $\Gamma_1$ is commutative
\begin{equation}
\label{compatibility-condition-for-2-morphisms}
\alxydim{}{\mathscr{P}_2 \circ s^*\mathscr{B} \ar@{=>}[d]_{\id \circ s^{*}\xi} \ar@{=>}[r]^{\eta} & t^*\mathscr{B}^\phi \circ \mathscr{P}_1 \ar@{=>}[d]^{t^{*}\xi^{\phi} \circ \id}
\\ \mathscr{P}_2 \circ s^*\mathscr{B}' \ar@{=>}[r]_{\eta'} & t^*\mathscr{B}'^\phi \circ \mathscr{P}_1.}
\end{equation}
\end{Definition}

Vertical and horizontal composition of 2-morphisms of super 2-line bundles over $(\Gamma,\phi)$ are given by the
vertical and horizontal composition in $\stwoLineBdl(\Gamma_0)$, respectively. Thus, super 2-line bundles over $(\Gamma,\phi)$ form a~bicategory. It is the central structure of this article.

\begin{Definition} \label{The bicategory of super 2-line bundles definition}
The bicategory of super 2-line bundles over a~graded Lie groupoid $(\Gamma,\phi)$ is denoted by $\stwoLineBdl(\Gamma,\phi)$.
\end{Definition}

\begin{Remark} \label{Inverse and adjoint of graded-equivariant 1-morphisms}
 Direct consequences of the definitions are the following statements about invertibility:
\begin{enumerate}[(a)]\itemsep=0pt
\item
 A 2-morphism $\xi\colon (\mathscr{B},\eta)\Rightarrow (\mathscr{B}',\eta')$ in $\stwoLineBdl(\Gamma,\phi)$ is invertible if and only if $\xi\colon \mathscr{B} \Rightarrow \mathscr{B}'$
is invertible in $\stwoLineBdl(\Gamma_0)$.

\item
A 1-morphism $(\mathscr{B},\eta)$ is invertible if and only if $\mathscr{B}$ is invertible in $\stwoLineBdl(\Gamma_0)$.

\item
A 1-morphism has $(\mathscr{B},\eta)$ has a~right (left)
 adjoint if and only if $\mathscr{B}$ has a~right (left) adjoint in $\stwoLineBdl(\Gamma_0)$.
\end{enumerate}
\end{Remark}

\begin{Remark}
\label{graded-equivariant-over-a-manifold}
If $\phi$ is trivial (i.e., $\phi=1$), we also just write
$\stwoLineBdl(\Gamma)$, and by inspection, we observe that this reproduces the homotopy limit construction~\eqref{homotopy-limit-equivariance}, i.e., the canonical definition of super 2-line bundles over Lie groupoids. In particular, when $\Gamma=X_{\rm dis}$ is a~smooth manifold, we have a~canonical equivalence of bicategories
$
\stwoLineBdl(X_{\rm dis})\cong \stwoLineBdl(X)$,
i.e., we reproduce the bicategory of super 2-line bundles over $X$ described in Section~\ref{section-super-2-line-bundles}.
\end{Remark}

\begin{Definition}\label{sub-bicategories}
We introduce several sub-bicategories
of $\stwoLineBdl(\Gamma,\phi)$:
\begin{enumerate}[(a)]\itemsep=0pt
\item
Restricting Definition~\ref{equivariant structure} along the fully faithful embedding $\sGrb\incl \stwoLineBdl$ of super bundle gerbes, we obtain the full sub-bicategory
\begin{equation*}
\sGrb(\Gamma,\phi) \subset \stwoLineBdl(\Gamma,\phi)
\end{equation*}
of super bundle gerbes over $(\Gamma,\phi)$.

\item
Restricting Definition~\ref{equivariant structure} along the fully faithful embedding $\cssAlgBdlbi\incl \stwoLineBdl$ of central simple super algebra bundles, we obtain the full sub-bicategory
\begin{equation*}
\cssAlgBdlbi(\Gamma,\phi) \subset \stwoLineBdl(\Gamma,\phi)
\end{equation*}
of central simple super algebra bundles over $(\Gamma,\phi)$.

\item
\label{linebdltriv}
Restricting Definition~\ref{equivariant structure} to the case where the super 2-line bundle $\mathscr{L}$ is the \emph{trivial} super 2-line bundle $\mathscr{I}$, we obtain the full sub-bicategory
\begin{equation*}
\stwoLineBdltriv(\Gamma,\phi)\subset \stwoLineBdl(\Gamma,\phi)
\end{equation*}
of ``graded-equivariant structures'' on the trivial super 2-line bundle.

\item
Restricting Definitions~\ref{equivariant structure},~\ref{equivariant 1-morphisms}, and~\ref{equivariant 2-morphisms}, as well as the sub-bicategories of (a), (b), and (c) to the non-super case, we obtain bicategories $\twoLineBdl(\Gamma,\phi)$, $\Grb(\Gamma,\phi)$, $\csAlgBdlbi(\Gamma,\phi)$, and $\twoLineBdltriv(\Gamma,\phi)$.
\end{enumerate}
\end{Definition}

All together, we have -- for each graded Lie groupoid $(\Gamma,\phi)$ -- a~commutative diagram
\begin{equation*}
\alxydim{}{& \twoLineBdltriv(\Gamma,\phi) \ar[dr]\ar[dl] \ar[d]\\\Grb(\Gamma,\phi) \ar[rd] \ar[d] & \stwoLineBdltriv(\Gamma,\phi) \ar[dr]\ar[dl] & \csAlgBdlbi(\Gamma,\phi) \ar[ld] \ar[d] \\ \sGrb(\Gamma,\phi) \ar[rd] & \twoLineBdl(\Gamma,\phi) \ar[d] & \cssAlgBdlbi(\Gamma,\phi) \ar[dl] \\ & \stwoLineBdl(\Gamma,\phi)}
\end{equation*}
of functors, in which all non-vertical functors are fully faithful.

\begin{Remark}
\label{trivial-things}
Similarly to Definition~\ref{sub-bicategories}\,\ref{linebdltriv} one can define bicategories $\sGrbtriv(\Gamma,\phi)$ and $\cssAlgBdlbi_{\mathrm{triv}}(\Gamma,\phi)$; however, as the embeddings of super bundle gerbes and central simple super algebra bundles into super 2-line bundles are fully faithful (see Propositions~\ref{embedding-of-bundle-gerbes} and~\ref{embedding-of-algebra-bundles}) these yield the same bicategory, i.e.,
\begin{equation*}
\sGrbtriv(\Gamma,\phi) =\cssAlgBdlbi_{\mathrm{triv}}(\Gamma,\phi)= \stwoLineBdltriv(\Gamma,\phi).
\end{equation*}
\end{Remark}

We also generalize the framing of super 2-line bundles recalled in Section~\ref{framing} from manifolds to graded Lie groupoids.

\begin{Definition} \label{definition: refinement of a~graded equivariant 2-line bundle}
 Let $(\Gamma,\phi)$ be a~graded Lie groupoid, and let $(\mathscr{L}_1,\mathscr{P}_1,\psi_1)$ and $(\mathscr{L}_2,\mathscr{P}_2,\psi_2)$
 be super 2-line bundles over $(\Gamma,\phi)$. A \emph{refinement} of super 2-line bundles over $(\Gamma,\phi)$
 between them consists of the following data: a~refinement $\mathscr{R}\colon \mathscr{L}_1 \to \mathscr{L}_2$
 of super 2-line bundles over $\Gamma_0$ together with a~2-isomorphism
 \[
 \eta \colon\ \mathscr{P}_2 \circ s^*F(\mathscr{R}) \Rightarrow t^*F(\mathscr{R})^\phi \circ \mathscr{P}_1
 \]
 over $\Gamma_1$, such that diagram~\eqref{equivariant 1-morphism: diagram} is commutative.
 Here, $F$ is the framing~\eqref{framing-of-2-line-bundles} of the bicategory of super 2-line bundles.
\end{Definition}

We denote the category of super 2-line bundles over $(\Gamma,\phi)$ and refinements by $\stwoLineBdlref(\Gamma,\phi)$.
\begin{Proposition} \label{proposition: graded equivariant 2-line bundles are framed}
 The bicategory $\stwoLineBdl(\Gamma,\phi)$ is framed under $\stwoLineBdlref(\Gamma,\phi)$.
\end{Proposition}
\begin{proof}
 The framing functor sends a~refinement $(\mathscr{R},\eta)$ to the 1-morphism $(F(\mathscr{R}),\eta)$. Combining~\cite[Proposition~3.5.2]{Kristel2020}
 and Remark~\ref{Inverse and adjoint of graded-equivariant 1-morphisms}, $(\mathscr{R},\eta)$ has a~right adjoint in $\stwoLineBdl(\Gamma,\phi)$.
\end{proof}

The framing restricts to all other seven bicategories in Definition~\ref{sub-bicategories} in a~straightforward way.

\begin{Remark}\label{Remark: the symmetric monoidal structure extends to the graded-equivariant case}
 The symmetric monoidal structure of super 2-line bundles described in \cite[Section~3.6]{Kristel2020} canonically extends from manifolds to graded Lie groupoids.
 This follows from an extension of the argument in~\cite{Kristel2020}, noting that pullbacks and complex conjugation are compatible with the symmetric monoidal structure of super
 2-line bundles over manifolds.
\end{Remark}

\subsection{Functoriality} \label{Subsection: Functoriality}

Graded Lie groupoids form a~category whose morphisms are \emph{even functors}, i.e., functors that strictly respect the grading.
Any even functor $F\colon (\Gamma,\phi) \to (\Gamma',\phi')$ clearly induces, via pullback, a~2-functor
\begin{equation*}
F^{*}\colon \ \stwoLineBdl(\Gamma',\phi') \to \stwoLineBdl(\Gamma,\phi).
\end{equation*}
 This way, super 2-line bundles form a~presheaf of bicategories over the category of graded Lie groupoids.
We prove in the next lemma that naturally isomorphic functors induce naturally equivalent pullback functors. In the ungraded case, this is
\cite[Proposition~6.2]{nikolaus2}.

\begin{Lemma}
\label{pullback-invariance-under-nat-iso}
Suppose $F,G\colon (\Gamma,\phi) \to (\Gamma',\phi')$ are even functors, and $\xi\colon F \Rightarrow G$ is an even smooth natural transformation, i.e., $\phi'(\xi_x)=1$ for each $x\in \Gamma_0$. Then, there is a~canonical natural equivalence $F^{*}\cong G^{*}$.
\end{Lemma}
\begin{proof}
Let $(\mathscr{L},\mathscr{P},\psi)$ be a~super 2-line bundle over $(\Gamma',\phi')$. We construct a~1-isomorphism $(\mathscr{B},\eta) \colon F^{*}(\mathscr{L},\mathscr{P},\psi) \to G^{*}(\mathscr{L},\mathscr{P},\psi)$. The 1-isomorphism $\mathscr{B}$ over $\Gamma_0$ is
\begin{equation*}
\alxydim{}{F_0^{*}\mathscr{L} =\xi^{*}s^{*}\mathscr{L} \ar[r]^-{\xi^{*}\mathscr{P}} & \xi^{*}(t^{*}\mathscr{L})^{\phi' \circ \xi} = G_0^{*}\mathscr{L},}
\end{equation*}
and the 2-morphism $\eta \colon G^{*}\mathscr{P} \circ s^*\mathscr{B} \Rightarrow (t^*\mathscr{B})^\phi \circ F^{*}\mathscr{P}$ is patched together from two applications of~$\psi$. This defines the component of a~natural equivalence $F^{*}\Rightarrow G^{*}$ at the object $(\mathscr{L},\mathscr{P},\psi)$. It is straightforward to provide the component at a~1-morphism, and to complete the proof.
\end{proof}

If an even functor $F\colon (\Gamma,\phi) \to (\Gamma',\phi')$ between graded Lie groupoids is a~strong equivalence (of the underlying Lie groupoids, see Appendix~\ref{bicategory-of-lie-groupoids}), then any inverse functor $G\colon \Gamma' \to \Gamma$ is automatically even, and as a~consequence of Lemma~\ref{pullback-invariance-under-nat-iso}, $F^{*}$ and $G^{*}$ establish an equivalence of bicategories
$
\stwoLineBdl(\Gamma',\phi') \cong \stwoLineBdl(\Gamma,\phi)$.

The above equivalence persists in fact when $F\colon \Gamma \to \Gamma'$ is only a~\emph{weak} equivalence (see Definition~\ref{Definition: weak equivalence of Lie groupoids}). More precisely, we have the following result.

\begin{Theorem}
\label{descent-for-graded-equivariant-2-line-bundles}
Suppose $F\colon (\Gamma,\phi) \to (\Gamma',\phi')$ is an even functor between graded Lie groupoids, such that the underlying functor $F\colon \Gamma \to \Gamma'$ is a~weak equivalence. Then,
\begin{equation*}
F^{*} \colon \ \stwoLineBdl(\Gamma',\phi') \to \stwoLineBdl(\Gamma,\phi).
\end{equation*}
is an equivalence of bicategories.
\end{Theorem}

\begin{proof}
For ungraded super 2-line bundles over Lie groupoids, this follows from Theorem~\ref{2-line-bundles-form-stack} and a~result of~\cite{nikolaus2} recalled in Appendix~\ref{appendix} as Proposition~\ref{Nikolaus-Schweigert weak equivalences}. Here, we generalize this to graded Lie groupoids.
By \cite[Proposition~5.7]{nikolaus2}, the weak equivalence $F\colon \Gamma \to \Gamma'$ factors through another Lie groupoid $\Omega$, via a~\emph{strong} equivalence $S\colon \Gamma \to \Omega$ (see Appendix~\ref{bicategory-of-lie-groupoids}) and a~\emph{covering functor} $T\colon \Omega \to \Gamma'$, i.e., $T$ is smoothly fully faithful in the sense of Definition~\ref{Definition: weak equivalence of Lie groupoids} and a~surjective submersion on the level of objects. In this situation, $\Omega$ becomes graded under $\omega: =\phi \circ S^{-1}$, in such a~way that $S$ and $T$ are even functors. Pulling back along $S$ and $S^{-1}$ induces an equivalence $\stwoLineBdl(\Gamma, \phi) \cong \stwoLineBdl(\Omega,\omega)$. It remains to show that $T^{*}$ induces an equivalence~${\stwoLineBdl(\Omega,\omega)\cong \stwoLineBdl(\Gamma',\phi')}$.

To establish an inverse to $T^{*}$, we intend to use descent along $T$, and for this, we need to work on the fibre product Lie groupoids (recalled in Definition~\ref{definition: fibre product of Lie groupoids})
\begin{equation*}
\Omega^{[k]} := \Omega \times_{\Gamma'} \Omega \times_{\Gamma'} \times\dots\times_{\Gamma'} \Omega.
\end{equation*}
We require two
facts about covering functors. The first is that covering functors are not only surjective submersions on the level of objects, but also -- automatically -- surjective submersions on the level of morphisms. Therefore, the fibre products $\Omega^{[k]}$ can be formed strictly, i.e., separately on the level of objects and on the level of morphisms. Second, the diagonal functors~${
\Delta^{k}\colon \Omega \to \Omega^{[k]}}
$
are strong equivalences \cite[Lemma~8.1\,(ii)]{nikolaus2}. Let $\pr_1,\pr_2\colon \Omega^{[2]} \to \Omega$ be the projection functors. We have $\pr_1\circ\Delta^2 = \id_{\Omega} = \pr_2\circ \Delta^2$. By \cite[Lemma~5.3]{nikolaus2}, it follows that $\pr_1$ and $\pr_2$ are strong equivalences, too, and both inverses of $\Delta^2$; in particular, there is a~natural isomorphism~${\pr_1\cong \pr_2}$.

We note that the fibre products $\Omega^{[k]}$ are canonically graded. Indeed, if $(\alpha_1,\alpha_2)\in \Omega_1 \times_{\Gamma'_1}\Omega_1$ is a~morphism in \smash{$\Omega^{[2]}$}, then we have, since $T$ is even,
\begin{equation*}
\omega(\alpha_1) = \phi'(T(\alpha_1))=\phi'(T(\alpha_2))=\omega(\alpha_2),
\end{equation*}
so that the grading is unambiguously defined through either of the morphisms. Moreover, all projection functors $\Omega^{[k]} \to \Omega^{[l]}$ and the diagonal functors $\Delta^{k}$ are even, and the natural transformation $\pr_1\cong \pr_2$ is even, too.

Now let $(\mathscr{L},\mathscr{P},\psi)$ be a~super 2-line bundle over $(\Omega,\omega)$. By Lemma~\ref{pullback-invariance-under-nat-iso}, there is a~canonical 1-isomorphism
$
 (\mathscr{B},\eta) \colon \pr_1^{*}(\mathscr{L},\mathscr{P},\psi) \to \pr_2^{*}(\mathscr{L},\mathscr{P},\psi)
$
over $\Omega^{[2]}$. Similarly, we obtain a~canonical 2-isomorphism
\begin{equation*}
\xi\colon \pr_{23}^{*}(\mathscr{B},\eta) \circ \pr_{12}^{*}(\mathscr{B},\eta) \Rightarrow \pr_{13}^{*}(\mathscr{B},\eta)
\end{equation*}
over $\Omega^{[3]}$, satisfying a~coherence condition over 4-fold fibre products. This establishes a~canonical descent structure on $(\mathscr{L},\mathscr{P},\psi)$, and it remains to perform the descent.

Separating structure, $\mathscr{L}$ is a~super 2-line bundle over $\Omega_0$, $\mathscr{B} \colon \pr_1^{*}\mathscr{L} \to \pr_2^{*}\mathscr{L}$ is a~1-isomorphism on $\Omega_0 \times_{\Gamma'_0} \Omega_0$, and $\xi$ is a~2-isomorphism on $\Omega_0 \times_{\Gamma'_0} \Omega_0 \times_{\Gamma'_0}\Omega_0$; hence, a~descent object for the 2-stack of super 2-line bundles. Thus, by Theorem~\ref{2-line-bundles-form-stack}, there exists a~super 2-line bundle~$\mathscr{L}'$ on $\Gamma_0'$ together with a~1-isomorphism $T_0^{*}\mathscr{L}'\cong \mathscr{L}$.
Next, we consider the 1-isomorphism~${\mathscr{P}\colon s^{*}\mathscr{L} \to (t^{*}\mathscr{L})^{\omega}}$ over $\Omega_1$, and the 2-isomorphism
\begin{equation*}
\alxydim{@C=4em}{\pr_1^{*}s^{*}\mathscr{L} \ar[r]^-{\pr_1^{*}\mathscr{P}} \ar[d]_{s^{*}\mathscr{B}} & \pr_1^{*}(t^{*}\mathscr{L})^{\omega}\ar[d]^{t^{*}\mathscr{B}^{\omega}} \\ \pr_2^{*}s^{*}\mathscr{L} \ar@{=>}[ur]|\eta \ar[r]_-{\pr_2^{*}\mathscr{P}} & \pr_2^{*}(t^{*}\mathscr{L})^{\omega}. }
\end{equation*}
Together with the compatibility condition~\eqref{compatibility-condition-for-2-morphisms} between $\eta$ and $\xi$, this means precisely that $\mathscr{P}$ descends along $T_1$ to a~1-isomorphism $\mathscr{P'}\colon s^{*}\mathscr{L}' \to (t^{*}\mathscr{L}')^{\phi'}$ over $\Gamma'_1$.
Finally, the compatibility condition~\eqref{equivariant 1-morphism: diagram} means that $\psi$ descends to a~2-isomorphism
\begin{equation*}
\psi' \colon\ (\operatorname{pr}_1^*\mathscr{P}')^{\operatorname{pr}_2^*\phi'} \circ \operatorname{pr}_2^*\mathscr{P}'
 \Rightarrow c^*\mathscr{P}',
\end{equation*}
which then automatically satisfies the associativity condition~\eqref{equivariant structure: diagram}. This shows that $(\mathscr{L}',\mathscr{P}',\psi')$ is a~super 2-line bundle over $(\Gamma',\phi')$. By construction, it is an essential preimage of $(\mathscr{L},\mathscr{P},\psi)$ under the pullback functor $T^{*}$.

This proves that $T^{*}$ is essentially surjective, which is the most important part of our theorem. Similarly, one proves that it is also fully faithful.
\end{proof}

\begin{Remark}
Theorem~\ref{descent-for-graded-equivariant-2-line-bundles} restricts to all sub-\emph{sheaves} of $\stwoLineBdl$ that we have described in Definition~\ref{sub-bicategories}, namely, to $\Grb$, $\sGrb$, and $\twoLineBdl$. It does not restrict to $\cssAlgBdlbi$ and versions thereof.
\end{Remark}

\subsection{Cohomological classification}

If $\Gamma$ is a~Lie groupoid and $\mathcal{F}$ is a~sheaf of abelian groups on the category of smooth manifolds, one defines the \v Cech cohomology groups $\check \h^k(\Gamma,\mathcal{F})$ using a~semi-simplicial open cover $\mathcal{U} = \bigl(\mathcal{U}^k\bigr)_{k\in \N}$ of the simplicial smooth manifold $N(\Gamma)$, i.e., \smash{$\mathcal{U}^k=\bigl(U^k_{i}\bigr)_{i\in I^k}$} is an open cover of $\Gamma_k$, and for each face map $f\colon \Gamma_k \to \Gamma_l$ there is a~map $\tilde f\colon I^k \to I^l$ such that \smash{$f\bigl(U^k_i\bigr)\subset U^{l}_{\tilde f(i)}$}.
Then, one can consider the double complex
\begin{equation*}
\alxydim{}{
\vdots& \vdots&
\\
\check C^1\bigl(\mathcal{U}^0,\mathcal{F}\bigr) \ar[u] \ar[r]_{\Delta} & \check C^1(\mathcal{U}^1,\mathcal{F}) \ar[u] \ar[r] & \hdots
\\
\check C^0\bigl(\mathcal{U}^0,\mathcal{F}\bigr) \ar[u]^{\delta} \ar[r]_{\Delta} & \check C^0(\mathcal{U}^1,\mathcal{F}) \ar[r] \ar[u]^{\delta}
& \hdots,}
\end{equation*}
whose columns are the usual \v Cech complexes of the open covers $\mathcal{U}^0$, $\mathcal{U}^1$, etc., and the differential~${\Delta\colon \check C^n\bigl(\mathcal{U}^l,\mathcal{F}\bigr) \to \check C^n\bigl(\mathcal{U}^{l+1},\mathcal{F}\bigr)}$ in the rows is defined by
\begin{equation*}
\Delta := \sum_{p=0}^{l+1} (-1)^{p} \partial^{*}_p,
\end{equation*}
where $\partial_p\colon \Gamma_{l+1} \to \Gamma_{l}$ is the $p$-th face map. Explicitly, $\partial_p\colon \Gamma_{1} \to \Gamma_0$ are $\partial_0=s$ and $\partial_1=t$, and~${\partial_p\colon \Gamma_{l+1} \to \Gamma_{l}}$ is given by
\begin{equation*}
\partial_p(\gamma_{0},\dots,\gamma_l) := \begin{cases}
(\gamma_{1},\dots,\gamma_{l}), & p=0 ,\\
(\gamma_0,\dots,\gamma_{p-1} \circ \gamma_{p},\dots,\gamma_l), & 0<p<l+1, \\
(\gamma_{0},\dots,\gamma_{l-1}), & p=l+1.
\end{cases}
\end{equation*}
For instance, \smash{$\Delta|_{ \check C^n(\mathcal{U}^0,\mathcal{F})}=s^{*}-t^{*}$}.
As usual, we form the total cohomology of this double complex with differential $D|_{\text{\v C}^p(\mathcal{U}^q,\mathcal{F})} := \delta + (-1)^{p}\Delta$.
There is a~partial order on simplicial open covers, and the direct limit over the total cohomology of the above double complex is called the \emph{\v Cech cohomology of the Lie groupoid $\Gamma$ with values in $\mathcal{F}$} and is denoted by $\check \h^n(\Gamma,\mathcal{F})$. We refer to~\cite{tu1} for more details and background on Lie groupoid cohomology. We note that the projection to the first column in the double complex induces a~map $p_0\colon \check \h^k(\Gamma,\mathcal{F}) \to \check\h^k(\Gamma_0,\mathcal{F})$ to the ordinary \v Cech cohomology of the manifold $\Gamma_0$. We also note the following special cases:
\begin{enumerate}[(i)]\itemsep=0pt
\item
If $X$ is a~smooth manifold, then $p_0\colon \check \h^k(X_{\rm dis},\mathcal{F}) \to \check \h^k(X,\mathcal{F})$ is an isomorphism \cite[Proposition~4.8]{tu1}.

\item
If $G$ is a~Lie group, and $\mathrm{B}G=\act \ast G$ is the corresponding Lie groupoid, $\check \h^k(\mathrm{B}G,\mathcal{F})$ is the (``Segal--Mitchison'', ``Brylinski'') smooth group cohomology of the Lie group $G$;
this follows, e.g., by comparing the constructions above with those in~\cite{brylinski3}.
\end{enumerate}
As usual, if $A$ is an abelian Lie group, we will denote by the same letter the corresponding sheaf of \emph{locally constant} $A$-valued functions, and by $\underline{A}$ the corresponding sheaf of \emph{smooth} $A$-valued functions.

\begin{Example}
\label{cech-cohomology-in-low-degrees}
Let us examine the \v Cech cohomology $\check \h^n(\Gamma,\Z_2)$ in lowest degrees:
\begin{enumerate}[(a)]\itemsep=0pt
\item
\label{cech-cohomology-in-low-degrees:a}
A class in $\check \h^0(\Gamma,\Z_2)$ is represented by a~family $\alpha=(\alpha_i)_{i\in I^0}$ of smooth maps $\alpha_i\colon U_i^0 \to \Z_2$. The cocycle condition $\delta \alpha=0$ requires $\alpha_i=\alpha_j$ on any intersection \smash{$U_i^0 \cap U_j^0$} of open sets of $\mathcal{U}^0$, and the cocycle condition $\Delta \alpha=0$ requires, for each open set $U_i^1$ of $\mathcal{U}^1$ that~\smash{$s^{*}\alpha_{s(i)}=t^{*}\alpha_{t(i)}$}. In other words, $\alpha$ is a~sign for each connected component of $\Gamma_0$, with signs coinciding whenever there is a~morphism between two components.

\item
\label{cech-cohomology-in-low-degrees:b}
A class in $\check\h^1(\Gamma,\Z_2)$ is represented by a~pair $(\epsilon,\sigma)$ of a~family $\epsilon=(\epsilon_{ij})_{i,j\in I^0}$, with $\epsilon_{ij}\colon U_i^0 \cap U_j^0 \to \Z_2$, and a~family $\sigma=(\sigma_i)_{i\in I^1}$ with $\sigma_i\colon U_i^1 \to \Z_2$, satisfying three conditions:
\begin{itemize}\itemsep=0pt
\item
$\delta \epsilon=0$: the usual \v Cech cocycle condition

\item
$\delta \sigma-\Delta \epsilon=0$, i.e., $s^{*}\epsilon_{s(i),s(j)}-t^{*}\epsilon_{t(i),t(j)}= \sigma_j-\sigma_i$ on $U^1_i \cap U_j^1$,

\item
$\Delta \sigma=0$, i.e., $\pr_2^{*}\sigma_{\pr_2(i)}-c^{*}\sigma_{c(i)}+\pr_1^{*}\sigma_{\pr_1(i)}=0$.
\end{itemize}
\end{enumerate}
\end{Example}

We note the following result, which is merely a~corollary of its graded version (to be proved below as Theorem~\ref{classification-of-graded-equivariant-two-line-bundles}), obtained by considering the case $\phi=1$ there.

\begin{Proposition}
\label{classification-of-equivariant-two-line-bundles}
Super $2$-line bundles over a~Lie groupoid $\Gamma$ are classified by the \v Cech cohomology of
$\Gamma$, in such a~way that there is a~canonical bijection
\begin{equation*}
\hc 0 (\stwoLineBdl(\Gamma)) \cong \check\h^0(\Gamma,\Z_2) \times \check\h^1(\Gamma,\Z_2) \times \check\h^2(\Gamma,\sheaf{\C^{\times}}), \end{equation*}
where $\hc 0$ denotes the set of isomorphism classes of objects.
Under this bijection, the subsets $\hc 0 \sGrb(\Gamma)$ and $\hc 0 (\cssAlgBdl(\Gamma))$ correspond to the subsets $\check \h^1(\Gamma,\Z_2) \times \check\h^2(\Gamma,\sheaf{\C^{\times}})$ and $\check \h^0(\Gamma,\Z_2) \times \check\h^1(\Gamma,\Z_2) \times p_0^{-1}\bigl(\operatorname{Tor} \check \h^2(\Gamma_0,\sheaf{\C^{\times}})\bigr)$, respectively.
 \end{Proposition}

The reduction of this result to $\Gamma=X_{\rm dis}$ reproduces Theorem~\ref{super-2-line-bundles}.

Our aim is now to extend \v Cech cohomology for Lie groupoids to graded Lie groupoids. For this purpose, we require that the sheaf $\mathcal{F}$ is equipped with a~$\Z_2$-action by morphisms of sheaves of abelian groups. We denote this action by $\varepsilon \lact f$, for $\varepsilon\in \{\pm 1\}$ and some $f\in \mathcal{F}(U)$.
Now we modify the setup described above in the following way. The row's differential $\Delta$ is replaced by a~differential $\Delta^{\phi}\colon \check C^n\bigl(\mathcal{U}^l,\mathcal{F}\bigr) \to \check C^n\bigl(\mathcal{U}^{l+1},\mathcal{F}\bigr)$ defined by
\begin{equation*}
\Delta^{\phi} :=\sum_{p=0}^{l+1} \varepsilon^{\phi}_{l+1,p}\lact (-1)^{p} \partial^{*}_p,
\end{equation*}
where the function $\varepsilon^{\phi}_{l,p}\colon \Gamma_{l} \to \Z_2=\{\pm 1\}$ is defined for $l>0$ and $0\leq p\leq l$ by
\begin{equation*}
\varepsilon^{\phi}_{l,p}(\gamma_1,\dots,\gamma_l) := \begin{cases}
1, & p< l ,\\
\phi(\gamma_{l}), &p=l.
\end{cases}
\end{equation*}
One can check that this satisfies $\Delta^{\phi}\circ \Delta^{\phi}=0$, and that $\Delta^{\phi}$ commutes with (the unchanged) \v Cech coboundary operator $\delta$.
Thus, we obtain, via the same procedure as in the ungraded case, the \emph{\v Cech cohomology $\check \h^n((\Gamma,\phi),\mathcal{F})$ of the graded Lie groupoid $(\Gamma,\phi)$}.

\begin{Remark}
\label{sfsdfdf}\quad
\begin{itemize}\itemsep=0pt
\item
If the $\mathcal{\Z}_2$-action on $\mathcal{F}$ is the trivial one, then the \v Cech cohomology of the graded Lie groupoid $(\Gamma,\phi)$ coincides with the one of the ungraded Lie groupoid $\Gamma$.

\item
Since the modification for the grading $\phi$ does not affect the \v Cech differential $\delta$, the projection $p_0\colon \check \h^n((\Gamma,\phi),\mathcal{F}) \to \check \h^n(\Gamma_0,\mathcal{F})$ is still well-defined.

\item
If $G$ is a~graded Lie group, i.e., a~Lie group equipped with a~homomorphism ${\varepsilon\colon G \to \Z_2}$, then $\check \h^n((\mathrm{B}G,\varepsilon),\mathcal{F})$ defines a~generalization of the (``Segal--Mitchison'', ``Brylinski'') \linebreak smooth group cohomology to graded Lie groups.
\end{itemize}
\end{Remark}

\begin{Example}
\label{reconstruction-of-equivariant-bundle-gerbes}
Let us examine the relevant \v Cech cohomology group $\check \h^2((\Gamma,\phi),\sheaf{\C^\times})$, where the $\Z_2$-action on $\sheaf{\C^{\times}}$ is given by complex conjugation.
A class in $\check \h^2((\Gamma,\phi),\sheaf{\C^\times})$ is represented by a~triple $(\mu,\eta,f)$ consisting of a~family $\mu=(\mu_{ijk})_{i,j,k\in I^0}$ of smooth maps $\mu_{ijk}\colon U^0_i\cap U^0_j \cap U^0_k \to \C^\times$, of a~family $\eta=(\eta_{ij})_{i,j\in I^1}$ of smooth maps $\eta_{ij}\colon U^1_i \cap U^1_j \to \C^\times$, and of a~family $f=(f_{i})_{i\in I^2}$ of smooth maps $f_i\colon U^2_i \to \C^\times$, satisfying the following four conditions:
\begin{itemize}\itemsep=0pt

\item
$\delta \mu=0$, the usual \v Cech cocycle condition.

\item
$\delta \eta+\Delta^{\phi} \mu=0$: for all $\gamma\in U^1_i\cap U^1_j \cap U^1_k$,
\begin{equation}
\label{graded-equivariant-cohomology:2}
\mu_{s(i),s(j),s(k)}(s(\gamma))\cdot \eta_{jk}(\gamma)\cdot\eta_{ik}(\gamma)^{-1}\cdot\eta_{ij} (\gamma)=\phi(\gamma)\lact \mu_{t(i),t(j),t(k)}(t(\gamma)).
\end{equation}

\item
$\delta f-\Delta^{\phi} \eta=0$: for all $(\gamma_1,\gamma_2)\in U^2_i\cap U^2_j$,
\begin{gather}
f_j(\gamma_1,\gamma_2)\cdot f_i(\gamma_1,\gamma_2)^{-1}\cdot \eta_{c(i),c(j)}(\gamma_1\circ \gamma_2)\nonumber\\
\qquad=\phi(\gamma_2)\lact \eta_{\pr_1(i),\pr_1(j)}(\gamma_1)\cdot\eta_{\pr_2(i),\pr_2(j)}(\gamma_2).\label{graded-equivariant-cohomology:3}
\end{gather}

\item
$\Delta^{\phi} f=0$: for $(\gamma_1,\gamma_2,\gamma_3)\in U_i^3$,
\begin{gather}
\label{graded-equivariant-cohomology:4}
f_{\partial_0(i)}(\gamma_2,\gamma_3) \cdot f_{\partial_2(i)}(\gamma_1,\gamma_2\circ \gamma_3)=f_{\partial_1(i)}(\gamma_1\circ\gamma_2,\gamma_3)\cdot \phi(\gamma_3)\lact f_{\partial_3(i)}(\gamma_1,\gamma_2).
\end{gather}
\end{itemize}
Under the usual procedure of classifying bundle gerbes by \v Cech cohomology (see, e.g., \cite{gawedzki8,stevenson1}), one can construct from the \v Cech cocycle $\mu$ a~bundle gerbe $\mathscr{G}$ over $\Gamma_0$. Moreover, one can construct from $\eta$ -- using~\eqref{graded-equivariant-cohomology:2} and that the complex conjugate cocycle constructs the complex conjugate bundle gerbe -- a~1-isomorphism $\mathscr{P}\colon s^{*}\mathscr{G} \to t^{*}\mathscr{G}^{\phi}$ over $\Gamma_1$. Finally, one can construct from~$f$~-- using~\eqref{graded-equivariant-cohomology:3} -- a~2-isomorphism $\psi\colon (\operatorname{pr}_1^*\mathscr{P})^{\operatorname{pr}_2^*\phi} \circ \operatorname{pr}_2^*\mathscr{P}
 \Rightarrow c^*\mathscr{P}$ over $\Gamma_2$, and~\eqref{graded-equivariant-cohomology:4} implies an equality
$
(c\times \id)^{*}\psi \bullet (\pr_{12}^{*}\psi^{\pr_3^{*}\phi} \circ \id) =(\id \times c)^{*}\psi \bullet (\id \circ \pr_{23}^{*}\psi)
$
of 2-isomorphisms over $\Gamma_3$, which is nothing but the commutative diagram~\eqref{equivariant structure: diagram}.
Summarizing, we have described the reconstruction of a~bundle gerbe over $(\Gamma,\phi)$ (see Definition~\ref{sub-bicategories}) from a~cocycle. It is standard to show that this establishes a~bijection
$
\check \h^2((\Gamma,\phi),\sheaf{\C^\times}) \cong \hc 0 \Grb(\Gamma,\phi)$;
in other words, $\check \h^2((\Gamma,\phi),\sheaf{\C^\times})$ classifies bundle gerbes over $(\Gamma,\phi)$.
\end{Example}

\begin{Theorem}
\label{classification-of-graded-equivariant-two-line-bundles}
Super $2$-line bundles over a~graded Lie groupoid $(\Gamma,\phi)$ are classified by the \v Cech cohomology of
$(\Gamma,\phi)$, namely, there is a~canonical bijection
\begin{equation*}
\hc 0 (\stwoLineBdl(\Gamma,\phi)) \cong \check\h^0((\Gamma,\phi),\Z_2) \times \check\h^1((\Gamma,\phi),\Z_2) \times \check\h^2((\Gamma,\phi),\sheaf{\C^\times}),
\end{equation*}
where the sheaves $\Z_2$ are equipped with the trivial involution, and the sheaf $\underline{\C^{\times}}$ is equipped with the involution given by complex conjugation. The subsets $\hc 0 (\sGrb(\Gamma,\phi))$ and $\hc 0 (\cssAlgBdl(\Gamma,\phi))$ correspond, respectively, to the subsets
\begin{align*}
&\check\h^1((\Gamma,\phi),\Z_2) \times \check\h^2((\Gamma,\phi),\sheaf{\C^{\times}})
\quand\\
&\check \h^0((\Gamma,\phi),\Z_2) \times \check\h^1((\Gamma,\phi),\Z_2)\times p_0^{-1}\bigl(\operatorname{Tor} \check \h^2(\Gamma_0,\sheaf{\C^{\times}})\bigr).
\end{align*}
\end{Theorem}

Here we recall, see Remark~\ref{sfsdfdf}, that $p_0^{-1}\bigl(\operatorname{Tor} \check \h^2(\Gamma_0,\sheaf{\C^{\times}})\bigr)$ denotes classes in $\check \h^2((\Gamma,\phi),\sheaf{\C^{\times}}))$ whose image under $p_0$ is torsion.

\begin{proof}
The main point is that extracting local data from a~super 2-line bundle $(\mathscr{L},\mathscr{P},\psi)$ over~$(\Gamma,\phi)$, as described for smooth manifolds in Remark~\ref{extracting-cocycle}, yields the correct cocycle data, in particular, the correct signs for the grading $\phi$. Indeed, we collect the data $(\alpha_i,\varepsilon_{ij},\mu_{ijk})$ extracted there for an open cover $\mathcal{U}^0$ of $\Gamma_0$. Next, we work with respect to an open cover~$\mathcal{U}^1$ on $\Gamma_1$. We first consider the restriction of $\mathscr{P}$ to an open set $U^1_i$ and consider the composite 1-isomorphism \smash{$\mathscr{P}_{i}: =t^{*}\mathscr{N}_{t(i)}^{\phi}\circ \mathscr{P}|_{U_i^1} \circ s^{*}\mathscr{N}_{s(i)}\colon s^{*}\mathcal{A}_{s(i)} \to t^{*}\mathcal{A}_{t(i)}^{\phi}$}. Since the inclusion of central simple algebra bundles into super 2-vector bundles is fully faithful (see Proposition~\ref{embedding-of-algebra-bundles}), $\mathscr{P}_i$ can be identified with an invertible \smash{$\bigl(t^{*}\mathcal{A}^{\phi}_{t(i)}\bigr)$}-\smash{$(s^{*}\mathcal{A}_{s(i)})$}-bimodule bundle $\mathcal{P}_i$. In particular, this implies \smash{$t^{*}\alpha_{t(i)}=s^{*}\alpha_{s(i)}$}, providing the remaining cocycle condition in Example~\ref{cech-cohomology-in-low-degrees}\,\ref{cech-cohomology-in-low-degrees:a}; thus, $\alpha$~gives the class $\check \h^0(\Gamma,\Z_2)=\check \h^0((\Gamma,\phi),\Z_2)$. Moreover, since, by construction, the algebras $A_i$ are fixed representatives of the Morita equivalence classes (see Remark~\ref{extracting-cocycle}), we have \smash{$A_{s(i)}=A_{t(i)}^{\phi_i}$}, where~$\phi_i$ is the value of $\phi$ over $U_i^1$.
We can assume that the bimodule bundle $\mathcal{P}_i$ is trivializable, $\mathcal{P}_i \cong U_i^1 \times P_{i}$, where $P_i$ is an invertible $\tilde A_i$-$\tilde A_i$-bimodule, where \smash{$\tilde A_i:=A_{t(i)}^{\phi_i}=A_{s(i)}$} for any $i\in I^1$. Thus, $P_i\cong-\sigma_{i}\tilde A_i$, for a~locally constant map $\sigma_i\colon U^1_i \to \Z_2$, where $-A$ stands for $\Pi A$, the grading reversal, see Remark~\ref{extracting-cocycle}.
On a~double overlap $U^1_{ij}:=U_i^1 \cap U^1_j$, we observe that the intertwiner $\phi$ of the 1-isomorphism $\mathscr{P}$ (see Definition~\ref{1-morphism-between-2-line-bundles}) induces an invertible intertwiner\looseness=1
\begin{equation*}
\phi_{ij}\colon\ \mathcal{P}_j \circ s^{*}\mathcal{M}_{s(i),s(j)} \to t^{*}\mathcal{M}_{t(i),t(j)}^{\phi_i} \circ \mathcal{P}_i,
\end{equation*}
where \smash{$\mathcal{M}_{ij}$} is the invertible bimodule bundle over \smash{$U^1_{ij}$} constructed in Remark~\ref{extracting-cocycle}. Under the trivializations $\mathcal{P}_i\cong U_i^1 \times P_i$ and $\mathcal{M}_{ij}\cong U^0_{ij}\times M_{ij}$ (where $M_{ij}$ is an $A$-$A$-bimodule, for $A=A_i=A_j$), we obtain two consequences:
\begin{enumerate}[(1)]\itemsep=0pt
\item
We obtain an invertible intertwiner between the typical fibres,
\begin{equation*}
M_{t(i),t(j)}^{\phi_i} \otimes_{\tilde A_i} P_i \cong P_j \otimes_{\tilde A_j} M_{s(i),s(j)}.
\end{equation*}
Using the isomorphisms $P_i \cong -\sigma_i\tilde A_i$ and $M_{ij}\cong\varepsilon_{ij}A_i=\varepsilon_{ij}A_j$ and the equality $\tilde A_i=A_{s(i)}$, this shows an equality $ \varepsilon_{t(i),t(j)}- \sigma_i=-\sigma_j + \varepsilon_{s(i),s(j)}$, which is the cocycle condition $\delta \sigma-\Delta \epsilon$ of~\ref{cech-cohomology-in-low-degrees:b}.

\item
Under above trivializations, the intertwiner $\phi_{ij}$ becomes a~smooth map $\eta_{ij}\colon U^1_{ij} \to \C^{\times}$. The commutative diagram~\eqref{diag:1Morphisms} for $\phi$ then shows precisely condition~\eqref{graded-equivariant-cohomology:2}, i.e., \smash{$\delta \eta+\Delta^{\phi} \mu=0$}.
\end{enumerate}

Finally, we work with a~cover $\mathcal{U}^2$ of $\Gamma_2$, and there on an open set $U^2_i$. The 2-isomorphism $\psi$ induces an isomorphism
\begin{equation*}
(\operatorname{pr}_1^*\mathcal{P}_{\pr_1(i)})^{\operatorname{pr}_2^*\phi_i} \circ \operatorname{pr}_2^*\mathcal{P}_{\pr_2(i)}
 \to c^*\mathcal{P}_{c(i)}.
\end{equation*}
Under the trivializations $P_i \cong -\sigma_i\tilde A_i$, it becomes an equality
\[
\pr_1^{*}\sigma_{\pr_1(i)}+\pr_2^{*}\sigma_{\pr_2(i)}=c^{*}\sigma_{c(i)},
\]
 which provides the cocycle condition $\Delta \sigma=0$ of Example~\ref{cech-cohomology-in-low-degrees}\,\ref{cech-cohomology-in-low-degrees:b}; thus showing that $(\varepsilon,\sigma)$ defines a~class in $\check \h^1(\Gamma,\Z_2)=\check\h^1((\Gamma,\phi),\Z_2)$. Moreover, we obtain a~smooth map $f_i\colon U_i^2 \to \C^{\times}$ satisfying~\eqref{graded-equivariant-cohomology:3}. Finally, the commutativity of~\eqref{equivariant structure: diagram} implies the remaining cocycle condition~${\Delta^{\phi} f=0}$, showing that $(\mu,\eta,f)$ defines a~class in $\check\h^2((\Gamma,\phi),\sheaf{\C^{\times}})$.

It is straightforward to generalize the reconstruction procedure explained in Example~\ref{reconstruction-of-equivariant-bundle-gerbes} for (ungraded) bundle gerbes to super 2-line bundles, and moreover, to show that this is inverse to the above extraction of local data.
\end{proof}

\section{Twistings and twisted vector bundles}
\label{twistings}

The following is the central definition and conveys the main point of this article.

\begin{Definition} \label{Definition of our twistings}
Let $(\Gamma,\phi)$ be a~graded Lie groupoid. A \emph{twisting} on $(\Gamma,\phi)$ is a~super 2-line bundle $(\mathscr{L},\mathscr{P},\psi)$ over $(\Gamma,\phi)$ in the sense of Definition~\ref{equivariant structure}.
\end{Definition}

The following is the corresponding, canonical definition of a~category of twisted vector bundles; these constitute the ingredients of the twisted K-groups belonging to the twistings of Definition~\ref{Definition of our twistings}.

\begin{Definition}
\label{Definition of our twisted vector bundles}
Let $(\mathscr{L},\mathscr{P},\psi)$ be a~twisting on $(\Gamma,\phi)$. An \emph{$(\mathscr{L},\mathscr{P},\psi)$-twisted super vector bundle} is a~1-morphism $(\mathscr{L},\mathscr{P},\psi) \to (\mathscr{I},\id,\id)$ in the bicategory $\stwoLineBdl(\Gamma,\phi)$. Twisted super vector bundles form the category
\begin{equation*}
\sVectBdl^{(\mathscr{L},\mathscr{P},\psi)}(\Gamma,\phi) := \Hom_{\stwoLineBdl(\Gamma,\phi)}((\mathscr{L},\mathscr{P},\psi), (\mathscr{I},\id,\id)).
\end{equation*}
\end{Definition}

\begin{Remark}
Unwinding Definition~\ref{Definition of our twisted vector bundles}, an $(\mathscr{L},\mathscr{P},\psi)$-twisted super vector bundle is a~1-mor\-phism of super 2-line bundles
$
 \mathscr{B} \colon \mathscr{L} \to \mathscr{I}
$
 over $\Gamma_0$ together with a~2-isomorphism of super 2-line bundles~${
 \eta \colon s^*\mathscr{B} \Rightarrow (t^*\mathscr{B})^\phi \circ \mathscr{P}}
$
 over $\Gamma_1$, such that the diagram
 \begin{equation*}
 \alxydim{@C=4em}{
\mathscr{B}_{s(\gamma_2)} \ar@{=>}[d]_{\id \circ \eta_{\gamma_2}} \ar@{=>}[r]^-{\eta_{\gamma_1\circ \gamma_2}} & \mathscr{B}^{\phi(\gamma_1\circ \gamma_2)}_{t(\gamma_1)} \circ \mathscr{P}_{\gamma_1\circ \gamma_2} \ar@{=>}[d]^{\id \circ \psi_{\gamma_1,\gamma_2}} \\ \mathscr{B}^{\phi(\gamma_2)}_{t(\gamma_2)} \circ \mathscr{P}_{\gamma_2} \ar@{=>}[r]_-{\eta_{\gamma_1}^{\phi(\gamma_2)} \circ \id} & \mathscr{B}_{t(\gamma_1)}^{\phi(\gamma_1)\phi(\gamma_2)} \circ \mathscr{P}_{\gamma_1}^{\phi(\gamma_2)}\circ \mathscr{P}_{\gamma_2}}
\end{equation*}
is commutative for all $(\gamma_1,\gamma_2)\in \Gamma_2$. If $\mathscr{L}=\mathscr{I}$ is the trivial 2-line bundle, then $\mathscr{B}$ is an ordinary super vector bundle.
If the twisting is completely trivial $(\mathscr{I},\id,\id)$, we obtain a~category~$\sVectBdl(\Gamma,\phi) := \sVectBdl^{(\mathscr{I},\id,\id)}(\Gamma,\phi)$ of (untwisted) super vector bundles over $(\Gamma,\phi)$.
\end{Remark}

The purpose of this paper is to show that the following list of versions of twistings that have appeared in the literature consists of special cases of Definitions~\ref{Definition of our twistings} and~\ref{Definition of our twisted vector bundles}.
\begin{enumerate}[(1)]\itemsep=0pt
\item
\label{case-algebra-bundles}
Central simple super algebra bundles~\cite{DK70,Karoubi1968}.
Central simple super algebra bundles for the field $\C$ over a~smooth manifold $M$ are a~special case of Definition~\ref{Definition of our twistings}, see Proposition~\ref{embedding-of-algebra-bundles}:
\begin{equation*}
\cssAlgBdl_\C(M)\subset \stwoLineBdl(M)\cong \stwoLineBdl(M_{\rm dis})= \stwoLineBdl(M_{\rm dis},1).
\end{equation*}
The twisted vector bundles for a~central simple super algebra bundle $\mathcal{A}$ are precisely the super $\mathcal{A}$-module bundles over $M$.

\item
\label{case-bundle-gerbes}
Bundle gerbes~\cite{murray}.
Bundle gerbes over a~smooth manifold $M$ are a~special case of Definition~\ref{Definition of our twistings}, as
\begin{equation*}
\Grb(M) \cong \twoLineBdl(M)\subset \stwoLineBdl(M)\cong \stwoLineBdl(M_{\rm dis})= \stwoLineBdl(M_{\rm dis},1),
\end{equation*}
see Propositions~\ref{embedding-of-bundle-gerbes} and Remark~\ref{graded-equivariant-over-a-manifold}. The twisted vector bundles are precisely the bundle gerbe modules.
Their role as twistings for K-theory has been described in
 \cite{bouwknegt1}. Bundle gerbes are of course just a~special case of the next item, \ref{case-equivariant-gerbes}.

\item
\label{case-equivariant-gerbes}
Equivariant bundle gerbes~\cite{nikolaus1,nikolaus2}.
Equivariant bundle gerbes over Lie groupoids $\Gamma$ are defined by the canonical homotopy limit construction outlined in~\eqref{homotopy-limit-equivariance}, and are included as
\begin{equation*}
\Grb(\Gamma)\cong \twoLineBdl(\Gamma) \subset \stwoLineBdl(\Gamma)=\stwoLineBdl(\Gamma,1),
\end{equation*}
see Definition~\ref{sub-bicategories} and Remark~\ref{graded-equivariant-over-a-manifold}.
A similar construction of equivariant bundle gerbes can be found in~\cite{Murray2017}. Precursors have appeared for the case when $\Gamma$ is an action groupoid, e.g., in~\cite{gawedzki2,gawedzki8,meinrenken1}.

\item
Chern--Simons twistings~\cite{freed6}.
Every multiplicative bundle gerbe over a~Lie group $G$ is also $G$-equivariant with respect to the conjugation action of $G$ on itself; this makes them a~subclass of equivariant bundle gerbes and a~special case of~\ref{case-equivariant-gerbes}. The twisted vector bundles of Definition~\ref{Definition of our twisted vector bundles} form (for a~finite group $G$) the modular tensor category that belongs to the corresponding Chern--Simons theory (i.e., Dijkgraaf--Witten theory). We discuss this in detail in Section~\ref{Cherns-Simons theory twistings}.

\item
\label{case-jandl-gerbes}
Jandl gerbes~\cite{schreiber1}.
Jandl gerbes live over a~real manifold, i.e., a~smooth manifold $M$ with an involution $\sigma\colon M \to M$. They are the same as bundle gerbes $(\mathscr{G},\mathscr{P},\psi)$ over the corresponding graded groupoid $\Gr(M,\sigma)$ (see~\ref{real-manifold}); thus,
\begin{equation*}
\mathscr{J}\mathrm{dl}\Grb(M,\sigma)\cong\Grb(\Gr(M,\sigma))=\twoLineBdl(\Gr(M,\sigma))\subset \stwoLineBdl(\Gr(M,\sigma)).
\end{equation*}
This is a~special case of~\ref{case-equivariant-jandl-gerbes}, which is investigated thoroughly in Section~\ref{equivariant-Jandl-gerbes}.

\item
Real gerbes of Gomi--Thiang~\cite{Gomi2021}.
These are Jandl gerbes $(\mathscr{G},\mathscr{P},\psi)$ where the surjective submersion of $\mathscr{G}$ is an open cover to which the involution lifts. This is a~subclass of~\ref{case-jandl-gerbes}.

\item
\label{Real-bundle-gerbes}
Real bundle gerbes of Hekmati--Murray--Stevenson--Vozzo~\cite{Hekmati}.
These are Jandl gerbes $(\mathscr{G},\mathscr{P},\psi)$ where the 1-isomorphism $\mathscr{P}$ is (induced from) a~refinement (see Definition~\ref{refinements-of-super-2-line-bundles}). This is another subclass of~\ref{case-jandl-gerbes}.

\item
\label{case-equivariant-jandl-gerbes}
Equivariant Jandl gerbes~\cite{gawedzki8}.
Equivariant Jandl gerbes live over a~smooth manifold $M$ with an action of a~finite graded group $(G,\varepsilon)$, i.e., a~group $G$ with a~homomorphism $\varepsilon\colon G \to \Z_2$, and we prove in Section~\ref{equivariant-Jandl-gerbes} that they are the same as bundle gerbes over the corresponding graded action groupoid~${\act M{(G,\varepsilon)}}$ (see Example~\ref{ex1}\,\ref{graded-action-groupoid}),
\begin{equation*}
\mathrm{GSW}_M(G,\varepsilon)\cong\Grb(\act M{(G,\varepsilon)}) \cong \twoLineBdl(\act M{(G,\varepsilon)}) \subset \stwoLineBdl(\act M{(G,\varepsilon)}).
\end{equation*}
For $(G,\varepsilon)=(\Z_2,\id)$, this is~\ref{case-jandl-gerbes}. For arbitrary $G$ and $\varepsilon=1$, it is a~special case of~\ref{case-equivariant-gerbes}.

\item
Moutuou's RG bundle gerbes~\cite{mohamedmoutuou2012}.
These are a~super version of~\ref{Real-bundle-gerbes}; by a~straightforward generalization of the discussion of~\ref{case-equivariant-jandl-gerbes} in Section~\ref{equivariant-Jandl-gerbes}, they embed into super 2-line bundles over
$\Gr(M,\sigma)$.

 \item
\label{case-FHT}
Twistings of Freed--Hopkins--Teleman~\cite{Freed2011a}.
Freed--Hopkins--Teleman's twistings live over (ungraded) Lie groupoids $\Gamma$, and are the same as (equivariant) super bundle gerbes,
$
\mathrm{FHT}(\Gamma) \cong \hc 1(\sGrb(\Gamma))\subset \hc 1(\stwoLineBdl(\Gamma))$.
This needs again more elaboration, which is carried out in Section~\ref{Freed--Hopkins--Teleman-twistings}.

\item
Twistings of \'Angel--G\'omez--Uribe~\cite{Angel2018}.
These live over action groupoids $\act XG$, and are a~special
case of~\ref{case-FHT}.

\item
\label{case-Freed-Moore-extensions}
Freed--Moore's twisted groupoid extensions~\cite{Freeda}.
These live over graded Lie groupoids $(\Gamma,\phi)$, and they are the same as super 2-line bundles~$(\mathscr{L},\mathscr{P},\psi)$ over $(\Gamma,\phi)$ whose underlying super 2-line bundle $\mathscr{L}$ is trivial; while the morphisms correspond to refinements
$
\mathrm{FM}(\Gamma,\phi)\cong \stwoLineBdltriv^{\text{ref}}(\Gamma,\phi) \subset \stwoLineBdl(\Gamma,\phi)
$.
This needs more elaboration, which is carried out in Section~\ref{freed-moore-extensions}. The twisted vector bundles of~\cite{Freeda} coincide precisely with our Definition~\ref{Definition of our twisted vector bundles}.

\item
Magnetic equivariant vector bundles of Serrano--Uribe--Xicotencatl~\cite{serrano2024rationalmagneticequivariantktheory}.
These are a~special case of our twisted vector bundles from Definition~\ref{Definition of our twisted vector bundles}, reduced to graded Lie groupoids induced from graded Lie groups (see Example~\ref{ex1}\,\ref{graded-Lie-group}), and reduced to trivial twistings. This is a~special case of~\ref{case-Freed-Moore-extensions}.

\item
\label{Moutuou}
Moutuou's real twistings~\cite{mohamedmoutuou2012}.
Moutuou's real twistings live over real groupoids $(\Gamma,\tau)$, i.e., Lie groupoids with an involutive functor $\tau$. They can be viewed as a~real analogue of Freed--Hopkins--Teleman's twistings of~\ref{case-Freed-Moore-extensions}. Indeed, associated to any real groupoid $(\Gamma,\tau)$ is a~graded Lie groupoid $\Gr(\Gamma,\tau)$, and every
real twisting over $(\Gamma,\tau)$ gives rise to a~super 2-line bundle over $\Gr(\Gamma,\tau)$. This is discussed thoroughly in Section~\ref{Moutuou Real twistings}.

\item
Berwick--Evans--Guo real twistings~\cite{BerwickEvans2024}.
These twistings are defined similarly to Moutuou's real twistings of~\ref{Moutuou} but using topological stacks instead of topological groupoids. A~comparison between both approaches can be found in \cite[Sections~3.4 and~3.5]{BerwickEvans2024}.

\item
\label{case-Freed}
Freed's twisted invertible algebra bundles~\cite{Freed2012a}.
These live over Lie groupoids $\Gamma$, and involve a~certain ``double cover'' by a~graded Lie groupoid. Under the assumption that the grading of the double cover is induced by a~grading $\phi$ on the Lie groupoid $\Gamma$,
Freed's twisted invertible algebra bundles can be identified with super 2-line bundles over $(\Gamma,\phi)$. This is discussed in detail in Section~\ref{Freed's invertible algebra bundle twistings}.

\item
\label{case-Distler--Freed--Moore-extensions}
Orientifold twistings of Distler--Freed--Moore~\cite{Distler2011,distler2010spin,Freed2009Slides}.
Distler--Freed--Moore's orientifold twistings are similar to~\ref{case-Freed}, but feature instead of algebra bundles a~certain $\Z$-valued grading function $d$. Under the same limitation as in~\ref{case-Freed},
these twistings yield super 2-line bundles over $(\Gamma,\phi)$. This is discussed in detail in Section~\ref{Distler Freed Moore twistings}.

\item
Gomi--Freed--Moore twistings~\cite{gomi2021freedmoorektheory}.
These live over graded groupoids $(\Gamma,\phi)$, and they are the same as the Distler--Freed--Moore-twistings of~\ref{case-Distler--Freed--Moore-extensions}
with two constraints: first, they satisfy our assumption that the grading of the double cover is induced from one on the base, and second, that the grading function~$d$ is trivial. This is a~special case of
\ref{case-Distler--Freed--Moore-extensions}, and yields all super bundle gerbes over~$(\Gamma,\phi)$.
\end{enumerate}

\section{Comparison work}\label{section-comparison-work}

In this section, we carry out the comparison work announced in Section~\ref{twistings}.

\subsection{Equivariant Jandl gerbes}
\label{equivariant-Jandl-gerbes}

The paper~\cite{gawedzki8} introduced a~general theory for B-fields in unoriented orbifold sigma models. The B-field is modelled by a~so-called equivariant Jandl gerbe, combining an orbifolding process with the treatment of unoriented worldsheets. The following notion was introduced.

\begin{Definition}
 An \emph{orientifold group} ($G,\varepsilon$) for a~smooth manifold $X$ is a~finite group $G$ acting smoothly from the left on $X$,
 together with a~group homomorphism $\varepsilon\colon G \to \mathbb{Z}_2$.
\end{Definition}

This embeds into the more general notion of graded Lie groupoids (see Definition~\ref{graded-groupoid}) by associating to an orientifold group $(G,\varepsilon)$ the action groupoid $\act XG$ together with the map $\phi_{\varepsilon}\colon G \times X \to \Z_2$, $(g,x)\mapsto \varepsilon(g)$.
We prove in Theorem~\ref{comparison-equivariant-jandl-gerbes} below that the equivariant Jandl gerbes of~\cite{gawedzki8} for an orientifold group $(G,\varepsilon)$ are the same as bundle gerbes over the graded Lie groupoid $(\act XG,\phi_{\varepsilon})$, and hence, in particular, 2-line bundles over that groupoid.

For the geometric part, we consider the bicategory $\hGrb(X) $ of hermitian bundle gerbes over~$X$, and therein the bigroupoid $\hGrb^{\times}(X)$ that contains only the invertible 1-morphisms and 2-morphisms. Forgetting the hermitian metrics yields the vertical functors in the following commutative diagram
\begin{equation*}
\alxydim{}{\hGrb^{\times}(X) \ar[r] \ar[d] & \hGrb(X) \ar[d] \\ \Grb^{\times}(X) \ar[r] & \Grb(X) \subset \stwoLineBdl(X). \hspace{-6em}}
\end{equation*}
In~\cite{gawedzki8}, the following constructions are made. Suppose $\mathscr{G}=(\pi,\mathcal{M},\mu)$ is a~hermitian bundle gerbe, i.e., $\mathcal{M}$ is a~hermitian line bundle and $\mu$ is an unitary bundle isomorphism. Then, $\mathcal{G}^{\dagger}:=\bigl(\pi,\mathcal{M}^{*},\mu^{tr-1}\bigr)$ denotes the (usual) dual bundle gerbe: $\mathcal{M}^{*}$ is the dual line bundle, equipped with the dual metric, and $\mu^{tr-1}$ is the inverse of the transpose of $\mu$. If $\mathscr{A}=(\zeta, \mathcal{P},\phi) \colon \mathscr{G}_1 \to \mathscr{G}_2$ is a~hermitian bundle gerbe isomorphism, then we set $\mathscr{A}^{\dagger}:=\bigl(\zeta,\mathcal{P}^{*},\phi^{tr-1}\bigr) \colon \mathscr{G}_1^{\dagger} \to \mathscr{G}_2^{\dagger}$. Thus, $\mathscr{A}^{\dagger}$ has the same surjective submersion $\zeta\colon Z \to Y_1 \times_X Y_2$ and the dual hermitian line bundle $\mathcal{P}^{*}$. We note that $\phi$ is a~unitary bundle isomorphism, over a~point $(z,z')\in Z^{[2]}$ with $\zeta(z)=:(y_1,y_2)$ and $\zeta(z')=:(y_1',y_2')$,
\begin{equation*}
\phi_{z,z'}\colon\ \mathcal{P}_{z'} \otimes (\mathcal{M}_1)_{y_1,y_1'} \to (\mathcal{M}_2)_{y_2,y_2'} \otimes \mathcal{P}_{z}.
\end{equation*}
Up to a~reordering of tensor factors, this coincides with the conventions of~\cite{gawedzki8}. Thus, $\phi^{tr-1}$ is an isomorphism
\begin{equation*}
\phi_{z,z'}^{tr-1}\colon\ \mathcal{P}^{*}_{z'} \otimes (\mathcal{M}_1^{*})_{y_1,y_1'} \to (\mathcal{M}_2^{*})_{y_2,y_2'} \otimes \mathcal{P}_{z}^{*},
\end{equation*}
and one checks that it satisfies the required coherence property so that $\mathscr{A}^\dagger$ is well-defined. Finally, if $\beta\colon \mathscr{A}_1 \Rightarrow \mathscr{A}_2$ is a~2-isomorphism, then $\beta^{\dagger} :=\beta^{tr-1}\colon \mathscr{A}_1^{\dagger} \Rightarrow \mathscr{A}_2^{\dagger}$ is a~new 2-isomorphism. These constructions define a~2-functor
\begin{equation*}
(..)^{\dagger}\colon\ \hGrb^{\times}(X) \to \hGrb^{\times}(X)
\end{equation*}
that is involutive under a~canonical natural isomorphism $(..)^{\dagger\dagger} \cong \id$.

The definitions in~\cite{gawedzki8} are based on the 2-functor $(..)^{\dagger}$. Before we recall those, we want to relate $(..)^{\dagger}$ to our complex conjugation functor $\overline{(..)}$ of~\eqref{complex-conjugation-2-functor}.

\begin{Lemma}
\label{GSW-dagger-and-conjugation}
The two $2$-functors
\begin{equation*}
\alxydim{}{\hGrb^{\times}(X) \to \Grb^{\times}(X) \subset \twoLineBdl(X) \ar[r]^-{\overline{(..)}} & \twoLineBdl(X)}
\end{equation*}
and
\begin{equation*}
\alxydim{}{\hGrb^{\times}(X) \ar[r]^-{(..)^{\dagger}} & \hGrb^{\times}(X) \to \Grb^{\times}(X) \subset \twoLineBdl(X)}
\end{equation*}
are canonically naturally isomorphic. In short: upon forgetting the hermitian metrics, the functor $(..)^{\dagger}$ is complex conjugation.
\end{Lemma}

\begin{proof}
We recall that an inner product $\left \langle -,- \right \rangle$ on a~finite-dimensional complex vector space $V$ (antilinear in the \emph{first} component) induces an isomorphism
\begin{equation}
\label{isomorphism-between-complex-conjugate-and-dual}
\overline{V} \to V^{*}\colon\ v \mapsto \left \langle v,- \right \rangle.
\end{equation}
Moreover, $\overline{V}$ is equipped with an inner product $\overline{\left \langle -,- \right \rangle}$, given by the complex conjugate of the original one. Also, $V^{*}$ is equipped with an inner product, which is usually induced from $\overline{V}$ along the isomorphism~\eqref{isomorphism-between-complex-conjugate-and-dual}, so that~\eqref{isomorphism-between-complex-conjugate-and-dual} actually becomes an \emph{isometric} isomorphism. Moreover, for an isometric isomorphism $\varphi\colon V \to W$, the relevant diagram
\begin{equation}
\label{diagram-duality}
\alxydim{}{\overline{V} \ar[d]_{\overline{\varphi}} \ar[r] & V^{*} \ar[d]^{\varphi^{tr-1}} \\ \overline{W} \ar[r] & W^{*} }
\end{equation}
is commutative. The isometric isomorphism~\eqref{isomorphism-between-complex-conjugate-and-dual} extends to hermitian vector bundles, and a~diagram analogous to~\eqref{diagram-duality} exists for unitary isomorphisms between hermitian vector bundles.

If $\mathscr{G}=(\pi,\mathcal{M},\mu)$ is a~hermitian bundle gerbe, then we obtain a~canonical isomorphism $\eta_{\mathscr{G}}\colon \overline{\mathscr{G}} \to \mathscr{G}^{\dagger}$, which is the refinement consisting of the identity $\id\colon Y \to Y$ and the unitary line bundle isomorphism $u_{\mathcal{M}}\colon \overline{\mathcal{M}} \to \mathcal{M}^{*}$. The commutativity of~\eqref{diagram-duality} shows the one of the relevant diagram~\eqref{condition-for-refinements}.

If $\mathscr{A}\colon \mathscr{G}_1 \to \mathscr{G}_2$ is a~1-isomorphism, then we construct a~2-isomorphism
\begin{equation*}
\alxydim{}{\overline{\mathscr{G}_1} \ar[d]_{\overline{\mathscr{A}}} \ar[r]^{\eta_{\mathscr{G}_1}} & \mathscr{G}_1^{\dagger} \ar@{<=}[dl]|*+{\eta_{\mathscr{A}}} \ar[d]^{\mathscr{A}^{\dagger}} \\ \overline{\mathscr{G}_2} \ar[r]_{\eta_{\mathscr{G}_2}} & \mathscr{G}_2^{\dagger}. }
\end{equation*}
We have to compute the two composites. Suppose $\mathscr{A}=(\zeta,\mathcal{P},\phi)$, so that $\overline{\mathscr{A}}=\bigl(\zeta,\overline{\mathcal{P}},\overline{\phi}\bigr)$ and $\mathscr{A}^{\dagger}=\bigl(\zeta,\mathcal{P}^{*},\phi^{tr-1}\bigr)$. Then, we get $\mathscr{A}^{\dagger} \circ \eta_{\mathscr{G}_1}=\bigl(\zeta, \mathcal{P}^{*},\phi^{tr-1} \circ (\id \otimes u_{\mathcal{M}_1})\bigr)$, with the morphism being
\begin{equation*}
\alxydim{@C=5em}{\mathcal{P}^{*}_{z'} \otimes \bigl(\overline{\mathcal{M}_1}\bigr)_{y_1,y_1'} \ar[r]^{\id \otimes u_{\mathcal{M}_1}} & \mathcal{P}^{*}_{z'} \otimes \bigl(\mathcal{M}_1^{*}\bigr)_{y_1,y_1'} \ar[r]^-{\phi_{z,z'}^{tr-1}} & \bigl(\mathcal{M}_2^{*}\bigr)_{y_2,y_2'} \otimes \mathcal{P}_{z}^{*},}
\end{equation*}
as well as $\eta_{\mathscr{G}_2} \circ \overline{\mathscr{A}}=\bigl(\zeta, \overline{\mathcal{P}},(u_{\mathcal{M}_2} \otimes \id) \circ \overline{\phi}\bigr)$, with the morphism being
\begin{equation*}
\alxydim{@C=5em}{ \overline{\mathcal{P}}_{z'} \otimes \bigl(\overline{\mathcal{M}_1}\bigr)_{y_1,y_1'} \ar[r]^-{\overline{\phi}_{z,z'}} & \bigl(\overline{\mathcal{M}_2}\bigr)_{y_2,y_2'} \otimes \overline{\mathcal{P}}_{z} \ar[r]^{u_{\mathcal{M}_2} \otimes \id} & \bigl(\mathcal{M}_2^{*}\bigr)_{y_2,y_2'} \otimes \overline{\mathcal{P}}_{z}.}
\end{equation*}
The claimed 2-morphism is given by
$
\eta_{\mathscr{A}}:= u_{\mathcal{P}}\colon \overline{\mathcal{P}} \to \mathcal{P}^{*}$.
The relevant commutative diagram is the outer shape of the diagram
\begin{equation*}
\alxydim{@C=5em}{
\overline{\mathcal{P}}_{z'} \otimes \bigl(\overline{\mathcal{M}_1}\bigr)_{y_1,y_1'} \ar[dr]_{u_{\mathcal{P} \otimes \mathcal{M}_1}} \ar[d]_{u_{\mathcal{P}} \otimes \id} \ar[r]^-{\overline{\phi}_{z,z'}} & \bigl(\overline{\mathcal{M}_2}\bigr)_{y_2,y_2'} \otimes \overline{\mathcal{P}}_{z} \ar[dr]^{u_{\mathcal{M}_2 \otimes \mathcal{P}}} \ar[r]^{u_{\mathcal{M}_2} \otimes \id} & \bigl(\mathcal{M}_2^{*}\bigr)_{y_2,y_2'} \otimes \overline{\mathcal{P}}_{z} \ar[d]^{\id \otimes u_{\mathcal{P}}}
\\
\mathcal{P}^{*}_{z'} \otimes \bigl(\overline{\mathcal{M}_1}\bigr)_{y_1,y_1'} \ar[r]_{\id \otimes u_{\mathcal{M}_1}} & \mathcal{P}^{*}_{z'} \otimes \bigl(\mathcal{M}_1^{*}\bigr)_{y_1,y_1'} \ar[r]_-{\phi_{z,z'}^{tr-1}} & \bigl(\mathcal{M}_2^{*}\bigr)_{y_2,y_2'} \otimes \mathcal{P}_{z}^{*},
}
\end{equation*}
which is commutative because all subdiagrams are commutative: the triangular diagrams because
\eqref{isomorphism-between-complex-conjugate-and-dual} is compatible with tensor products, and the one in the middle is~\eqref{diagram-duality}. One can now check that $\eta_{\mathscr{G}}$ and $\eta_{\mathscr{A}}$ are the components of a~natural isomorphism $\eta\colon \overline{(..)}|_{\hGrb^{\times}(X)} \Rightarrow (..)^{\dagger}$.
\end{proof}

By Lemma~\ref{GSW-dagger-and-conjugation}, we can replace occurrences of $(..)^{\dagger}$ in~\cite{gawedzki8} consistently by complex conjugation. Next we recall the following definition from~\cite{gawedzki8}. Given an orientifold group $(G,\varepsilon)$ and a~bundle gerbe $\mathscr{G}$ over $X$, we set for $g \in G$
\smash{$
 g \mathscr{G} := \bigl(g^{-1}\bigr)^*\mathscr{G}^{\varepsilon(g)}$},
 where $g^{-1}\colon X \to X$ is the diffeomorphism assigned to $g^{-1}$ via the action of $G$ on $X$, and
 $(..)^{\varepsilon(g)}$ refers to the convention fixed in Section~\ref{graded-equivariant-definitions}, i.e., it means nothing when $\varepsilon(g)=1$ and it means complex conjugation when $\varepsilon(g)=-1$.
For a~1-isomorphism $\mathscr{A} \colon \mathscr{G} \to \mathscr{H}$ between bundle gerbes, we define
\[
 g \mathscr{A} := \bigl(g^{-1}\bigr)^*\mathscr{A}^{\varepsilon(g)}\colon\ g\mathscr{G} \to g\mathscr{H}.
\]
Finally, for a~2-isomorphism $\beta\colon \mathscr{A} \Rightarrow \mathscr{A}'$, we define
\[
 g \beta :=\bigl(g^{-1}\bigr)^*\beta^{\varepsilon(g)}
 \colon\ g \mathscr{A} \Rightarrow g \mathscr{A}'.
\]
These definitions promote each group element $g\in G$ to a~functor $F_g\colon \Grb(X) \to \Grb(X)$ that is covariant at all levels and satisfies $F_{g_2} \circ F_{g_1}\cong F_{g_2g_1}$. The following is \cite[Definition~2.9]{gawedzki8}.

\begin{Definition}
\label{GSW-equivariant-structure-on-bundle-gerbe}
 Let $(G,\varepsilon)$ be an orientifold group for $X$ and let $\mathscr{G}$ be a~bundle gerbe over $X$.
 A \emph{$(G,\varepsilon)$-equivariant structure on $\mathscr{G}$} consists of 1-isomorphisms
$
 \mathscr{A}_g\colon \mathscr{G} \to g\mathscr{G}
$
 for each $g \in G$, and of 2-isomorphisms
$
 \varphi_{g_1,g_2} \colon g_1\mathscr{A}_{g_2} \circ \mathscr{A}_{g_1} \to \mathscr{A}_{g_1 g_2}
$
 for each pair $g_1,g_2 \in G$, such that the diagram
 \[
 \begin{aligned} \begin{xy}
 \xymatrix{
 g_1 g_2 \mathscr{A}_{g_3} \circ g_1 \mathscr{A}_{g_2} \circ \mathscr{A}_{g_1} \ar@{=>}[rr]^-{\id \circ \varphi_{g_1,g_2}} \ar@{=>}[d]_-{g_1 \varphi_{g_2,g_3} \circ \id} & & g_1 g_2 \mathscr{A}_{g_3} \circ\mathscr{A}_{g_1 g_2} \ar@{=>}[d]^-{\varphi_{g_1g_2, g_3}}\\
 g_1 \mathscr{A}_{g_2 g_3} \circ \mathscr{A}_{g_1} \ar@{=>}[rr]_-{\varphi_{g_1,g_2g_3}} & & \mathscr{A}_{g_1 g_2 g_3}
 }
 \end{xy} \end{aligned}
 \]
 is commutative for all $g_1,g_2,g_3 \in G$. We call the triple $(\mathscr{G},\mathscr{A},\varphi)$ a~\emph{$(G,\varepsilon)$-equivariant bundle gerbe}.
\end{Definition}
The following is \cite[Definition~2.11]{gawedzki8}.
\begin{Definition}
 Let $(G, \varepsilon)$ be an orientifold group for $X$ and let $\bigl(\mathscr{G}^1,\mathscr{A}^1,\varphi^1\bigr)$ and $\bigl(\mathscr{G}^2,\mathscr{A}^2,\varphi^2\bigr)$
 be $(G,\varepsilon)$-equivariant bundle gerbes. A \emph{$(G,\varepsilon)$-equivariant 1-morphism} $\bigl(\mathscr{G}^1, \mathscr{A}^1,\varphi^1\bigr) \to \bigl(\mathscr{G}^2, \mathscr{A}^2,\varphi^2\bigr)$
 is a~1-morphism $\mathscr{B} \colon \mathscr{G}^1 \to \mathscr{G}^2$ of the underlying bundle gerbes together with a~family of 2-isomorphisms
$
 \eta_g \colon g \mathscr{B} \circ \mathscr{A}_g^1 \Rightarrow \mathscr{A}_g^2 \circ \mathscr{B}$,
 one for each $g \in G$, such that the diagram
 \[
 \begin{aligned} \begin{xy}
 \xymatrix{
 g_1g_2 \mathscr{B} \circ g_1\mathscr{A}_{g_2}^1 \circ \mathscr{A}^1_{g_1} \ar@{=>}[rr]^-{\id\circ \varphi_{g_1,g_2}^1 } \ar@{=>}[d]_-{g_1 \eta_{g_2}\circ \id} & & g_1g_2\mathscr{B} \circ \mathscr{A}^1_{g_1g_2} \ar@{=>}[dd]^-{\eta_{g_1g_2}} \\
 g_1 \mathscr{A}_{g_2}^2 \circ g_1 \mathscr{B} \circ \mathscr{A}_{g_1}^1 \ar@{=>}[d]_-{\id\circ \eta_{g_1}}& & \\
 g_1\mathscr{A}_{g_2}^2 \circ \mathscr{A}^2_{g_1} \circ \mathscr{B} \ar@{=>}[rr]_-{\varphi_{g_1,g_2}^2\circ \id} & & \mathscr{A}^2_{g_1g_2} \circ \mathscr{B}
 }
 \end{xy} \end{aligned}
 \]
 of 2-isomorphisms is commutative for all $g_1,g_2 \in G$.
\end{Definition}
We also recall the notion of 2-morphisms from~\cite{gawedzki8}.
\begin{Definition}
 Suppose we have two $(G,\varepsilon)$-equivariant bundle gerbes $\bigl(\mathscr{G}^1,\mathscr{A}^1,\psi^1\bigr)$ and $\bigl(\mathscr{G}^2,\mathscr{A}^2,\allowbreak\psi^2\bigr)$
 and two equivariant 1-morphisms $(\mathscr{B},\eta)$ and $(\mathscr{B}',\eta')$ between them. A \emph{$(G,\varepsilon)$-equivariant $2$-morphism}
$(\mathscr{B},\eta) \Rightarrow (\mathscr{B}',\eta')$
 is a~2-morphism $\xi\colon \mathscr{B} \Rightarrow \mathscr{B}'$ which is compatible with the 2-isomorphisms $\eta_g$ and $\eta_g'$ in the sense that the diagram
 \[
 \begin{aligned} \begin{xy}
 \xymatrix{
 g\mathscr{B} \circ \mathscr{A}^1_g\ar@{=>}[rr]^-{\eta_g} \ar@{=>}[d]_-{g\xi \circ \id} & & \mathscr{A}^2_g \circ \mathscr{B} \ar@{=>}[d]^-{\id\circ\xi} \\
 g\mathscr{B}' \circ \mathscr{A}^1_g \ar@{=>}[rr]_-{\eta'_g}& & \mathscr{A}^2_g \circ \mathscr{B}'
 }
 \end{xy} \end{aligned}
 \]
 of 2-morphisms commutes for all $g \in G$.
\end{Definition}

$(G,\varepsilon)$-equivariant bundle gerbes in the sense of~\cite{gawedzki8} form a~bicategory that we denote here for one moment by $\mathrm{GSW}_X(G,\varepsilon)$.

\begin{Theorem}
\label{comparison-equivariant-jandl-gerbes}
Let $(G,\varepsilon)$ be an orientifold group for a~smooth manifold $X$, and let $(\act XG,\phi_{\varepsilon})$ be the associated graded Lie groupoid.
 The bicategory of $(G,\varepsilon)$-equivariant bundle gerbes of~{\rm\cite{gawedzki8}} is equivalent to the bicategory of bundle gerbes over $(\act XG,\phi_{\varepsilon})$,
\begin{equation*}
 \mathrm{GSW}_X(G,\varepsilon)\cong \Grb(\act XG,\phi_{\varepsilon}) \cong \twoLineBdl(\act XG,\phi_{\varepsilon})\subset \stwoLineBdl(\act XG,\phi_{\varepsilon}).
\end{equation*}
\end{Theorem}

\begin{proof}
The proof consists mainly of a~conversion between the different conventions.
If $\mathscr{G}$ is a~bundle gerbe over $X$, and $(\mathscr{A}_{g},\varphi_{g_1,g_2})$ is a~$(G,\varepsilon)$-equivariant structure on $\mathscr{G}$ in the sense of Definition~\ref{GSW-equivariant-structure-on-bundle-gerbe}, then, using that $G$ is finite,
\begin{equation*}
\mathscr{L}:= \mathscr{G}
,\qquad
\mathscr{P}|_{\{g\} \times X} :=\mathscr{A}_{g^{-1}}
\quand
\psi|_{\{g_1\}\times \{g_2\} \times X }:= \varphi_{g_2^{-1},g_1^{-1}}
\end{equation*}
is a~bundle gerbe over $(\act XG,\phi_{\varepsilon})$ in the sense of Definition~\ref{sub-bicategories}.
Next, if $(\mathscr{B},\eta_g) \colon \bigl(\mathscr{G}^1, \mathscr{A}^1,\varphi^1\bigr) \to \bigl(\mathscr{G}^2, \mathscr{A}^2,\varphi^2\bigr)$ is a~$(G,\varepsilon)$-equivariant 1-morphism, then
$
 \mathscr{B}$, \smash{$\eta|_{\{g\} \times X} := \eta_{g^{-1}}^{-1}$}
is a~1-morphism of the associated bundle gerbes $\bigl(\mathscr{L}^1,\mathscr{P}^1,\psi^1\bigr) \to \bigl(\mathscr{L}^2,\mathscr{P}^2,\psi^2\bigr)$ over $(\act XG,\phi_{\varepsilon})$
in the sense of Definition~\ref{sub-bicategories}. Similarly, if $\xi\colon (\mathscr{B},\eta_g) \to (\mathscr{B}',\eta_g')$ is a~$(G,\varepsilon)$-equivariant 2-morphism,
then $\xi\colon (\mathscr{B},\eta) \to (\mathscr{B}',\eta')$ is a~2-morphism of the associated 1-morphisms of bundle gerbes over $(\act XG,\phi_{\varepsilon})$
in the sense of Definition~\ref{sub-bicategories}. This defines a~functor $\mathrm{GSW}_X(G,\varepsilon) \to \Grb(\act XG,\phi_{\varepsilon})$ which is now easily checked to be an equivalence.
\end{proof}

\begin{Remark}
\cite[Proposition~3.1]{gawedzki8} gives a~cohomological classification of equivariant Jandl gerbes (with hermitian metrics and connections), by a~``twisted-equivariant'' version of smooth Deligne cohomology. Stripping off the connection data, the cocycle data for an equivariant Jandl gerbe satisfies conditions (3.3), (3.8), (3.10) and (3.11) in~\cite{gawedzki8}. Under slight change of conventions, one can check that these coincide with the conditions for a~2-cocycle in $\check\h^2((\act XG,\phi_{\varepsilon}),\sheaf{\C^{\times}})$ listed in Example~\ref{reconstruction-of-equivariant-bundle-gerbes}. Indeed, under the equivalence of Theorem~\ref{comparison-equivariant-jandl-gerbes}, Example~\ref{reconstruction-of-equivariant-bundle-gerbes} is precisely the classification of equivariant Jandl gerbes (without metrics and connections).
\end{Remark}

\begin{Remark}
Descent for equivariant Jandl gerbes has already been considered in~\cite{gawedzki8}.
We consider an orientifold group $(G,\varepsilon)$ for a~smooth manifold $X$, in such a~way that the normal subgroup $G^{0}:=\mathrm{ker}(\epsilon)$ acts without fixed points. Thus, the quotient $X': =X/G^{0}$ is again a~smooth manifold, for which the remaining group $G':=G/G^{0}$ together with its induced grading~${\epsilon'\colon G' \to \Z_2}$ is an orientifold group. Then, \cite[Theorem~3.1]{gawedzki8} proves an equivalence
\begin{equation}
\label{equivalence-GSW}
\mathrm{GSW}_X(G,\epsilon)\cong \mathrm{GSW}_{X'}\bigl(G',\varepsilon'\bigr).
\end{equation}
In our current framework, this result follows from Theorem~\ref{descent-for-graded-equivariant-2-line-bundles} and our equivalence of Theorem~\ref{comparison-equivariant-jandl-gerbes}.
To see this, we consider the graded Lie groupoids $(\act XG,\phi_\varepsilon)$ and $(\act {X'}{G'},\phi_{\varepsilon'})$, and the obvious projection functor
$
(\act XG,\phi_\varepsilon) \to (\act {X'}{G'},\phi_{\varepsilon'})$.
This functor is even and one can easily prove that it is a~weak equivalence, in fact, a~covering functor. Hence, Theorem~\ref{descent-for-graded-equivariant-2-line-bundles} implies the equivalence~\eqref{equivalence-GSW}. \end{Remark}

\subsection{Freed--Moore's twisted groupoid extensions}
\label{freed-moore-extensions}

Freed and Moore introduced in~\cite{Freeda} so-called twisted groupoid extensions to describe symmetry protected phases in condensed matter physics.
We show (see Theorem~\ref{The main theorem of the Freed-Moore chapter}) that Freed and Moore's twisted groupoid extensions are the same as
graded-equivariant super 2-line bundles whose underlying 2-line bundle is trivial; hence, they are particular graded-equivariant super 2-line bundles.

We recall that we denote by $\overline{\mathcal{L}}$ the complex conjugate of a~complex line bundle $\mathcal{L}$.
We fix the following notation. Let $\mathcal{L}$ be a~line bundle and $\phi=\pm 1$. Then,
\begin{equation} \label{complex conjugation according to sign}
 \mathcal{L}^\phi = \begin{cases}
 \mathcal{L}, & \phi = +1, \\
 \overline{\mathcal{L}}, & \phi = -1.
 \end{cases}
\end{equation}

\begin{Definition}
\label{Definition: twisted super extension}
 Let $(\Gamma,\phi)$ be a~graded Lie groupoid.
 A \emph{$\phi$-twisted super extension of $\Gamma$} is a~pair $(\mathcal{L},\lambda)$ consisting
 of a~super line bundle $\mathcal{L}$ over $\Gamma_1$ together with an isomorphism
 \[
 \lambda\colon \ \operatorname{pr}_1^*\mathcal{L}^{\operatorname{pr}_2^*\phi} \otimes \operatorname{pr}_2^*\mathcal{L}
 \to c^*\mathcal{L}
 \]
 of super line bundles over $\Gamma_2$, such that the diagram
 \begin{equation} \label{Definition: twisted super extension diagram}
 \begin{gathered}
 \begin{xy}
 \xymatrix@C=1em{
 \bigl(\mathcal{L}_{\gamma_1}^{\phi(\gamma_2)} \otimes \mathcal{L}_{\gamma_2}\bigr)^{\phi(\gamma_3)} \otimes \mathcal{L}_{\gamma_3} \ar[rr] \ar[d]_-{\lambda_{\gamma_1,\gamma_2}^{\phi(\gamma_3)} \otimes \id} && \mathcal{L}_{\gamma_1}^{\phi(\gamma_2) \cdot \phi(\gamma_3)} \otimes \bigl(\mathcal{L}_{\gamma_2}^{\phi(\gamma_2)} \otimes \mathcal{L}_{\gamma_3}\bigr) \ar[d]^-{\id \otimes \lambda_{\gamma_2,\gamma_3}}\\
 \mathcal{L}_{\gamma_1\circ \gamma_2}^{\phi(\gamma_3)} \otimes \mathcal{L}_{\gamma_3} \ar[dr]_{\lambda_{\gamma_1\circ \gamma_2,\gamma_3}} && \mathcal{L}_{\gamma_1}^{\phi(\gamma_2 \circ \gamma_3) } \otimes \mathcal{L}_{\gamma_2\circ \gamma_3} \ar[dl]^{\lambda_{\gamma_1,\gamma_2\circ \gamma_3}} \\
 & \mathcal{L}_{\gamma_1\circ \gamma_2\circ \gamma_3}
 }
 \end{xy}
 \end{gathered}
 \end{equation}
 commutes for each composable triple of morphisms $(\gamma_1,\gamma_2,\gamma_3)\in \Gamma_3$.
\end{Definition}

Definition~\ref{Definition: twisted super extension} coincides with \normalfont\cite[Definition~7.23\,(iii)]{Freeda} under the dropping of hermitian structures and some notational reformulations. First, our notion of a~graded Lie groupoid coincides with Freed--Moore's notion of a~``groupoid with a~homomorphism to ${\pm 1}$''.
Second, our super line bundle $\mathcal{L}$ is described there as a~line bundle with a~map $\vartheta\colon \Gamma_1 \to \{\pm 1\}$.
Third, our notation~$\mathcal{L}^{\phi}$ (defined in~\eqref{complex conjugation according to sign}) coincides with Freed--Moore's notation $^{\phi}\mathcal{L}$ defined in \cite[equation~(7.9)]{Freeda}.
 Then, the isomorphism $\lambda$ is exactly as in \cite[equation~(7.10)]{Freeda}, and the diagram is precisely \cite[equation~(7.11)]{Freeda}.

\begin{Example} For any graded Lie groupoid $(\Gamma,\phi)$, there is a~trivial $\phi$-twisted super extension~$I_{\Gamma,\phi}$ consisting of the trivial line bundle $\C \times \Gamma_1$, considered as purely even, and of the isomorphism~$\lambda_{triv}$ given in each fibre by $ \C^\phi \otimes \C \cong \C \otimes \C \to \C$, using the canonical isomorphism~$\C\cong \overline\C$ over points with $\phi=-1$.
\end{Example}

\begin{Example}\label{twisted-groupoid-extension-of-BG}
 Let $(G,\varepsilon)$ be a~graded Lie group (see Example~\ref{ex1}\,\ref{graded-Lie-group}); we consider the associated graded groupoid $(\mathrm{B}G,\phi)$ with $\phi: =\varepsilon$. A
 $\phi$-twisted super extension of $\mathrm{B}G$ is the same as a~group extension
$
 1 \longrightarrow \U(1) \longrightarrow \hat{G} \longrightarrow G
 \longrightarrow 1$,
 which for all $\hat g \in \hat{G}$ and $\lambda \in U(1)$ satisfies
 \[
 \hat g\lambda =
 \begin{cases}
 \lambda \hat g, \quad \phi(\hat g) = +1, \\
 \overline{\lambda} \hat g, \quad \phi(\hat g) = -1
 \end{cases}
 \]
 together with a~Lie group homomorphism $\vartheta \colon G \to \{\pm 1\}$.
 The trivial $\phi$-twisted super extension~$I_{\mathrm{B}G,\phi}$ corresponds to the semi-direct product $G \ltimes \U(1)$,
 where $G$ acts on $\U(1)$ via $\phi\colon G \to \Z_2 \cong \Aut(\U(1))$, and the constant Lie group homomorphism $\vartheta = 1$.
 Twisted group extensions are discussed in \cite[Section~3.4]{Braun2002} and \cite[Section~1]{Freeda}. The super version thereof appears in
 \cite[Section~3]{Freeda}.
\end{Example}

Freed--Moore arrange $\phi$-twisted super extensions of $\Gamma$ in a~category under the following notion of morphisms, which we copy from \cite[Definition~7.27]{Freeda}.
We call these morphisms ``refinements'' for reasons that will become clear later.
\begin{Definition} \label{The category of twisted super extensions}
Let $(\Gamma,\phi)$ be a~graded Lie groupoid and let $(\mathcal{L}_1,\lambda_1)$ and $(\mathcal{L}_2,\lambda_2)$ be $\phi$-twisted super extensions of $\Gamma$.
A \emph{refinement} $(\mathcal{L}_1,\lambda_1) \to (\mathcal{L}_2,\lambda_2)$ is an isomorphism
$
 \eta \colon \mathcal{L}_1 \to \mathcal{L}_2
$
 of super line bundles over $\Gamma_1$ such that the diagram
 \begin{equation}
 \begin{aligned} \begin{xy} \label{Diagram: Freed moore refinements}
 \xymatrix{
 \mathcal{L}_1|_{\gamma_1}^{\phi(\gamma_2)} \otimes \mathcal{L}_1|_{\gamma_2} \ar[d]_-{\eta_{\gamma_1}^{\phi(\gamma_2)} \otimes \eta_{\gamma_2}} \ar[rr]^-{\lambda_1|_{\gamma_1,\gamma_2}} && \mathcal{L}_1|_{\gamma_1\circ\gamma_2} \ar[d]^-{\eta_{\gamma_1\circ \gamma_2}} \\
 \mathcal{L}_2|_{\gamma_1}^{\phi(\gamma_2)} \otimes \mathcal{L}_2|_{\gamma_2} \ar[rr]_-{\lambda_2|_{\gamma_1,\gamma_2}} && \mathcal{L}_2|_{\gamma_1\circ \gamma_2}\\
 }
 \end{xy} \end{aligned}
 \end{equation}
commutes for each pair $\gamma_1,\gamma_2\in \Gamma_1$ of composable morphisms.
\end{Definition}

We denote the category of Freed--Moore's $\phi$-twisted super extensions of $\Gamma$ by $\mathrm{FM}(\Gamma,\phi)$.

Interestingly, twisted super extensions also form a~bicategory under a~natural and more general kind of morphisms.
These more general morphisms do not appear in~\cite{Freeda} but will be useful in order to understand the relation to graded-equivariant super 2-line bundles,
and also clarify the role of twisted vector bundles, see Proposition~\ref{Freed-Moore twisted vector bundles}. However, there is an untwisted (i.e., $\phi=1$)
version of these more general morphisms in \cite[Section~2]{Freed2011a}.
The situation is fully analogous to the setting of algebras, which may be regarded as a~category (with algebra homomorphisms), but also as a~bicategory (with bimodules).

\begin{Definition} \label{1-morphisms in the bicategory of twisted super extensions}
Let $(\Gamma,\phi)$ be a~graded Lie groupoid and $(\mathcal{L}_1,\lambda_1)$ and $(\mathcal{L}_2,\lambda_2)$ be $\phi$-twisted super extensions of $\Gamma$.
A \emph{$1$-morphism} $(\mathcal{L}_1,\lambda_1) \to (\mathcal{L}_2,\lambda_2)$ is a~pair $(\mathcal{W},\rho)$ consisting of a~super vector bundle $\mathcal{W}$ over $\Gamma_0$
 and an isomorphism
 \[
 \rho\colon \ \mathcal{L}_2 \otimes s^*\mathcal{W} \to t^*\mathcal{W}^\phi \otimes \mathcal{L}_1
 \]
of super vector bundles over $\Gamma_1$ such that the diagram
 \begin{equation} \label{FM generalized twisted vector bundles category diagram}
 \hspace{-1.1em}\alxydim{@C=-1.8em}{\operatorname{pr}_1^*\mathcal{L}_2^{\operatorname{pr}_2^*\phi} \otimes \operatorname{pr}_2^*\bigl(\mathcal{L}_2 \otimes s^{*}\mathcal{W}\bigr)\ar[r] \ar[d]_{\id \otimes \pr_2^{*}\rho} & \bigl(\operatorname{pr}_1^*\mathcal{L}_2^{\operatorname{pr}_2^*\phi} \otimes \operatorname{pr}_2^*\mathcal{L}_2 \bigr)\otimes \pr_2^{*}s^{*}\mathcal{W} \ar[rr]^-{\lambda_2 \otimes \id} &\hspace{5em}& c^{*}\mathcal{L}_2 \otimes \pr_2^{*}s^{*}\mathcal{W} \ar[d] \\
 \pr_1^{*}\mathcal{L}_2^{\pr_2^{*}\phi} \otimes \pr_2^{*}\bigl(t^{*}\mathcal{W}^{\phi} \otimes \mathcal{L}_1\bigr) \ar[d] &&& c^{*}(\mathcal{L}_2 \otimes s^{*}\mathcal{W}) \ar[dd]^{c^{*}\rho} \\
 \pr_1^{*}\bigl(\mathcal{L}_2\otimes s^{*}\mathcal{W}\bigr)^{\pr_2^{*}\phi} \otimes \pr_2^{*}\mathcal{L}_1 \ar[d]_{\pr_1^{*}\rho^{\pr_2^{*}\phi} \otimes \id} &&& \\
 \pr_1^{*}\bigl(t^{*}\mathcal{W}^{\phi} \otimes \mathcal{L}_1\bigr)^{\pr_2^{*}\phi} \otimes \pr_2^*\mathcal{L}_1 \ar[d] &&& c^{*}\bigl(t^{*}\mathcal{W}^{\phi} \otimes \mathcal{L}_1\bigr) \ar[d]\\ \pr_1^{*}t^{*}\mathcal{W}^{\pr_1^{*}\phi \cdot \pr_2^{*}\phi} \otimes \bigl(\pr_1^{*}\mathcal{L}_1^{\pr_2^{*}\phi} \otimes \pr_2^{*}\mathcal{L}_1\bigr) \ar[rrr]_-{\id \otimes \lambda_1} &&& c^{*}t^{*}\mathcal{W}^{c^{*}\phi} \otimes c^{*}\mathcal{L}_1 }
\hspace{-0.7em}
 \end{equation}
 of vector bundle morphisms over $\Gamma_2$ is commutative.
\end{Definition}

Let $(\mathcal{W},\rho) \colon (\mathcal{L}_1,\lambda_1) \to (\mathcal{L}_2,\lambda_2)$ and $(\mathcal{W}',\rho') \colon (\mathcal{L}_2,\lambda_2) \to (\mathcal{L}_3,\lambda_3)$ be 1-morphisms.
The composition $(\mathcal{W}',\rho') \circ (\mathcal{W},\rho) \colon (\mathcal{L}_1,\lambda_1) \to (\mathcal{L}_3,\lambda_3)$ is the 1-morphism given by
the super vector bundle $\mathcal{W}' \otimes \mathcal{W}$ over $\Gamma_0$ and the isomorphism of vector bundles
 \begin{equation*}
 \alxydim{@C=2em}{
 \mathcal{L}_3 \otimes s^*(\mathcal{W}'\otimes \mathcal{W}) \ar[rr]^-{\rho' \otimes \id} & &
 (t^*\mathcal{W}')^{\phi} \otimes \mathcal{L}_2 \otimes s^*\mathcal{W} \ar[rr]^-{\id \otimes \rho} & &
 t^*(\mathcal{W}' \otimes \mathcal{W})^\phi \otimes \mathcal{L}_1.
 }
 \end{equation*}
Moreover, if $(\mathcal{W},\rho)$ and $(\mathcal{W}',\rho')$ are both 1-morphisms $(\mathcal{L}_1,\lambda_1) \to (\mathcal{L}_2,\lambda_2)$, a~2-morphism $\xi\colon (\mathcal{W},\rho)\Rightarrow (\mathcal{W}',\rho')$
is a~morphism $\xi\colon \mathcal{W} \to \mathcal{W}'$ of super vector bundles over $\Gamma_0$
that is compatible with the isomorphisms ${\rho}$ and ${\rho}'$ in the sense that the following diagram commutes
\begin{equation*}
 \alxydim{}{\mathcal{L}_2 \otimes s^*\mathcal{W} \ar[d]_{\id \otimes s^{*}\xi} \ar[r]^{\rho} & t^*\mathcal{W}^\phi \otimes \mathcal{L}_1 \ar[d]^{t^{*}\xi^{\phi} \otimes \id}\\
 \mathcal{L}_2 \otimes s^*\mathcal{W}' \ar[r]_{\rho'} & t^*\mathcal{W}'^\phi \otimes \mathcal{L}_1.}
\end{equation*}
It is straightforward to complete this structure to a~bicategory, which we denote by $\mathrm{FM}^{bi}(\Gamma,\phi)$, the bicategory of super extensions of the graded Lie groupoid $(\Gamma,\phi)$.

\begin{Proposition}\label{Proposition: Framing of Freed-Moore definitions}
 The bicategory $\mathrm{FM}^{bi}(\Gamma,\phi)$ is framed under Freed--Moore's category of twisted super extensions $\mathrm{FM}(\Gamma,\phi).$
\end{Proposition}
\begin{proof}
 Given a~refinement $\eta \colon (\mathcal{L}_1,\lambda_1) \to (\mathcal{L}_2,\lambda_2)$, we associate to it the 1-morphism $(\mathcal{W},\rho) \colon \allowbreak (\mathcal{L}_1,\lambda_1) \to (\mathcal{L}_2,\lambda_2)$
 whose super vector bundle $\mathcal{W}$ over $\Gamma_0$ is the trivial one and whose map $\rho\colon \mathcal{L}_2 \to \mathcal{L}_1$ is $\eta^{-1}$.
This obviously preserves the composition. Moreover, it is clear that~$(\mathcal{W},\rho)$ is invertible and hence admits an adjoint.
\end{proof}

Now we are in position to establish the relation between the bicategory $\mathrm{FM}^{bi}(\Gamma,\phi)$
of $\phi$-twisted super extensions and our bicategory of twistings $\stwoLineBdl(\Gamma,\phi)$.
\begin{Theorem} \label{The main theorem of the Freed-Moore chapter}
 Let $(\Gamma,\phi)$ be a~graded Lie groupoid. There is a~canonical equivalence of framed bicategories
 \begin{equation} \label{The main theorem of the Freed-Moore chapter diagram}
\alxydim{}{\mathrm{FM}(\Gamma,\phi) \ar[d]_{\cong} \ar[r] & \mathrm{FM}^{bi}(\Gamma,\phi) \ar[d]^{\cong} \\ \stwoLineBdltriv^{\mathrm{ref}}(\Gamma,\phi) \ar[r] & \stwoLineBdltriv(\Gamma,\phi)}
\end{equation}
 between the framed bicategory of $\phi$-twisted super extensions of $\Gamma$ and the framed bicategory
 of super $2$-line bundles over $(\Gamma,\phi)$ whose underlying super $2$-line bundle is trivial.
\end{Theorem}

\begin{proof}
 In order to define the vertical functors on object level, we compare Definitions~\ref{Definition: twisted super extension} and~\ref{equivariant structure}.
 First, we note that trivial super 2-line bundles canonically pull back and complex conjugate to trivial super 2-line bundles.
 Then, the data of the 1-isomorphism $\mathscr{P} \colon \mathscr{I} \to \mathscr{I}$ over $\Gamma_1$ (cf.\ \eqref{equivariant structure: 1-morphism})
 reduces to a~super line bundle $\mathcal{L}$ over $\Gamma_1$
 and the 2-isomorphism $\psi\colon \operatorname{pr}_1^*\mathscr{P}^{\operatorname{pr}_2^*\phi} \circ \operatorname{pr}_2^*\mathscr{P}
 \Rightarrow c^*\mathscr{P}$ (cf.\ \eqref{equivariant structure: 2-morphism}) reduces to an isomorphism
 \smash{$\lambda_{\gamma_1, \gamma_2} \colon \mathcal{L}_{\gamma_1}^{\phi(\gamma_2)} \otimes \mathcal{L}_{\gamma_2} \to \mathcal{L}_{\gamma_1 \circ \gamma_2}$} of super line bundles so that
 \eqref{equivariant structure: diagram} becomes~\eqref{Definition: twisted super extension diagram}. In other words, a~triple
 $(\mathscr{I},\mathscr{P},\psi)$ is the same as a~$\phi$-twisted super extension $(\mathcal{L},\lambda)$ of $\Gamma$. Hence, we define
 the vertical functors in diagram~\eqref{The main theorem of the Freed-Moore chapter diagram} on object level by mapping $(\mathcal{L},\lambda) \to (\mathscr{I},\mathcal{L},\lambda)$.

 A similar comparison can be done on the level of refinements, 1- and 2-morphisms; we will not spell this out in too much detail, but we write down the definitions
 that should be compared. To define the functor $\mathrm{FM}(\Gamma,\phi) \to \stwoLineBdltriv^{\mathrm{ref}}(\Gamma,\phi)$ on the level of morphisms, we note that, for trivial underlying 2-line bundles,
 the data of Definition~\ref{definition: refinement of a~graded equivariant 2-line bundle} reduces to a~2-isomorphism~${\eta \colon \mathscr{P}_2 \Rightarrow \mathscr{P}_1}$
 which, by the description above, is the same as an isomorphism of super line bundles $\mathcal{L}_2 \to \mathcal{L}_1$ over $\Gamma_1$ that is -- due to~\eqref{equivariant 1-morphism: diagram} -- compatible
 with the maps $\lambda$ in the sense of~\eqref{Diagram: Freed moore refinements}.
 In other words, $\eta^{-1} \colon \mathcal{L}_1 \to \mathcal{L}_2$ is precisely a~refinement of $\phi$-twisted super extensions of $\Gamma$
 in the sense of Definition~\ref{The category of twisted super extensions}.

 To define the functor $\mathrm{FM}^{bi}(\Gamma,\phi)\to \stwoLineBdltriv(\Gamma,\phi)$ on the level of 1-morphisms, we note that Definition~\ref{equivariant 1-morphisms}
 reduces -- for trivial underlying 2-line bundles -- to Definition~\ref{1-morphisms in the bicategory of twisted super extensions}. More precisely,
 the 1-morphism $\mathscr{B} \colon \mathscr{I} \to \mathscr{I}$ over~$\Gamma_0$ of~\eqref{equivariant 1-morphism: 1-morphism} is the same as a~super vector bundle~$\mathcal{W}$ over $\Gamma_0$, and the 2-isomorphism $\eta \colon \mathscr{P}_2 \circ s^*\mathscr{B} \Rightarrow t^*\mathscr{B}^\phi \circ \mathscr{P}_1$
 of~\eqref{equivariant 1-morphism: 2-morphism} reduces to an isomorphism~${\rho\colon \mathcal{L}_2 \otimes s^*\mathcal{W} \to t^*\mathcal{W}^\phi \otimes \mathcal{L}_1}$
 of super vector bundles over $\Gamma_1$ such that~\eqref{equivariant 1-morphism: diagram} becomes~\eqref{FM generalized twisted vector bundles category diagram}.
 Similar reasoning holds on the level of 2-morphisms.

 It is also easy to check that both constructions result in functors and moreover in an equivalence of categories and an equivalence of bicategories, respectively.
\end{proof}

\begin{Remark}
Theorem~\ref{The main theorem of the Freed-Moore chapter} together with our cohomological classification Theorem~\ref{classification-of-graded-equivariant-two-line-bundles} shows that $\phi$-twisted super extensions (up to 1-isomorphisms in the sense of Definition~\ref{1-morphisms in the bicategory of twisted super extensions}) are classified by the kernel of the homomorphism
\begin{equation*}
p_0 \times p_0\colon\ \check\h^1(\Gamma,\Z_2) \times \check\h^2((\Gamma,\phi),\sheaf{\C^{\times}}) \to \h^1(\Gamma_0,\Z_2) \times \check\h^2(\Gamma_0,\sheaf{\C^{\times}}),
\end{equation*}
since this homomorphism projects to the underlying super 2-line bundle over $\Gamma_0$,
which is trivial in Theorem~\ref{The main theorem of the Freed-Moore chapter}. In particular, in the case of an action groupoid $(\act XG,\phi_{\varepsilon})$ for a~graded Lie group $(G,\varepsilon)$ acting on a~smooth manifold $X$, we obtain
\begin{gather*}
\hc 0 \mathrm{FM}^{bi}(\act XG,\phi_{\varepsilon})\\
\qquad=\mathrm{ker}\bigl(p_0 \times p_0\colon \check\h^1(\act XG,\Z_2) \times \check\h^2\bigl((\act XG,\phi_{\varepsilon}),\sheaf{\C^{\times}}\bigr) \to \h^1(X,\Z_2) \times \check\h^2\bigl(X,\sheaf{\C^{\times}}\bigr)\bigr);
\end{gather*}
which is in line with a~rigorous combination of the ideas sketched in \cite[Remark~7.28]{Freeda} and the discussion of the ungraded ($\phi=1$) case in
\cite[Section~2.2.1]{Freed2011a}.
Even more special, if
 $\Gamma=\mathrm{B}G$ (see Example~\ref{twisted-groupoid-extension-of-BG}), we get
\begin{equation*}
\hc 0 \mathrm{FM}^{bi}(\mathrm{B}G,\phi) \cong \check\h^1(\mathrm{B}G,\Z_2) \times \check\h^2\bigl((\mathrm{B}G,\phi),\sheaf{\C^{\times}}\bigr).
\end{equation*}
\end{Remark}

In the remainder of this section, we
embed the twisted super vector bundles of Freed--Moore~\cite[Definition~7.23]{Freeda} into our framework.

\begin{Definition}
 Let $(\Gamma,\phi)$ be a~graded Lie groupoid and let $(\mathcal{L},\lambda)$ a~$\phi$-twisted super extension of $\Gamma$.
 A \emph{$(\mathcal{L},\lambda)$-twisted super vector bundle} over $\Gamma$ is a~super vector bundle $\mathcal{W}$ over $\Gamma_0$
 together with an isomorphism of super vector bundles
\begin{equation} \label{structural map of a~FM vector bundle}
\rho\colon\ (\mathcal{L} \otimes s^{*}\mathcal{W})^{\phi} \to t^{*}\mathcal{W}
\end{equation}
over $\Gamma_1$ such that the following diagram
\begin{equation} \label{FM vector bundles diagram}
\begin{gathered}
 \begin{xy}
 \xymatrix@C=0em{
 \bigl(\bigl(\mathcal{L}_{\gamma_1}^{\phi(\gamma_2)} \otimes \mathcal{L}_{\gamma_2 }\bigr) \otimes \mathcal{W}_{s(\gamma_1 \circ \gamma_2)}\bigr)^{\phi(\gamma_1 \circ \gamma_2)} \ar@{=}[rr] \ar[d]_-{(\lambda_{\gamma_1,\gamma_2}\otimes \id)^{\phi(\gamma_1 \circ \gamma_2)}}
 && \bigl(\mathcal{L}_{\gamma_1} \otimes (\mathcal{L}_{\gamma_2 } \otimes \mathcal{W}_{s(\gamma_2)})^{\phi( \gamma_2)}\bigr)^{\phi(\gamma_1)} \ar[d]^-{(\id \otimes \rho_{\gamma_2})^{\phi(\gamma_1)}}\\
 (\mathcal{L}_{\gamma_1\circ \gamma_2} \otimes \mathcal{W}_{s(\gamma_1 \circ \gamma_2)})^{\phi(\gamma_1 \circ \gamma_2)} \ar[dr]_{\rho_{\gamma_1\circ \gamma_2}} &&
 (\mathcal{L}_{\gamma_1} \otimes \mathcal{W}_{s(\gamma_1)})^{\phi(\gamma_1)} \ar[dl]^{\rho_{\gamma_1}} \\
 & \mathcal{W}_{t(\gamma_1 \circ \gamma_2)}
 }
 \end{xy}
\end{gathered}
\end{equation}
of vector bundle morphisms commutes for each pair $(\gamma_1,\gamma_2) \in \Gamma_2$.
\end{Definition}

\begin{Remark} \label{FM vector bundles remark}
Comparing this with \cite[Definition~7.23\,(iv)]{Freeda}, we see that~\eqref{structural map of a~FM vector bundle} coincides with \cite[equation~(7.24)]{Freeda} and diagram~\eqref{FM vector bundles diagram} is what is demanded as an analog of diagram \cite[equation~(7.16)]{Freeda}.
The requirement that ``the map $\rho$ is even or odd according to $\vartheta(\gamma)$'' is implemented here automatically as $\rho$ is a~super vector bundle homomorphism, and $\vartheta$ is the grading of $\mathcal{L}$ in our setting.
There are two minor differences to~\cite{Freeda}:
 \begin{itemize}\itemsep=0pt
 \item The map $\rho$ in~\cite{Freeda} is only required to be linear,
 not an isomorphism; however, since we are talking about equivariance, after all, our choice appears quite naturally.
 \item We have formulated all definitions throughout the paper in a~smooth setting, whereas the regularity is partly not clearly stated in loc.cit.
 \end{itemize}

\end{Remark}

The following definition establishes the natural kind of morphism between $(\mathcal{L},\lambda)$-twisted vector bundles; this notion does not appear in~\cite{Freeda}.
\begin{Definition} \label{The category of twisted super vector bundles}
 Let $(\Gamma,\phi)$ be a~graded Lie groupoid, let $(\mathcal{L},\lambda)$ be a~$\phi$-twisted super extension of $\Gamma$, and let
 $(\mathcal{W},\rho)$ and $(\mathcal{W}',\rho')$ be $(\mathcal{L},\lambda)$-twisted super vector bundles. A morphism $(\mathcal{W},\rho) \to (\mathcal{W}',\rho')$ is a~morphism of super vector bundles $\xi\colon \mathcal{W} \to \mathcal{W}'$ over $\Gamma_0$
 that is compatible with $\rho$ and $\rho'$ in the sense that the following diagram is commutative
 \[
 \begin{aligned} \begin{xy}
 \xymatrix{
 (\mathcal{L} \otimes s^{*}\mathcal{W})^{\phi} \ar[rr]^-{\rho} \ar[d]_-{(\id \otimes s^*\xi)^{\phi}} && t^{*}\mathcal{W} \ar[d]^-{t^*\xi}\\
 (\mathcal{L} \otimes s^{*}\mathcal{W}')^{\phi} \ar[rr]_-{\rho'} && t^{*}\mathcal{W}'.
 }
 \end{xy} \end{aligned}
 \]
We denote the category of $(\mathcal{L},\lambda)$-twisted super vector bundles by $\mathrm{FM}\text{-}\sVectBdl^{(\mathcal{L},\lambda)}(\Gamma,\phi)$.
\end{Definition}

\begin{Proposition} \label{Freed-Moore twisted vector bundles}
 There is an equivalence of categories
 \[
 \mathrm{FM}\text{-}\sVectBdl^{(\mathcal{L},\lambda)}(\Gamma,\phi) \cong \Hom_{\mathrm{FM}^{bi}(\Gamma,\phi)}((\mathcal{L},\lambda),I_{\Gamma,\phi}).
 \]
\end{Proposition}
\begin{proof}
 This is clear by comparing definitions and slight rephrasing.
 The isomorphism~\eqref{structural map of a~FM vector bundle} can equivalently be understood
 as an isomorphism
$
 \mathcal{L} \otimes s^{*}\mathcal{W} \to (t^{*}\mathcal{W})^{\phi}
$
 which by duality is then equivalent to an isomorphism
$
 s^{*}\mathcal{W}^* \to (t^{*}\mathcal{W}^*)^{\phi} \otimes \mathcal{L}$.
 Hence, we assign $\bigl(\mathcal{W}^*,(\rho^*)^{-1}\bigr) \in \Hom_{\mathrm{FM}^{bi}(\Gamma,\phi)}((\mathcal{L},\lambda),I_{\Gamma,\phi})$
 to the object \smash{$(\mathcal{W},\rho) \in \mathrm{FM}\text{-}\sVectBdl^{(\mathcal{L},\lambda)}(\Gamma,\phi)$}. It is straightforward to compare the respective compatibility diagrams and
 extend this assignment into an equivalence of categories.
\end{proof}

\begin{Corollary}
\label{twisted-super-vector-bundles-FM}
There is an equivalence of categories
\begin{equation*}
\mathrm{FM}\text{-}\sVectBdl^{(\mathcal{L},\lambda)}(\Gamma,\phi) \cong \sVectBdl^{(\mathscr{I},\mathcal{L},\lambda)}(\Gamma,\phi) .
\end{equation*}
Thus, the twisted vector bundles of Freed--Moore are precisely our $(\mathscr{I},\mathcal{L},\lambda)$-twisted super vector bundles in the sense of Definition~{\rm\ref{Definition of our twisted vector bundles}}.
\end{Corollary}

\begin{proof}
We combine Proposition~\ref{Freed-Moore twisted vector bundles} and Theorem~\ref{The main theorem of the Freed-Moore chapter}, using the fact that $\stwoLineBdltriv(\Gamma,\phi)$ is a~full sub-bicategory of $\stwoLineBdl(\Gamma,\phi)$.
\end{proof}

\subsection{Freed--Hopkins--Teleman's twistings}
\label{Freed--Hopkins--Teleman-twistings}

We show in
Theorem~\ref{The main theorem of the Freed Hopkins Teleman chapter} that the twistings of Freed, Hopkins, and Teleman are equivariant super bundle gerbes, and hence, in particular, equivariant super 2-line bundles.
To start with, the twistings of~\cite{Freed2011a} are based on the notion of a~central super extension of a~Lie groupoid.

\begin{Definition} \label{Definition: central super extension}
Let $\Gamma$ be a~Lie groupoid. The bicategory of \emph{central super extensions of $\Gamma$} is defined as
$
\sExt(\Gamma) := \mathrm{FM}^{bi}(\Gamma,1)$,
i.e., a~central super extension of $\Gamma$ is a~$\phi$-twisted super extension of $\Gamma$ in the sense of
 Definition~\ref{Definition: twisted super extension}, for the trivial grading $\phi=1$.
\end{Definition}

This definition has been given originally in \cite[Section~2]{Freed2011a}; here we present it -- for economical reasons -- as a~special case of our bicategory extension of Freed--Moore's category of $\phi$-twisted super extensions, carried out in Section~\ref{freed-moore-extensions}.

\begin{Remark}
\label{sExtref}
Likewise, as the bicategory extension $\mathrm{FM}^{bi}(\Gamma,\phi)$ is framed under Freed--Moore's original category $\mathrm{FM}(\Gamma,\phi)$, the bicategory $\sExt(\Gamma)$ is framed under the category $\sExt^\mathrm{ref}(\Gamma) := \mathrm{FM}(\Gamma,1)$.

\end{Remark}

The following definition is a~smooth version of Freed--Hopkins--Teleman's definition of twistings, see \cite[Section~2.3]{Freed2011a}.

\begin{Definition}
\label{Freed--Hopkins--Teleman's category of twistings}
 Let $\Gamma$ be a~Lie groupoid. An \emph{FHT-twisting} on $\Gamma$ is a~quadruple $(\Lambda, T,\mathcal{L},\lambda)$ consisting of a~Lie groupoid $\Lambda$, a~weak equivalence $T \colon \Lambda \to \Gamma$ and a~central super extension $(\mathcal{L},\lambda)$ of $\Lambda$.
 A morphism $(\Lambda_1,T_1,\mathcal{L}_1,\lambda_1) \to (\Lambda_2,T_2,\mathcal{L}_2,\lambda_2)$ is an equivalence class $[\Omega, S,\mathcal{W},\rho]$ represented by a~Lie groupoid $\Omega$, a~weak equivalence $S\colon \Omega \to \Lambda_1 \times_{\Gamma} \Lambda_2$ and of
 a~1-morphism
 \begin{equation*}
 (\mathcal{W},\rho) \colon\ S^{*}\pr_1^{*}(\mathcal{L}_1,\lambda_1) \to S^{*}\pr_2^{*}(\mathcal{L}_2,\lambda_2)
 \end{equation*}
 in $\sExt(\Omega)$.
 The equivalence relation identifies $(\Omega,S,\mathcal{W},\rho)$ with $(\Omega',S',\mathcal{W}',\rho')$ if there exists a~smooth functor $G\colon \Omega \to \Omega'$,
 a~smooth natural transformation $S \cong S' \circ G$, and a~2-isomorphism
 \begin{equation*}
\alxydim{@C=5em}{S^{*}\pr_1^{*}(\mathcal{L}_1,\lambda_1) \ar[d] \ar[r]^{(\mathcal{W},\rho)} & S^{*}\pr_2^{*}(\mathcal{L}_2,\lambda_2) \ar@{=>}[dl] \ar[d] \\ G^{*}S'^{*}\pr_1^{*}(\mathcal{L}_1,\lambda_1) \ar[r]_{G^{*}(\mathcal{W}',\rho')}& G^{*}S'^{*}\pr_2^{*}(\mathcal{L}_2,\lambda_2) }
\end{equation*}
in $\sExt(\Omega)$, where the vertical 1-morphisms are the ones induced by the natural isomorphism $S \cong S' \circ G$ under a~procedure analogous to our Lemma~\ref{pullback-invariance-under-nat-iso}.
 We denote the category of FHT-twistings on the Lie groupoid $\Gamma$ by $\mathrm{FHT}(\Gamma)$.
\end{Definition}

We recall the fibre product of Lie groupoids used above in Definition~\ref{definition: fibre product of Lie groupoids}.

The composition of two morphisms $[\Omega,S,\mathcal{W},\rho] \colon (\Lambda_1,T_1,\mathcal{L}_1,\lambda_1) \to (\Lambda_2,T_2,\mathcal{L}_2,\lambda_2)$ and
$\big[\tilde\Omega,\tilde{S},\allowbreak\smash{\tilde{\mathcal{W}},\tilde{\rho}}\big] \colon (\Lambda_2, T_2,\mathcal{L}_2,\lambda_2) \to (\Lambda_3,T_3,\mathcal{L}_3,\lambda_3)$ in $\mathrm{FHT}(\Gamma)$ consists of the fibre product Lie groupoid $\Omega \times_{\Lambda_2} \tilde{\Omega}$, the weak equivalence
\begin{equation*}
 \xymatrix{
 \Omega \times_{\Lambda_2} \tilde{\Omega} \ar[rr]^-{S \times_{\Lambda_2} \tilde{S}} && (\Lambda_1 \times_{\Gamma} \Lambda_2) \times_{\Lambda_2} (\Lambda_2 \times_{\Gamma} \Lambda_3)
 \ar[rr]^-{\pr_{14}} && \Lambda_1 \times_{\Gamma} \Lambda_3,
 }
\end{equation*}
and of the 1-morphism given by the composition of the pullbacks of $(\mathcal{W},\rho)$ and $\bigl(\tilde{\mathcal{W}},\tilde{\rho}\bigr)$ to $\Omega \times_{\Lambda_2} \tilde{\Omega}$.
The identity 1-morphism of an object $(\Lambda,T,\mathcal{L},\lambda) \in \mathrm{FHT}(\Gamma)$ can be represented by the diagonal functor $\Delta \colon \Lambda \to \Lambda \times_{\Gamma} \Lambda$
together with the trivial 1-isomorphism $\Delta^*\pr_1^{*}(\mathcal{L},\lambda) \cong \Delta^*\pr_2^{*}(\mathcal{L},\lambda)$.

\begin{Remark}
Our Definition~\ref{Freed--Hopkins--Teleman's category of twistings} is a~smooth version of the original definition of~\cite{Freed2011a}.
More precisely, in that original version, the groupoids $\Lambda$ and $\Omega$ are \emph{topological} groupoids, and the central super extension $(\mathcal{L},\lambda)$ is a~\emph{continuous} version of smooth central super extensions. Moreover, the functors $T$ and $S$ are \emph{local} equivalences instead of weak equivalences.
 A continuous functor~${F \colon \Gamma \to \Lambda}$ between topological groupoids is called a~\emph{local equivalence} \cite[Definition A.4]{Freed2011a} if it is
 topologically essentially surjective (this is analogous to the first point of Definition~\ref{Definition: weak equivalence of Lie groupoids}, with ``surjective submersion'' replaced by ``locally split map'') and if the induced map on hom-sets
 \begin{equation}
 \label{homeomorphism-on-hom-sets}
 F_{x,x'} \colon\ \hom_{\Gamma}\bigl(x,x'\bigr) \to \hom_{\Lambda}\bigl(F_0(x),F_0\bigl(x'\bigr)\bigr)
 \end{equation}
 is a~homeomorphism for all objects $x,x' \in \Gamma_0$, where $\hom_{\Gamma}(x,x')\subset \Gamma_1$ and analogously $\hom_{\Lambda}(F_0(x),F_0(x'))$ carry the subspace topologies. One can show that weak equivalences are always local equivalences, but the converse may not hold.
Indeed, it is claimed in \cite[Remark~A.5 and Example~A.7]{Freed2011a} that a~continuous functor is a~\emph{local} equivalence if and only if it induces equivalences on stalks between the topological stacks represented by the groupoids; whereas a~functor is a~\emph{weak} equivalence if and only if it induces an equivalence between the stacks.
Denoting the category of original twistings, as
defined in
 \cite{Freed2011a}, by $\mathrm{FHT}^{\rm top}(\Gamma)$; our Definition~\ref{Freed--Hopkins--Teleman's category of twistings} yields a~subcategory
\begin{equation}
\label{weak-FHT-twists-in-all-FHT-twists}
\mathrm{FHT}(\Gamma) \subset \mathrm{FHT}^{\rm top}(\Gamma).
\end{equation}
This subcategory inclusion~\eqref{weak-FHT-twists-in-all-FHT-twists} is essentially surjective. This follows from the classification results: \cite[Corollary~2.25]{Freed2011a} shows a~bijection
\begin{equation*}
\hc 0 \mathrm{FHT}^{\rm top}(\Gamma) \cong \check \h^1(\Gamma,\Z_2) \times \check\h^2(\Gamma,\sheaf{\C^{\times}}).
\end{equation*}
On the other hand, Theorem~\ref{The main theorem of the Freed Hopkins Teleman chapter} and Proposition~\ref{classification-of-equivariant-two-line-bundles}
imply that $\mathrm{FHT}(\Gamma)$ has the same classification. Thus, for a~Lie groupoid $\Gamma$, our smooth version of FHT-twistings coincides with the original ones.
\end{Remark}

\begin{Theorem}
\label{The main theorem of the Freed Hopkins Teleman chapter}
 Let $\Gamma$ be a~Lie groupoid. There is an equivalence of categories
$
\mathrm{FHT}(\Gamma) \cong \hc 1(\sGrb(\Gamma))\subset \hc 1(\stwoLineBdl(\Gamma))
$
between the category of Freed--Hopkins--Teleman twistings and the homotopy 1-category of super bundle gerbes over $\Gamma$.
\end{Theorem}

\begin{proof}
The equivalence is split into two functors
\begin{equation*}
\hc 1(\sGrb(\Gamma)) \to\hc 1\bigl(\sGrbtriv^{+}(\Gamma)\bigr) \stackrel{\mathcal{F}}\to\mathrm{FHT}(\Gamma)
\end{equation*}
which are both equivalences. The category in the middle is the homotopy category of the bicategory $\sGrbtriv^{+}(\Gamma)$, which is obtained by considering the presheaf of bicategories $\sGrbtriv$ of trivial super bundle gerbes, stackifying it, and evaluating it on the Lie groupoid $\Gamma$. Below we describe it explicitly.
By definition of a~bundle gerbe (see~\cite{nikolaus2}), there is an equivalence~${
\sGrb \cong \sGrbtriv^{+}}
$
between 2-stacks on the category of smooth manifolds, which extends to an equivalence of their extensions to Lie groupoids. This shows the first equivalence. The functor~$\mathcal{F}$ is constructed below; the proof that it is an equivalence is given in Proposition~\ref{The FHT functor: essential surjectivity}.
\end{proof}

The construction of the functor $\mathcal{F}$ faces two main challenges:
\begin{enumerate}[(1)]\itemsep=0pt
\item
The translation between the language of bundle gerbes and the language of central super extensions.

\item
The transformation of surjective submersions
and their fibre products into weak equivalences and their fibre products.
\end{enumerate}
The first challenge is taken by the reduction of Theorem~\ref{The main theorem of the Freed-Moore chapter} to ungraded Lie groupoids, which we state (using Remark~\ref{trivial-things}) as the following corollary.

\begin{Corollary}
\label{central-super-extensions-and-line-2-bundles}
The bicategories of
central super extensions and of trivial super bundle gerbes over a~Lie groupoid $\Gamma$ are canonically equivalent,
$
\sExt(\Gamma) \cong \sGrbtriv(\Gamma)$.
\end{Corollary}

\begin{Remark}
The equivalence of Corollary~\ref{central-super-extensions-and-line-2-bundles} was obtained separately by Mertsch
\cite[Proposition~B.3.7]{Mertsch2020}.
In \cite[Theorem~B.3.13]{Mertsch2020}, a~further comparison between the category of FHT-twistings and super bundle gerbes
was obtained for global quotient groupoids $\Gamma=\act XG$. The main insight of Mertsch was that the definition of FHT-twistings
resembles the plus construction described in {\cite[Proposition 9.1]{nikolaus2}}. The idea of our construction is similar; the main difference is that
we do not restrict to global quotient groupoids, and use weak equivalences instead of local equivalences. According to \cite[Lemma B.3.10]{Mertsch2020}, the notions of local and weak equivalences coincide for certain global quotient groupoids. \end{Remark}

The second challenge is taken manually in course of the construction of the functor $\mathcal{F}$.
We first recall a~concrete description of the bicategory $\sGrbtriv^{+}(\Gamma)$ obtained in \cite[Proposition 9.1]{nikolaus2}.
We use the notion of coverings groupoids $\pi\colon \Gamma^{\pi}\to \Gamma$ from Definition~\ref{Example : covering groupoid and weak equivalences}, for a~Lie groupoid~$\Gamma$ and a~surjective submersion $\pi\colon Y \to \Gamma_0$. We also use the terminology that a~\emph{refinement} of a~surjective submersion $\pi\colon Y \to M$ is another surjective submersion $\zeta\colon Z \to M$ together with a~surjective submersion $Z \to Y$ such that the diagram
\begin{equation*}
\alxydim{}{Z \ar[dr]_{\zeta} \ar[rr] && Y \ar[dl]^{\pi} \\ &M}
\end{equation*}
is commutative. In that notation, the description of $\sGrbtriv^{+}(\Gamma)$ is the following:
\begin{enumerate}[(i)]\itemsep=0pt
\item

 An object is a~quadruple $(Y,\pi,\mathscr{P},\psi)$, where $\pi \colon Y \to \Gamma_0$ is a~surjective submersion and $(\mathscr{I},\mathscr{P},\psi)$ is an object in $\sGrbtriv(\Gamma^{\pi})$ (in the triple notation of Definition~\ref{equivariant structure}).

 \item A 1-morphism $(Y,\pi,\mathscr{P},\psi) \to \bigl(Y',\pi',\mathscr{P}',\psi'\bigr)$ is again a~quadruple $(Z,\zeta,\mathscr{B},\eta)$, where
 $\zeta\colon Z \to \Gamma_0$ is a~refinement of $Y \times_{\Gamma_0} Y' \to \Gamma_0$ and
\begin{equation*}
 (\mathscr{B},\eta) \colon\ \zeta_{Y}^{*}(\mathscr{I},\mathscr{P},\psi) \to \zeta_{Y'}^{*}(\mathscr{I},\mathscr{P}',\psi')
\end{equation*}
is a~1-morphism in $\sGrbtriv\bigl(\Gamma^{\zeta}\bigr)$ (in the notation of Definition~\ref{equivariant 1-morphisms}). Here, we have denoted the projections from $Z$ to $Y$ and $Y'$ by $\zeta_Y$ and $\zeta_{Y'}$, respectively; the passage from refinements to the corresponding weak equivalences can be visualized as follows:
\begin{equation*}
\alxydim{@C=3em@R=3em@!0}{& Z \ar@/_1pc/[ddl]_{\zeta_Y} \ar@/^1pc/[ddr]^{\zeta_{Y'}} \ar[d] &&&& \Gamma^{\zeta} \ar@/_1pc/[ddl]_{\zeta_Y} \ar@/^1pc/[ddr]^{\zeta_{Y'}} \\ & Y \times_{\Gamma_0} Y' \ar[dr]\ar[dl]\\Y \ar[dr]_{\pi} && Y' \ar[dl]^{\pi'} & \rightsquigarrow & \Gamma^{\pi} && \Gamma^{\pi'} \\ & \Gamma_0}
\end{equation*}

\item
\label{explanation-of-composition}
The composition
\begin{equation*}
\alxydim{@C=6em}{(Y,\pi,\mathscr{P},\psi) \ar[r]^-{(Z,\zeta,\mathscr{B},\eta)} & \bigl(Y',\pi',\mathscr{P}',\psi'\bigr)\ar[r]^-{(Z',\zeta',\mathscr{B}',\eta')} & (Y'',\pi'',\mathscr{P}'',\psi'')}
\end{equation*}
of 1-morphisms is defined over the fibre product $\tilde Z := Z \times_{Y'} Z'$, considered as a~refinement of $Y \times_{\Gamma_0} Y'' \to \Gamma_0$, together with the 1-morphism $\big(\tilde\zeta_{Z'}\big)^{*}(\mathscr{B}',\eta') \circ \big(\tilde\zeta_Z\big)^{*}(\mathscr{B},\eta)$.
The involved refinements are depicted in the following diagrams:
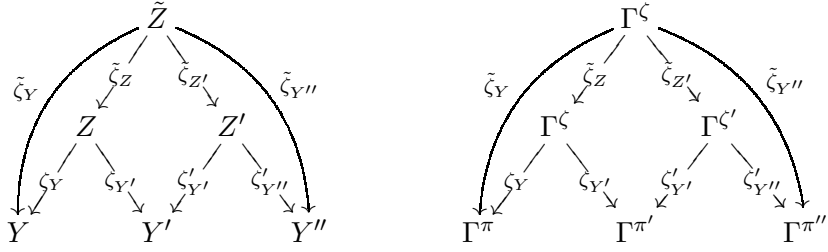
\begin{equation*}
\alxydim{@C=1em}{
&& \tilde Z \ar@/_1.5pc/[ddll]_{\tilde\zeta_Y} \ar@/^1.5pc/[ddrr]^{\tilde\zeta_{Y''}} \ar[dr]|{\tilde\zeta_{Z'}}\ar[dl]|{\tilde\zeta_{Z}} &&&&&&& \Gamma^{\tilde\zeta} \ar@/_1.5pc/[ddll]_{\tilde\zeta_Y} \ar@/^1.5pc/[ddrr]^{\tilde\zeta_{Y''}} \ar[dr]|{\tilde\zeta_{Z'}}\ar[dl]|{\tilde\zeta_{Z}}
\\
& Z \ar[dl]|{\zeta_Y}\ar[dr]|{\zeta_{Y'}} && Z' \ar[dl]|{\zeta'_{Y'}}\ar[dr]|{\zeta'_{Y''}} &&&&& \Gamma^{\zeta} \ar[dl]|{\zeta_Y} \ar[dr]|{\zeta_{Y'}} && \Gamma^{\zeta'} \ar[dl]|{\zeta'_{Y'}} \ar[dr]|{\zeta'_{Y''}}
\\
Y && Y' && Y'' &&& \Gamma^{\pi} && \Gamma^{\pi'} && \Gamma^{\pi''}
}
\end{equation*}

\item
\label{identity-morphisms-explained}
The identity morphism of an object $(Y,\pi,\mathscr{P},\psi)$ is defined over the fibre product $Z:= Y \times_{\Gamma_0} Y$, considered as a~refinement $\zeta\colon Z \to \Gamma_0$ of itself; its 1-morphism over \smash{$\Gamma^{\zeta}_0=Z$} is~${\mathscr{B} := i^{*}\mathscr{P}}$, where $i\colon Z \to \Gamma^{\pi}_1$ is defined by $(y_1,y_2) \mapsto (y_2,\id,y_1)$, and its 2-isomorphism $\eta$ over \smash{$\Gamma^{\zeta}_1=Y \times_{\Gamma_0} Y \times_{\Gamma_0}\Gamma_1 \times_{\Gamma_0}Y \times_{\Gamma_0} Y$} is
\smash{$
\eta := i_{2}^{*}\psi^{-1} \circ i_{1}^{*}\psi
$}
where the two maps $\smash{\Gamma^{\zeta}_1} \to \Gamma^{\pi}_2$ are given by
\begin{align*}
& i_{1}(y_1,y_2,\gamma,y_4,y_5) := ((y_2,\gamma,y_5),(y_5,\id,y_4)) \qquad \text{and} \\
& i_{2}(y_1,y_2,\gamma,y_4,y_5) := ((y_2,\id,y_1),(y_1,\gamma,y_4)).
\end{align*}

\item
\label{2-morphisms-in-grb-triv}
Let \smash{$(Z,\zeta,\mathscr{B},\eta),\bigl(\tilde Z,\tilde{\zeta},\tilde{\mathscr{B}},\tilde\eta\bigr) \colon (Y,\pi,\mathscr{P},\psi) \to \bigl(Y',\pi',\mathscr{P}',\psi'\bigr)$} be 1-morphisms.
 A 2-morphism \smash{$(Z,\zeta,\mathscr{B},\eta) \Rightarrow \bigl(\tilde Z,\tilde{\zeta},\tilde{\mathscr{B}},\tilde\eta\bigr)$} is an equivalence class of triples $(W,\omega,\xi)$ consisting of
 a~refinement $\omega\colon W \to \Gamma_0$ of \smash{$Z \times_{Y \times_{\Gamma_0} Y'} \tilde Z\to \Gamma_0$} and of a~2-morphism
$
 \xi\colon \omega_Z^{*}(\mathscr{B},\eta) \Rightarrow \omega_{\tilde Z}^{*}\bigl(\tilde {\mathscr{B}},\tilde\eta\bigr)
$
in $\sGrbtriv(\Gamma^{\omega})$;
 the various refinements and corresponding weak equivalences are depicted in the following commutative diagram:
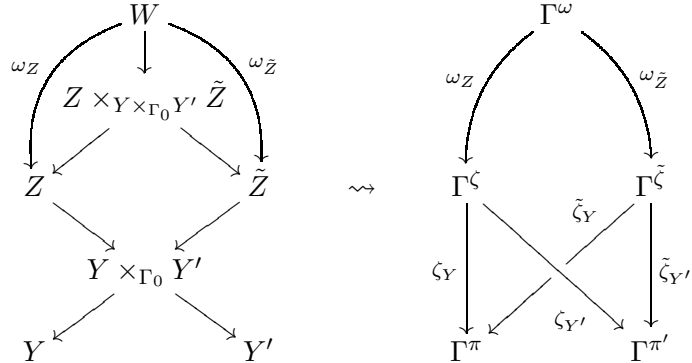
\begin{equation*}
\alxydim{@C=0em@R=1.3em}{ & W \ar@/^1.5pc/[ddr]^{\omega_{\tilde Z}} \ar@/_1.5pc/[ddl]_{\omega_{Z}} \ar[d] &&&&\Gamma^{\omega} \ar@/_1pc/[ddl]_{\omega_Z}\ar@/^1pc/[ddr]^{\omega_{\tilde Z}} \\ & Z \times_{Y \times_{\Gamma_0} Y'} \tilde Z\ar[dr] \ar[dl] \\ Z \ar[dr] && \tilde Z \ar[dl]& \qquad\rightsquigarrow\qquad & \Gamma^{\zeta} \ar[ddrr]_>>>>>>{\zeta_{Y'}} \ar[dd]_{\zeta_Y} & \hspace{4em}& \Gamma^{\tilde\zeta} \ar[ddll]_<<<<<<{\tilde\zeta_{Y}}|!{[ll];[dd]}\hole \ar[dd]^{\tilde\zeta_{Y'}}\\ & Y \times_{\Gamma_0} Y' \ar[dr]\ar[dl]\\Y && Y' & & \Gamma^{\pi} && \Gamma^{\pi'} }
\end{equation*}
The equivalence relation identifies $(W,\omega,\xi) \sim \bigl(\hat W,\hat{\omega},\hat{\xi}\bigr)$ if there is a~refinement of the fibre product of $W$ and $\hat W$ over \smash{$Z \times_{Y \times_{\Gamma_0} Y'} \tilde Z$}, on which the pullbacks of $\xi$ and \smash{$\hat{\xi}$} coincide. Since we are only interested in the homotopy 1-category, the details of this equivalence relation, as well as horizontal and vertical composition of 2-morphisms, are not relevant.
\end{enumerate}

\begin{Remark}
Under this description of the bicategory $\sGrbtriv^{+}(\Gamma)$, the functor $\sGrb(\Gamma) \to \sGrbtriv^{+}(\Gamma)$ is the following. If $(\mathscr{G},\mathscr{P},\psi)$ is a~super bundle gerbe over $\Gamma$, where $\mathscr{G}=(Y,\pi,\mathcal{L},\mu)$ is its super bundle gerbe over $\Gamma_0$,
then $\pi^{*}(\mathscr{G},\mathscr{P},\psi)$ is a~super bundle gerbe over $\Gamma^{\pi}$, containing the super bundle gerbe $\pi^{*}\mathscr{G}$ over $\Gamma^{\pi}_0=Y$. The surjective submersion of that bundle gerbe has a~canonical section, which induces a~canonical trivialization $\mathscr{T}\colon \pi^{*}\mathscr{G} \to \mathscr{I}$. This, in turn, induces a~1-isomorphism $\pi^{*}(\mathscr{G},\mathscr{P},\psi) \cong (\mathscr{I},\mathscr{P}',\psi')$ in $\sGrb(\Gamma^{\pi})$, whose codomain is an object in~$\sGrbtriv(\Gamma^{\pi})$. Thus, the quadruple $(Y,\pi,\mathscr{P}',\psi')$ is an object in $\sGrbtriv^{+}(\Gamma)$.
\end{Remark}

Now we are in position to construct the functor
$
\mathcal{F}\colon\ \hc 1\bigl(\sGrbtriv^{+}(\Gamma)\bigr) \to \mathrm{FHT}(\Gamma)$.
 To an object in $(Y,\pi,\mathscr{P},\psi) \in \hc 1\bigl(\sGrbtriv^{+}(\Gamma)\bigr)$ we assign the object $(\Gamma^\pi,\pi,\mathcal{L},\lambda)\in\mathrm{FHT}(\Gamma)$ consisting of
 the weak equivalence $\pi\colon \Gamma^\pi \to \Gamma$ and the
 central super extension $(\mathcal{L},\lambda)$ of $\Gamma^{\pi}$ that corresponds to the super bundle gerbe $(\mathscr{I},\mathscr{P},\psi)$ over $\Gamma^{\pi}$ under the equivalence of Corollary~\ref{central-super-extensions-and-line-2-bundles}.
 To a~morphism
 \begin{equation*}
 [(Z,\zeta,\mathscr{B},\eta)] \colon\ (Y,\pi,\mathscr{P},\psi) \to \bigl(Y',\pi',\mathscr{P}',\psi'\bigr)
 \end{equation*}
 in $\hc 1\bigl(\sGrbtriv^{+}(\Gamma)\bigr)$ we assign the equivalence class
 $\big[\Gamma^{\zeta},S,\mathcal{W},\rho\big]$, where $S\colon \Gamma^{\zeta} \to \Gamma^\pi \times_{\Gamma} \Gamma^{\pi'}$ is the weak equivalence obtained from the refinement map $Z \to Y \times_{\Gamma_0} Y'$ and the weak equivalence of Lemma~\ref{Common refinements of common refinements induce weak equivalences NEU}, and $(\mathcal{W},\rho)$ is the 1-morphism between central super extension that corresponds to $(\mathscr{B},\eta)$ under the equivalence of Corollary~\ref{central-super-extensions-and-line-2-bundles}.

 First, we need to check that this is well-defined. Indeed, let $(Z,\zeta,\mathscr{B},\eta)$ and $\bigl(\tilde Z,\tilde \zeta,\tilde{\mathscr{B}},\tilde \eta\bigr)$ be two representatives of the same 1-morphism $\hc 1\bigl(\sGrbtriv^{+}(\Gamma)\bigr)$, i.e., there exists a~2-isomorphism $(W,\omega,\xi)$ between them. The refinement diagram in~\ref{2-morphisms-in-grb-triv} above
yields, together with Lemma~\ref{Common refinements of common refinements induce weak equivalences NEU}, a~strictly commutative diagram of weak equivalences
\begin{equation*}
\alxydim{}{& \Gamma^{\omega}\ar[dd]^{U} \ar[dl]_{F} \ar[dr]^{\tilde F} \\ \Gamma^{\zeta} \ar[dr]_-{S} && \Gamma^{\tilde\zeta} \ar[dl]^-{\tilde S} \\ & \Gamma^{\pi} \times_{\Gamma} \Gamma^{\pi'} .}
\end{equation*}
Moreover, $\xi$ is a~2-isomorphism $\xi\colon F^{*}(\mathscr{B},\eta) \Rightarrow \tilde F^{*}\bigl(\tilde{\mathscr{B}},\tilde \eta\bigr)$.
The equivalence relation imposed in Definition~\ref{Freed--Hopkins--Teleman's category of twistings} on the morphisms of $\mathrm{FHT}(\Gamma)$ identifies then
\begin{equation*}
\bigl(\Gamma^{\zeta},S,\mathcal{W},\rho\bigr) \sim (\Gamma^{\omega},U,F^{*}(\mathcal{W},\rho)) \sim \bigl(\Gamma^{\omega},U,\tilde F^{*}(\tilde {\mathcal{W}},\tilde \rho)\bigr)\sim \bigl(\Gamma^{\tilde\zeta},\tilde S,\tilde {\mathcal{W}},\tilde \rho\bigr);
\end{equation*}
this shows that our functor is well-defined.

 Next we prove that identity morphisms go to identities.
 Let $(Y,\pi,\mathscr{P},\psi)\in \hc 1\bigl(\sGrbtriv^{+}(\Gamma)\bigr)$, with $\mathcal{F}((Y,\pi,\mathscr{P},\psi))=:(\Gamma^\pi,\pi,\mathcal{L},\lambda)$. Let \smash{$\id_{(Y,\pi,\mathscr{P},\psi)}= (Z,\zeta,\mathscr{B},\eta)$} be the identity morphism, as explained in~\ref{identity-morphisms-explained}.
Thus, \smash{$\mathcal{F}(\id_{(Y,\pi,\mathscr{P},\psi)})=\big[\Gamma^{\zeta},S,\mathcal{W},\rho\big]$}, where $S$ is the weak equivalence $\Gamma^{\zeta} \to \Gamma^\pi \times_{\Gamma} \Gamma^{\pi}$ of Lemma~\ref{Common refinements of common refinements induce weak equivalences NEU}. We consider the functor $G\colon \Gamma^{\pi} \to \Gamma^{\zeta}$ defined on objects by $y \mapsto (y,y)$, and analogously on morphisms. Let $i$, $i_1$, $i_2$ be the maps
considered in~\ref{identity-morphisms-explained}; then, we have $(i \circ G)(y)=(y,\id,y)=\id_y$ and thus, via Remark~\ref{normalized graded-equivariant structures}, $G^{*}i^{*}\mathscr{P}=\id^{*}\mathscr{P}\cong \id_{\mathscr{I}}$. Moreover, we have $(i_1 \circ G)(y_1,\gamma,y_2) =((y_1,\gamma,y_2),(y_2,\id,y_2))$, and similarly, $(i_2 \circ G)(y_1,\gamma,y_2) =((y_1,\id,y_1),(y_1,\gamma,y_2))$, showing that
$
G^{*}(\mathscr{B},\eta)\cong (\id_{\mathscr{I}},\id_{\mathscr{P}})=\id_{(\mathscr{I},\mathscr{P},\psi)}$.
Under the equivalence of Corollary~\ref{central-super-extensions-and-line-2-bundles}, this implies $G^{*}(\mathcal{W},\rho) \cong \id_{(\mathcal{L},\lambda)}$ in $\sExt(\Gamma^{\pi})$. The diagram
 \begin{equation*}
 \xymatrix{
 \Gamma^{\pi} \ar[dr]_-{\Delta} \ar[rr]^-{G} && \Gamma^{\zeta} \ar[dl]^-{S} \\
 & \Gamma^{\pi} \times_{\Gamma} \Gamma^{\pi} &
 }
 \end{equation*}
is strictly commutative; hence, we have the following equality of morphisms in $\mathrm{FHT}(\Gamma)$
\begin{equation*}
\mathcal{F}(\id_{(Y,\pi,\mathscr{P},\psi)})=\big[\Gamma^{\zeta},S,\mathcal{W},\rho\big]=[\Gamma^{\pi},\Delta,\id_{(\mathcal{L},\lambda)}]=\id_{\mathcal{F}((Y,\pi,\mathscr{P},\psi))}.
\end{equation*}
This shows that $\mathcal{F}$ preserves identities.

 It remains to show that $\mathcal{F}$ is compatible with composition. Suppose
\begin{equation*}
\alxydim{@C=6em}{(Y,\pi,\mathscr{P},\psi) \ar[r]^-{(Z,\zeta,\mathscr{B},\eta)} & \bigl(Y',\pi',\mathscr{P}',\psi'\bigr)\ar[r]^-{(Z',\zeta',\mathscr{B}',\eta')} & (Y'',\pi'',\mathscr{P}'',\psi'')}
\end{equation*}
are 1-morphisms and $\bigl(\tilde Z,\tilde \zeta, \tilde{\mathscr{B}},\tilde\eta\bigr)$ is their composition as defined in~\ref{explanation-of-composition}.
We compute
\begin{equation*}
\mathcal{F}([Z',\zeta',\mathscr{B}',\eta']) \circ \mathcal{F}([Z,\zeta,\mathscr{B},\eta])=\big[\Gamma^{\zeta'},S',\mathcal{W}',\rho'\big] \circ \big[\Gamma^{\zeta},S,\mathcal{W},\rho\big]
\end{equation*}
in $\mathrm{FHT}(\Gamma)$, which consists of the fibre product Lie groupoid \smash{$\Gamma^{\zeta} \times_{\Gamma^{\pi'}} \Gamma^{\zeta'}$}, the weak equivalence
\begin{equation*}
 \xymatrix{
 \Gamma^{\zeta} \times_{\Gamma^{\pi'}} \Gamma^{\zeta'} \ar[rr]^-{S \times_{\Gamma^{\pi'}} S'} && \bigl(\Gamma^{\pi} \times_{\Gamma} \Gamma^{\pi'}\bigr) \times_{_{\Gamma^{\pi'}} } \bigl(\Gamma^{\pi'} \times_{\Gamma} \Gamma^{\pi''}\bigr)
 \ar[rr]^-{\pr_{14}} && \Gamma^{\pi} \times_{\Gamma} \Gamma^{\pi''},
 }
\end{equation*}
and of the 1-morphism given by the composition of the pullbacks of $(\mathcal{W},\rho)$ and $(\mathcal{W}',\rho')$ to \smash{$\Gamma^{\zeta} \times_{\Gamma^{\pi'}} \Gamma^{\zeta'}$}.
On the other hand, we set \smash{$\big[\Gamma^{\tilde\zeta},\tilde S,\tilde {\mathcal{W}},\tilde\rho\big] := \mathcal{F}\bigl(\big[\tilde Z,\tilde \zeta, \tilde{\mathscr{B}},\tilde\eta\big]\bigr)$}. Lemma~\ref{Common refinements of common refinements induce weak equivalences NEU} provides a~weak equivalence \smash{$\Gamma^{\tilde\zeta} \to \Gamma^{\zeta} \times_{\Gamma^{\pi'}} \Gamma^{\zeta'}$} under which the 1-morphisms coincide. This implies
\begin{equation*}
\big[\Gamma^{\tilde\zeta},\tilde S,\tilde {\mathcal{W}},\tilde\rho\big] = \big[\Gamma^{\zeta'},S',\mathcal{W}',\rho'\big] \circ \big[\Gamma^{\zeta},S,\mathcal{W},\rho\big]
\end{equation*}
in $\mathrm{FHT(\Gamma)}$, showing that $\mathcal{F}$ respects the composition. It remains to prove the following.

\begin{Proposition}
\label{The FHT functor: essential surjectivity}
\label{The FHT functor: fullness}
\label{The FHT functor: faithfulness}
The functor $\mathcal{F}$ is an equivalence of categories.
\end{Proposition}
\begin{proof}
We prove the claim in the usual three parts. We start by proving that $\mathcal{F}$ is essentially surjective.
 Let $(\Lambda, T,\mathcal{L},\lambda)$ be a~FHT-twisting. According to Proposition~\ref{Lemma: Extracting covering groupoids from weak equivalences}, there is a~surjective submersion $\pi\colon Y \to \Gamma_0$ and a~strong equivalence $F \colon \Gamma^{\pi} \to \Lambda $ such that $T \circ F \cong P^{\pi}$.
 Let $(\mathscr{I},\mathscr{P},\psi)\in \sGrbtriv(\Gamma^{\pi})$ correspond to the central super extension $F^{*}(\mathcal{L},\lambda)$ of $\Gamma^{\pi}$, under the equivalence of Corollary~\ref{central-super-extensions-and-line-2-bundles}.
 Then, $(Y,\pi,\mathscr{P},\psi)$ is an object in $\hc 1\bigl(\sGrbtriv^{+}(\Gamma)\bigr)$, which we claim to be an essential preimage of $(\Lambda, T,\mathcal{L},\lambda)$ under $\mathcal{F}$.
 By definition of $\mathcal{F}$, we have $\mathcal{F}(Y,\pi,\mathscr{P},\psi)=(\Gamma^{\pi},\pi,F^*(\mathcal{L},\lambda))$. It remains to prove that
 $(\Gamma^{\pi},\pi,F^*(\mathcal{L},\lambda)) \cong (\Lambda, T,\mathcal{L},\lambda)$ in $\mathrm{FHT}(\Gamma)$. Indeed, an isomorphism between them is given by the Lie groupoid $\Omega := \Gamma^\pi \times_\Gamma \Gamma^{\pi}$,
 the weak equivalence
 \begin{equation*}
 S:= \id \times_{\Gamma} F \colon \ \Gamma^\pi \times_\Gamma \Gamma^{\pi} \to \Gamma^\pi \times_\Gamma \Lambda,
 \end{equation*}
 constructed via Lemma~\ref{Lemma: weak equivalences induce weak equivalences of fibre products},
 and the isomorphism
$
 S^*\pr_{1}^* F^*(\mathcal{L},\lambda) \cong S^*\pr_{2}^* (\mathcal{L},\lambda)
$
 of central super extensions of $\Omega$, which is induced from the equality
$F \circ \pr_1 \circ S = \pr_2 \circ S$ of functors~${\Omega \to \Lambda}$.

Next we show that the functor $\mathcal{F}$ is full.
 Let $(Y,\pi,\mathscr{P},\psi)$ and $\bigl(Y',\pi',\mathscr{P}',\psi'\bigr)$ be objects in~$\sGrbtriv^{+}(\Gamma)$ and let
\begin{equation*}
 [\Omega,S,\mathcal{W},\rho] \colon\ \bigl(\Gamma^\pi,\pi,\mathcal{L},\lambda\bigr)\to \bigl(\Gamma^{\pi'},{\pi'},\mathcal{L}',\lambda'\bigr)
\end{equation*}
be a~morphism between their images in $\mathrm{FHT}(\Gamma)$. By definition, $S\colon \Omega \to \Gamma^{\pi} \times_\Gamma \Gamma^{\pi'}$ is a~weak equivalence.
By Proposition~\ref{Lemma: Extracting covering groupoids from weak equivalences}, there exists a~surjective submersion
\begin{equation*}
\psi\colon\ W \to \bigl(\Gamma^{\pi} \times_\Gamma \Gamma^{\pi'}\bigr)_0=Y \ttimes{\pi}{s} \Gamma_1 \ttimes{t}{\pi'} Y',
\end{equation*}
and a~strong equivalence $F\colon \bigl(\Gamma^{\pi} \times_\Gamma \Gamma^{\pi'}\bigr)^{\psi} \to \Omega$ such that $S \circ F \cong {\psi}$. Thus, we have
\begin{equation*}
[\Omega,S,\mathcal{W},\rho] = \big[\bigl(\Gamma^{\pi} \times_\Gamma \Gamma^{\pi'}\bigr)^{\psi},{\psi},F^{*}(\mathcal{W},\rho)\big].
\end{equation*}
Next we consider the weak equivalence $H\colon \Gamma^{\tilde\pi} \to \Gamma^{\pi} \times_\Gamma \Gamma^{\pi'}$ of Lemma~\ref{Common refinements of common refinements induce weak equivalences NEU}, where $\tilde\pi$ is the projection $\tilde\pi\colon Y \times_{\Gamma_0} Y' \to \Gamma_0$, and pull back $\psi$ along $H$, see Lemma~\ref{pullback-of-coverings}. This yields a~surjective submersion $\zeta\colon Z \to Y \times_{\Gamma_0} Y'$ and a~weak equivalence \smash{$G\colon \bigl(\Gamma^{\tilde\pi}\bigr)^{\zeta} \to \bigl(\Gamma^{\pi} \times_\Gamma \Gamma^{\pi'}\bigr)^{\psi}$}. Noting that \smash{$\bigl(\Gamma^{\tilde\pi}\bigr)^{\zeta} = \Gamma^{\tilde\pi \circ \zeta}$}, we obtain another equality
\begin{equation*}
\big[\bigl(\Gamma^{\pi} \times_\Gamma \Gamma^{\pi'}\bigr)^{\psi},{\psi},F^{*}(\mathcal{W},\rho)\big] = \big[\Gamma^{\tilde\pi \circ \zeta},H \circ {\zeta},G^{*}F^{*}(\mathcal{W},\rho)\big].
\end{equation*}
The right-hand side is the image under $\mathcal{F}$ of the 1-morphism $(Z,\zeta,\mathscr{B},\eta)$, where $(\mathscr{B},\eta)$ corresponds to the 1-morphism of central super extensions $G^{*}F^{*}(\mathcal{W},\rho)$ under the equivalence of Corollary~\ref{central-super-extensions-and-line-2-bundles}.

Finally, we prove that
the functor $\mathcal{F}$ is faithful.
 Let
 \begin{equation*}
 (Z,\zeta,\mathscr{B},\eta),\bigl(\tilde Z,\tilde{\zeta},\tilde{\mathscr{B}},\tilde\eta\bigr) \colon\ (Y,\pi,\mathscr{P},\psi) \to \bigl(Y',\pi',\mathscr{P}',\psi'\bigr)
 \end{equation*}
 be morphisms in \smash{$\sGrbtriv^{+}(\Gamma)$}
 such that their images under $\mathcal{F}$ coincide, i.e., $\big[\Gamma^{\zeta},S,\mathcal{W},\rho\big] =\smash{\big[\Gamma^{\tilde \zeta},\tilde S,\tilde {\mathcal{W}},\tilde \rho\big]}$ in $\mathrm{FHT}(\Gamma)$. This means that there is a~smooth functor $G$, a~smooth
 natural transformation $\kappa$ fitting into the diagram
 \begin{equation}
 \label{FHT-proof-first-diagram}
 \begin{gathered}
 \xymatrix@C=1em{
 \Gamma^{\zeta} \ar[rr]^-{G} \ar[dr]_--{S} \ar@{}[dr]_--{}="1" && \Gamma^{\tilde{\zeta}} \ar@{=>}"1"|*+{\kappa} \ar[dl]^-{\tilde{S}} \\
 &\Gamma^{\pi} \times_{\Gamma} \Gamma^{\pi'}, &
 }
 \end{gathered}
 \end{equation}
 and a~2-isomorphism \smash{$\phi \colon G^*\bigl(\tilde{\mathcal{W}},\tilde{\rho}\bigr) \Rightarrow (\mathcal{W},\rho)$} in $\sExt\bigl(\Gamma^{\zeta}\bigr)$, under the identifications $S^{*}\pr_1^{*}(\mathcal{L},\lambda)\cong G^{*}\tilde S^{*}\pr_1^{*}(\mathcal{L},\lambda)$ and $S^{*}\pr_2^{*}(\mathcal{L}',\lambda')\cong G^{*}\tilde S^{*}\pr_2^{*}(\mathcal{L}',\lambda')$ provided by $\kappa$. We consider the fibre product \smash{$W:= Z \times_{Y \times_{\Gamma_0} Y'} \tilde Z$} with its projection $\omega\colon W \to \Gamma_0$ as the trivial refinement of itself. By Lemma~\ref{Common refinements of common refinements induce weak equivalences NEU}, we have a~weak equivalence \smash{$\Gamma^{\omega} \cong \Gamma^{\zeta} \times_{\Gamma^{\pi} \times_{\Gamma} \Gamma^{\pi'}} \Gamma^{\tilde\zeta}$}, and thus, a~smooth natural transformation
\begin{equation}
\label{FHT-proof-second-diagram}
\alxydim{@C=0.5em}{&\Gamma^{\omega} \ar[dr]^{\omega_{\tilde Z}}\ar[dl]_{\omega_Z} \\ \Gamma^{\zeta} \ar[dr] && \Gamma^{\tilde\zeta} \ar@{=>}[ll]|*+{\lambda} \ar[dl] \\ & \Gamma^{\pi} \times_{\Gamma} \Gamma^{\pi'}.}
\end{equation}
Now, the diagram in~\eqref{FHT-proof-first-diagram} is a~sub-diagram of the one in~\eqref{FHT-proof-second-diagram}, and
by Remark~\ref{Lemma: facts about weak equivalences}\,\ref{2-out-of-3}, all functors in both diagrams are weak equivalences. Using Remark~\ref{Lemma: facts about weak equivalences}\,\ref{Weak equivalence filling}, this implies that there is also a~smooth natural transformation
in the complementary sub-diagram,
 \begin{equation*}
 \xymatrix{
 & \Gamma^{\omega} \ar[dr]^-{\omega_{\tilde Z}}="1" \ar@{=>}"1";[dl]|*+{\beta} \ar[dl]_-{\omega_Z} &\\
 \Gamma^{\zeta} \ar[rr]_-{G} && \Gamma^{\tilde{\zeta}}.
 }
 \end{equation*}
The pullback of $\phi$ along $\omega_Z$ hence gives a~2-isomorphism in $\sExt(\Gamma^{\omega})$ between $\omega_{\tilde Z}^{*}\bigl(\tilde{\mathcal{W}},\tilde{\rho}\bigr)$ and $\omega_Z^{*}(\mathcal{W},\rho)$. Under the equivalence of Corollary~\ref{central-super-extensions-and-line-2-bundles}, it corresponds to a~2-isomorphism in $\sGrbtriv(\Gamma^{\omega})$ between $\bigl(\tilde{\mathscr{B}},\tilde\eta\bigr)$ and $(\mathscr{B},\eta)$. This proves that
\[
[Z,\zeta,\mathscr{B},\eta]=\big[\tilde Z,\tilde{\zeta},\tilde{\mathscr{B}},\tilde\eta\big]
\]
 in $\hc 1 \bigl(\sGrbtriv^{+}(\Gamma)\bigr)$.
\end{proof}

\subsection{Moutuou's real twistings}
\label{Moutuou Real twistings}

In this section, we compare our twistings to the real twistings of Moutuou, see~\cite{moutuou2011twistings,mohamedmoutuou2012}, which are defined over real groupoids. We associate to each real Lie groupoid a~graded Lie groupoid,
and lift this to associate to each real twisting a~graded-equivariant super 2-line bundle.

We recall that a~\emph{real space} is a~topological space $X$ together with a~continuous involution~${\tau \colon X \to X}$. A \emph{real map} between real spaces is a~continuous map that commutes with the involutions.
We will work here in a~smooth setting, where $X$ is a~smooth manifold and~$\tau$ is smooth, and call this a~\emph{real manifold}.
The notion of a~real vector bundle was introduced by Atiyah in~\cite{Atiyah1966Reality}. A \emph{real vector bundle} over a~real manifold $(X,\tau)$ is a~complex vector bundle~${\pi \colon E \to X}$ in which $E$ is a~real manifold such that $\pi$ is real, and, for all $x \in X$, the involution~${E_x \to E_{\tau(x)}}$ is anti-linear.
A morphism of real vector bundles is a~real morphism of vector bundles.
We denote the category of real vector bundles over $(X,\tau)$ by $\RVectBdl(X,\tau)$. Before we continue with real twistings, we remark that real vector bundles are included in our setting, namely, in Definition~\ref{Definition of our twisted vector bundles}.
\begin{Proposition}
For $(X,\tau)$ a~real manifold, there is a~canonical isomorphism of categories
\begin{equation*}
\mathrm{R}\VectBdl(X,\tau) \cong \VectBdl(\act X\Z_2,\phi),
\end{equation*}
where $\act X\Z_2$ is the action groupoid
associated to the involution $\tau \colon X \to X$, and $\phi$ is the grading induced by the projection $X \times \Z_2 \to \Z_2$.
\end{Proposition}
\begin{proof}
 An object in $\VectBdl(\act X\Z_2,\phi)=\Hom_{\twoLineBdl(\act X\Z_2,\phi)}((\mathscr{I},\id,\id), (\mathscr{I},\id,\id))$ consists according to Definition~\ref{equivariant 1-morphisms} of a~vector bundle $\mathcal{W}$ over $X$ together with a~linear isomorphism
 $\mathcal{W}_x \to \overline{\mathcal{W}_{\tau(x)}}$ such that the composition $\mathcal{W}_x \to \overline{\mathcal{W}_{\tau(x)}} \to \mathcal{W}_x$ is the identity. Put differently, we have an antilinear
 map $\mathcal{W}_x \to \mathcal{W}_{\tau(x)}$ which turns $\mathcal{W}$ into a~real vector bundle. It is straightforward to complete this to the claimed isomorphism.
\end{proof}

Next we continue with twistings, and first recall from~\cite{mohamedmoutuou2012} that a~\emph{real structure} on a~Lie groupoid $\Gamma$ is a~smooth functor $\tau \colon \Gamma \to \Gamma$ such that
$\tau^2 = \id$ holds. A Lie groupoid together with a~real structure is called a~\emph{real Lie groupoid}.
Moreover, a~\emph{real functor} between real Lie groupoids $(\Gamma,\tau) \to (\Gamma',\tau')$ is a~smooth functor $F \colon \Gamma \to \Gamma'$ such that
$\tau' \circ F = F \circ \tau$ holds. We denote the category of real Lie groupoids and real functors by $\Rliegrpd$.

The definition of real twistings of Moutuou carries out a~real version of Freed--Hopkins--Teleman's construction recalled in Section~\ref{Freed--Hopkins--Teleman-twistings}.
We start by giving the following definition, which is a~slight reformulation of \cite[Definition 2.5.1]{mohamedmoutuou2012}.
We recall that a~central super extension of a~groupoid $\Gamma$ is a~$\phi$-twisted super extension of $\Gamma$ in the
sense of
 Definition~\ref{Definition: twisted super extension} (for $\phi=1$), and that central super extensions of $\Gamma$ form a~category $\sExt^{\mathrm{ref}}(\Gamma) := \mathrm{FM}(\Gamma,1)$, see Remark~\ref{sExtref}.
\begin{Definition}
 A \emph{real central super extension} of $(\Gamma,\tau)$ consists of a~central super extension~${(\mathcal{L},\lambda)}$ of $\Gamma$ together with a~real structure $\beta\colon \mathcal{L} \to \overline{\tau^*\mathcal{L}}$ on the super line bundle $\mathcal{L}$
 that is compatible with $\lambda$ in the sense that the diagram
 \begin{equation*}
 \xymatrix{
 \mathcal{L}_{\gamma_1} \otimes \mathcal{L}_{\gamma_2} \ar[rr]^-{\lambda_{\gamma_1,\gamma_2}} \ar[d]_-{\beta_{\gamma_1} \otimes \beta_{\gamma_2}} && \mathcal{L}_{\gamma_1 \circ \gamma_2} \ar[d]^-{\beta_{\gamma_1 \circ \gamma_2}} \\
 \overline{\mathcal{L}}_{\tau(\gamma_1)} \otimes \overline{\mathcal{L}}_{\tau(\gamma_2)} \ar[rr]_-{\overline{\lambda}_{\tau(\gamma_1),\tau(\gamma_2)}} && \overline{\mathcal{L}}_{\tau(\gamma_1 \circ \gamma_2)}
 }
 \end{equation*}
 commutes for all composable $\gamma_1,\gamma_2 \in \Gamma_1$.
 If $(\mathcal{L},\lambda,\beta)$ and $(\mathcal{L}',\lambda',\beta')$ are real central super extensions of $(\Gamma,\tau)$, then an \emph{isomorphism} $(\mathcal{L},\lambda,\beta) \to (\mathcal{L}',\lambda',\beta')$
 is a~refinement $\eta \colon (\mathcal{L},\lambda) \to (\mathcal{L}',\lambda)$ of central super extensions that is compatible
 with $\beta$ in the sense that the diagram
 \begin{equation*}
 \xymatrix{
 \mathcal{L}_{\gamma} \ar[rr]^-{\eta_\gamma} \ar[d]_{\beta_\gamma} && \mathcal{L}_{\gamma}' \ar[d]^{\beta'_\gamma} \\
 \overline{\mathcal{L}}_{\tau(\gamma)} \ar[rr]_-{\overline{\eta}_{\tau(\gamma)}} && \overline{\mathcal{L}'}_{\tau(\gamma)}
 }
 \end{equation*}
 commutes for all $\gamma \in \Gamma_1$.
\end{Definition}

We denote the category of real central super extensions of $(\Gamma,\tau)$ by $\RsExt(\Gamma,\tau)$.
By construction, it comes equipped with a~functor
$
\RsExt(\Gamma,\tau) \to \sExt^{\mathrm{ref}}(\Gamma)
$
to the category of central super extensions, see Remark~\ref{sExtref}; this functor forgets the real structure $\beta$.

Let $(\Gamma,\tau)$ be a~real Lie groupoid, $(Y,\rho)$ be a~real manifold, and $\pi \colon Y \to \Gamma_0$ be a~real surjective submersion. The covering groupoid $\Gamma^\pi$
of Definition~\ref{Example : covering groupoid and weak equivalences} attains a~real structure $\tau^\pi \colon \Gamma^\pi \to \Gamma^\pi$ given on objects by
$y \mapsto \rho(y)$ and on morphisms by $(y,\gamma,y') \mapsto (\rho(y),\tau(\gamma),\rho(y'))$. The projection functor $\pi \colon \Gamma^\pi \to \Gamma$ is a~real functor.

\begin{Definition} \label{Definition: Real Morita equivalence}
 A real functor $F\colon (\Gamma',\tau') \to (\Gamma,\tau)$ between real Lie groupoids is called \emph{real weak equivalence} if there exists a~real surjective
 submersion $\pi \colon Y \to \Gamma_0$ and a~real strong equivalence $S \colon (\Gamma^\pi,\tau^\pi) \to (\Gamma',\tau')$ such that $F \circ S \cong \pi$.
\end{Definition}

\begin{Remark}\label{Real weak equivalences have underlying weak equivalences}
In Definition~\ref{Definition: Real Morita equivalence}, a~real strong equivalence is an isomorphism in the category~$\Rliegrpd$.
 Real weak equivalences between real Lie groupoids are in particular weak equivalences between the underlying Lie groupoids; this follows
 from Proposition~\ref{Lemma: Extracting covering groupoids from weak equivalences}.
\end{Remark}

\begin{Definition}
\label{Moutuou-Real-twisting}
 A \emph{Moutuou real twisting} of a~real groupoid $(\Gamma,\tau)$ consists of another real groupoid $(\Lambda,\tau_{\Lambda})$, a~real weak equivalence
 $F\colon (\Lambda,\tau_{\Lambda}) \to (\Gamma,\tau)$, and a~real central super extension~${(\mathcal{L},\lambda,\beta)}$ of $(\Lambda,\tau_{\Lambda})$. We denote
 the set of Moutuou real twistings by $\Mou(\Gamma,\tau)$.
\end{Definition}

\begin{Remark} \label{Remark: Real twistings have underlying FHT twistings}
 Any Moutuou real twisting of a~real groupoid $(\Gamma,\tau)$ has an underlying FHT-twisting of $\Gamma$ (see Definition~\ref{Freed--Hopkins--Teleman's category of twistings})
 using the forgetful functor $\RsExt(\Gamma,\tau) \to \sExt^{\mathrm{ref}}(\Gamma)$ and Remark~\ref{Real weak equivalences have underlying weak equivalences}.
\end{Remark}

\begin{Remark}
 Definition~\ref{Moutuou-Real-twisting} gives -- a~priori -- only a~special case of Moutuou's twistings, where~${F\colon (\Lambda,\tau_{\Lambda}) \to (\Gamma,\tau)}$ is allowed to be a~\emph{real Morita equivalence} \cite[Definitions~2.4.1 and~2.6.3]{mohamedmoutuou2012}.
In general, real Morita equivalences, are \emph{spans} of so-called \emph{real generalized morphisms}. However, if
\begin{equation*}
(\Lambda,\tau_{\Lambda}) \stackrel G\leftarrow (\Omega,\omega) \stackrel F\to (\Gamma,\tau)
\end{equation*}
is a~span of real generalized morphisms, then one can first replace $(\Omega,\omega)$ by another real Lie groupoid such that the real generalized morphisms $G$ and $F$ become real weak equivalences in the sense of Definition~\ref{Definition: Real Morita equivalence} \cite[Proposition~2.4.9]{mohamedmoutuou2012}.
In a~second step, if a~general Moutuou real twisting comes with a~real central super extension $(\mathcal{L},\lambda,\beta)$ of $(\Lambda,\tau_{\Lambda})$, then one may pull back~${(\mathcal{L},\lambda,\beta)}$ along $G$ to obtain one of the more special twistings we defined in Definition~\ref{Moutuou-Real-twisting}. We believe that -- under a~suitable notion of isomorphisms between Moutuou twistings, analogous to Definition~\ref{Freed--Hopkins--Teleman's category of twistings} -- the two steps outlined above will result in isomorphisms, and show that every general Moutuou twisting is isomorphic to one defined above.
\end{Remark}

In order to compare Moutuou real twistings to super 2-line bundles, we first associate to each real Lie groupoid a~graded Lie groupoid, via the following construction. First of all, we recall that a~strict action of a~finite group $G$ on a~Lie groupoid $\Gamma$ consists of smooth functors~${g_* \colon \Gamma \to \Gamma}$, for each $g\in G$, such
 that $1_* = \id_{\Gamma}$ and $(hg)_* = h_* g_*$ for all $g,h \in G$. The \emph{semi-direct product} $\Gamma \rtimes G$ has the manifold of objects $\Gamma_0$ and the manifold of morphisms $\Gamma_1 \times G$, with~${s(\gamma,g) := g_{*}^{-1}(s(\gamma))}$ and $t(\gamma,g) := t(\gamma)$, and composition $(\delta,h) \circ (\gamma,g) := (\delta \circ (h_{*}(\gamma)), hg)$.
Now let $(\Gamma,\tau)$ be a~real groupoid. We regard $\tau$ as a~$\Z_2$-action on $\Gamma$, and consider the corresponding semidirect product~${\Gamma \rtimes_{\tau} \Z_2}$, which we equip with a~grading $\phi_{\tau}$ by projection to the second factor. We denote this graded Lie groupoid by $\Gr(\Gamma,\tau)$. Likewise, to a~real functor $F \colon (\Gamma,\tau) \to (\Gamma',\tau')$ we assign the smooth
functor $\Gamma \rtimes_{\tau} \Z_2 \to \Gamma' \rtimes_{\tau'} \Z_2$, defined on objects by $x \mapsto F(x)$ and on morphisms by~${(\gamma,g) \mapsto (F(\gamma),g)}$. This functor is even with respect to the gradings $\phi_{\tau}$ and $\phi_{\tau'}$.
We have the following result.
\begin{Lemma}
\label{lemma: Real Morita equivalences give rise to even weak equivalences between the associated semidirect product groupoids}
 The assignment $(\Gamma,\tau) \mapsto \Gr(\Gamma,\tau) : =(\Gamma \rtimes_{\tau} \Z_2,\phi_{\tau})$ extends to a~faithful functor
 \begin{equation*}
 \Gr\colon\ \Rliegrpd \to\Grliegrpd,
 \end{equation*}
 and it sends real weak equivalences to even weak equivalences.
\end{Lemma}

\begin{proof}
 Let $(\Lambda,\tau_{\Lambda})$ and $(\Gamma,\tau)$ be real Lie groupoids and $F \colon (\Lambda,\tau_{\Lambda}) \to (\Gamma,\tau)$ a~real weak equivalence.
 Thus, there exists a~real surjective
 submersion $\pi \colon Y \to \Gamma_0$ and a~real strong equivalence~${S \colon (\Gamma^\pi,\tau^\pi) \to (\Lambda,\tau_{\Lambda})}$ such that $F \circ S \cong \pi$.
Since functors send isomorphisms to isomorphisms, our functor sends real strong equivalences to even strong equivalences. Moreover, it is straightforward to show that the image of $\pi\colon \Gamma^{\pi} \to \Gamma$,
$
 (\Gamma^\pi \rtimes_{\tau^\pi} \Z_2,\phi_{\tau^{\pi}}) \to (\Gamma \rtimes_{\tau} \Z_2,\phi_{\tau})$,
 is a~weak equivalence. Via the 2-out-of-3 property of weak equivalences, this shows that the image of $F$ is a~weak equivalence.
\end{proof}

\begin{Remark}
 The functor $\Gr$ is not essentially surjective. Indeed, consider the group homomorphism $\Z_4 \to \Z_2$, $x \mapsto x \mod 2$ which defines a~graded Lie groupoid
 $(\mathrm{B}\Z_4,\bmod 2)$. Since the group extension
$
 \xymatrix{
 0 \ar[r] & \Z_2 \ar[r] & \Z_4 \ar[r] & \Z_2 \ar[r] & 0
 }
$
 is non-split, $\Z_4$ is not a~semidirect product.
\end{Remark}

\begin{Remark}
\label{graded-groupoid-functor-inclusion}
If $(\Gamma,\tau)$ is a~real Lie groupoid, then the underlying groupoid can be obtained from the associated graded Lie groupoid as the kernel of its grading,
$
\Gamma = \Gr(\Gamma,\tau)^{\rm even}$.
In particular, there is an even functor $(\Gamma,1) \to \Gr(\Gamma,\tau)$.
\end{Remark}

\begin{Remark}
 Regarding a~real manifold $(X,\tau)$ as a~real Lie groupoid $(X_{\rm dis},\tau_{\rm dis})$ with only identity morphisms, the graded Lie groupoid $\Gr(X_{\rm dis},\tau_{\rm dis})$ coincides with the graded action groupoid $\mathrm{Gr}(X,\tau) =\act X{(\Z_2,\id)}$ associated to $(X,\tau)$ considered in Example~\ref{ex1}\,\ref{real-manifold}.
\end{Remark}

Next, we compare real central super extensions and super 2-line bundles. Let $(\Gamma,\tau)$ be a~real Lie groupoid and $(\mathcal{L},\lambda,\beta)$ a~real central super extension of
$(\Gamma,\tau)$. We construct a~super 2-line bundle over the associated graded Lie groupoid $\Gr(\Gamma,\tau)$. The underlying super 2-line bundle over $(\Gamma \rtimes_{\tau} \Z_2)_0 = \Gamma_0$ is $\mathscr{I}$, the trivial one.
We define a~super line bundle over $(\Gamma \rtimes_{\tau} \Z_2)_1$ by setting
\[
 \tilde{\mathcal{L}}_{(\gamma,x)} :=
 \begin{cases}
 \mathcal{L}_\gamma & \text{if} \ x=+1 ,\\
 \mathcal{L}_{\tau(\gamma)} &\text{if} \ x=-1.
 \end{cases}
\]
It corresponds to the required 1-isomorphism
$\mathscr{P} \colon s^*\mathscr{I}=\mathscr{I} \to \mathscr{I} =(t^*\mathscr{I})^\phi$ over $(\Gamma \rtimes_{\tau} \Z_2)_1$. Next we construct the 2-isomorphism $\psi$, which corresponds to a~super line bundle isomorphism $\tilde\lambda$ over~${(\Gamma \rtimes_{\tau} \Z_2)_2}$, fibre-wise
\begin{equation*}
 (\tilde{\lambda}_\beta)_{(\gamma,x),(\gamma',x')} \colon\ \overline{\tilde{\mathcal{L}}_{(\gamma,x)}}^{x'} \otimes \tilde{\mathcal{L}}_{(\gamma',x')} \to \tilde{\mathcal{L}}_{(\gamma \circ (x \cdot \gamma'),xx')}.
\end{equation*}
Considering the four possible signs of $x$ and $x'$ separately, $\tilde\lambda_{\beta}$ is defined in the following ways:
\begin{itemize}\itemsep=0pt
\item
$\xymatrix@C=2em{
 \mathcal{L}_{\gamma_1} \otimes \mathcal{L}_{\gamma_2} \ar[rr]^-{\lambda_{\gamma_1,\gamma_2}} && \mathcal{L}_{\gamma_1 \circ \gamma_2},}$ $x=+1$, $x'=+1$,

\item
$\xymatrix@C=2em{
 \mathcal{L}_{\tau(\gamma_1)} \otimes \mathcal{L}_{\gamma_2} \ar[rr]^-{\lambda_{\tau(\gamma_1),\gamma_2}} && \mathcal{L}_{\tau(\gamma_1 \circ \tau(\gamma_2))},
 }$ $x=-1$, $x'=+1$,

\item
$\xymatrix@C=2em{
 \overline{\mathcal{L}_{\gamma_1}} \otimes \mathcal{L}_{\tau(\gamma_2)} \ar[rr]^-{\overline{\beta_{\gamma_1}} \otimes \id} && \mathcal{L}_{\tau(\gamma_1)} \otimes \mathcal{L}_{\tau(\gamma_2)}
 \ar[rr]^-{\lambda_{\tau(\gamma_1),\tau(\gamma_2)}} && \mathcal{L}_{\tau(\gamma_1 \gamma_2)},
 }$ $x=+1$, $x'=-1$,

\item
$\xymatrix@C=2em{
 \overline{\mathcal{L}_{\tau(\gamma_1)}} \otimes \mathcal{L}_{\tau(\gamma_2)} \ar[rr]^-{\beta_{\gamma_1}^{-1} \otimes \id} && \mathcal{L}_{\gamma_1} \otimes \mathcal{L}_{\tau(\gamma_2)}
 \ar[rr]^-{\lambda_{\gamma_1, \tau(\gamma_2)}} && \mathcal{L}_{\gamma_1 \tau(\gamma_2)},
 }$ $x=-1$, $x'=-1$.
\end{itemize}
One checks that these pieces patch together to an isomorphism of super line bundles over ${(\Gamma \rtimes_{\tau} \Z_2)_2}$, using the compatibility of $\lambda$ with the grading and the real structure $\beta$.
Further, one deduces the associativity property of $\tilde{\lambda}_{\beta}$ from the associativity condition of $\lambda$. This finalizes the construction of a~super 2-line bundle $(\mathscr{I},\mathscr{P},\psi)$ over $\Gr(\Gamma,\tau)$. We find the following result.
\begin{Theorem}
\label{Lemma: Real central extensions are super 2-line bundles}
\label{Proposition: Moutou real extensions embed into trivial super 2-line bundles}
 The assignment described above extends to an essentially injective and faithful functor
$
 \RsExt(\Gamma,\tau) \to \stwoLineBdltriv^{\mathrm{ref}}(\Gr(\Gamma,\tau))
$
 from the category of real central super extensions of~$(\Gamma,\tau)$ to the category of trivial super $2$-line bundles over
 $\Gr(\Gamma,\tau)$ with refinements as morphisms.
\end{Theorem}
\begin{proof}
 Let $\theta \colon (\mathcal{L},\lambda,\beta) \to (\mathcal{L}',\lambda',\beta')$ be an isomorphism of real central super extensions. We construct a~refinement $(\mathscr{R},\eta) \colon \bigl(\mathscr{I},\tilde{\mathcal{L}},\tilde{\lambda}_{\beta}\bigr) \to \bigl(\mathscr{I},\tilde{\mathcal{L}'},\tilde{\lambda'}_{\beta'}\bigr)$. The refinement
 $\mathscr{R} \colon \mathscr{I} \to \mathscr{I}$ of super 2-line bundles over $(\Gamma \rtimes_{\tau} \Z_2)_0 = \Gamma_0$ is taken to be the identity morphism
 of $\mathscr{I}$. The 2-isomorphism $\eta$ is then defined fibrewise over $(\gamma,x) \in (\Gamma \rtimes_{\tau} \Z_2)_1$ as
 \[
 \eta \colon\ \tilde{\mathcal{L}'} \to \tilde{\mathcal{L}}, \qquad \eta_{(\gamma,x)} :=
 \begin{cases}
 \theta_{\gamma}^{-1} \colon \mathcal{L}'_\gamma \to \mathcal{L}_\gamma & \text{if} \ x=+1, \\
 \theta_{\tau(\gamma)}^{-1} \colon \mathcal{L}'_{\tau(\gamma)} \to \mathcal{L}_{\tau(\gamma)} & \text{if} \ x=-1.
 \end{cases}
 \]
 The compatibility of $\theta$ with $\lambda$, $\lambda'$ and $\beta$, $\beta'$ ensures that this yields a~refinement $(\mathscr{R},\eta)$.
 One checks that the assignment is functorial.

It is clear that the functor is faithful, so it remains to show that it is essentially injective.
 Let $(\mathcal{L},\lambda,\beta)$, $(\mathcal{L}',\lambda',\beta')$ be real central super extensions of $(\Gamma,\tau)$ such that there exists an isomorphism
 \smash{$\eta \colon \bigl(\mathscr{I},\tilde{\mathcal{L}},\tilde{\lambda}_{\beta}\bigr) \to \bigl(\mathscr{I},\tilde{\mathcal{L}'},\tilde{\lambda'}_{\beta'}\bigr)$}. We construct an isomorphism
 $(\mathcal{L},\lambda,\beta) \to (\mathcal{L}',\lambda',\beta')$ as follows.
 The isomorphism $\mathcal{L} \to \mathcal{L}'$ over $\Gamma_1$ is obtained by pulling back the isomorphism $\eta \colon \tilde{\mathcal{L}} \to \tilde{\mathcal{L}'}$ over $(\Gamma \rtimes_{\tau} \Z_2)_1$
 along the inclusion $\iota_1 \colon \Gamma_1 \to (\Gamma \rtimes_{\tau} \Z_2)_1$, noting that \smash{$(\iota_1^*\tilde{\mathcal{L}})_{\gamma} = \tilde{\mathcal{L}}_{(\gamma,1)} = \mathcal{L}_\gamma$} by definition
 and similar for $\tilde{\mathcal{L}'}$. The compatibility of $\iota_1^*\eta$ with $\lambda$, $\lambda'$ is then encoded in the pullback of the compatibility diagram
 of $\eta$ with $\tilde{\lambda}_{\beta}$, $\tilde{\lambda'}_{\beta'}$ along the inclusion $\iota_2 \colon \Gamma_2 \to (\Gamma \rtimes_{\tau} \Z_2)_2$. It remains to
 check that $\iota_1^*\eta$ is compatible with the real structures. Here we use that the real structures on $\mathcal{L}$, $\mathcal{L}'$ are essentially encoded in
 $\tilde{\lambda}_{\beta}, \tilde{\lambda'}_{\beta'}$. In detail, we recover the real structure on~$\mathcal{L}$ from $\tilde{\lambda}_{\beta}$ (and the symmetric braiding)
 as follows:
 \begin{align*}
 \mathcal{L}_\gamma
 &= \mathcal{L}_{(\gamma,1)} \cong \overline{\mathcal{L}_{(\id,-1)}} \otimes \mathcal{L}_{(\tau(\gamma),-1)}\cong \overline{\mathcal{L}_{(\id,-1)}} \otimes \overline{\mathcal{L}_{(\tau(\gamma),1)}} \otimes \mathcal{L}_{(\id,-1)} \\
 &\cong \overline{\mathcal{L}_{(\id,-1)}} \otimes \mathcal{L}_{(\id,-1)} \otimes \overline{\mathcal{L}_{(\tau(\gamma),1)}} \cong \mathcal{L}_{(\id,1)} \otimes \overline{\mathcal{L}_{(\tau(\gamma),1)}} \cong \overline{\mathcal{L}_{(\tau(\gamma),1)}}= \overline{\mathcal{L}_{\tau(\gamma)}}.
 \end{align*}
 The compatibility of $\iota_1^*\eta$ with the real structures then follows from the compatibility of $\eta$ with~$\tilde{\lambda}_{\beta}$,~$\tilde{\lambda'}_{\beta'}$.
 This shows that $\iota_1^*\eta \colon (\mathcal{L},\lambda,\beta) \to (\mathcal{L}',\lambda',\beta')$ is an isomorphism of real central super extensions.
\end{proof}

In combination with our Theorem~\ref{The main theorem of the Freed-Moore chapter}, we obtain the following result.

\begin{Corollary}
There is a~canonical essentially injective and faithful functor
\begin{equation*}
\RsExt(\Gamma,\tau)\to \mathrm{FM}(\Gr(\Gamma,\tau))
\end{equation*}
from the category of real central super extensions
 to the category of Freed--Moore's twisted groupoid extensions. \end{Corollary}

Next, we compare Moutuou's real twistings to ours: we associate to each Moutuou real twisting of a~real groupoid~$(\Gamma,\tau)$ a~super 2-line bundle over the associated graded Lie groupoid
$\Gr(\Gamma,\tau)=(\Gamma \rtimes_{\tau} \Z_2, \phi_{\tau})$ as follows. Let $(\Lambda,\tau_{\Lambda},F,\mathcal{L},\lambda,\beta)$ be a~real twisting of $(\Gamma,\tau)$.
According to Theorem~\ref{Lemma: Real central extensions are super 2-line bundles},
the real central super extension $(\mathcal{L},\lambda,\beta)$ of $(\Lambda,\tau_{\Lambda})$ gives rise to a~super 2-line bundle
$
 \bigl(\mathscr{I},\tilde{\mathcal{L}},\tilde{\lambda}_{\beta}\bigr) \in \stwoLineBdl(\Lambda \rtimes_{\tau_{\Lambda}} \Z_2,\phi_{\tau_{\Lambda}})$.
Combining Lemma~\ref{lemma: Real Morita equivalences give rise to even weak equivalences between the associated semidirect product groupoids}
with Theorem~\ref{descent-for-graded-equivariant-2-line-bundles}, the real weak equivalence $F\colon (\Lambda,\tau_{\Lambda}) \to (\Gamma,\tau)$ induces an equivalence of bicategories
$\stwoLineBdl(\Lambda \rtimes_{\tau_{\Lambda}} \Z_2,\phi_{\tau_{\Lambda}}) \cong \stwoLineBdl(\Gamma \rtimes_{\tau} \Z_2,\phi_{\tau})$. The image~$(\mathscr{\mathcal{L}},\mathscr{P},\psi)$
of $\bigl(\mathscr{I},\tilde{\mathcal{L}},\tilde{\lambda}_{\beta}\bigr)$ under this equivalence is a~super 2-line bundle over $(\Gamma \rtimes_{\tau} \Z_2,\phi_{\tau})$. Thus, we have the following result.

\begin{Theorem}
\label{main-theorem-moutuou}
Every Moutuou real twisting in the sense of Definition~{\rm\ref{Moutuou-Real-twisting}} over a~real groupoid ${(\Gamma,\tau)}$ gives rise to a~super $2$-line bundle over the corresponding graded Lie groupoid $\Gr(\Gamma,\tau)$, in such a~way that the diagram
\begin{equation*}
 \xymatrix{
 \Mou(\Gamma,\tau) \ar[d]^-{} \ar[rr]^-{} && \hc 0(\stwoLineBdl(\Gr(\Gamma,\tau))) \ar[d] \\
 \hc 0(\mathrm{FHT}(\Gamma)) \ar[rr]^-{} && \hc 0(\stwoLineBdl(\Gamma))
 }
\end{equation*}
commutes. Here, the left vertical arrow is described in Remark~{\rm\ref{Remark: Real twistings have underlying FHT twistings}}, the right vertical arrow is pullback along
the even functor $(\Gamma,1) \to \Gr(\Gamma,\tau)$ of Remark~{\rm\ref{graded-groupoid-functor-inclusion}},
and the lower horizontal arrow is induced from Theorem~{\rm\ref{The main theorem of the Freed Hopkins Teleman chapter}}.
\end{Theorem}

\begin{Remark}
We have not attempted to construct a~full bicategory of Moutuou's real twistings here; we expect that a~proper formulation of
this bicategory of real twistings would allow us to extend the previous construction to a~functor. We note that there is a~notion of equivalence among such twistings
defined in \cite[Definition 2.6.4]{mohamedmoutuou2012}, which we have not included in our comparison since it is rather intricate.
\end{Remark}

\subsection{Chern--Simons twistings} \label{Cherns-Simons theory twistings}
In this section, we review,
based on \cite[Section~4]{freed6}, how Chern--Simons theory gives rise to equivariant twistings and a~corresponding theory of twisted equivariant vector bundles, and we show that both fit into our general setting.

Let $G$ be a~compact Lie group. For maps $G^k \to G^l$ between products of $G$ with itself, we use a~self-explaining notation that combines projections, multiplication, and inversion; e.g., $p_{1,2\cdot3,4^{-1}}\colon G^5 \to G^3$ denotes the map $(g_1,\dots,g_5) \mapsto \bigl(g_1,g_2g_3,g^{-1}_4\bigr)$.
Moreover, the pullback along such a~map will be denoted by $p_{1,2\cdot 3,4^{-1}}^{*}(\dots)=(\dots)_{1,2\cdot 3,4^{-1}}$. The following definition is~\cite[Definition~5.1]{carey4}.

\begin{Definition} \label{definition multiplicative bundle gerbe}
 A \emph{multiplicative bundle gerbe over $G$} is a~bundle gerbe $\mathscr{G}$ over $G$ together with a~1-isomorphism
$
 \mathscr{K} \colon \mathscr{G}_1 \otimes
 \mathscr{G}_2 \to \mathscr{G}_{1\cdot 2}
$
 over $G \times G$,
 and a~2-isomorphism
 \begin{equation*}
 \alxydim{}{
 \mathscr{G}_1 \otimes \mathscr{G}_2 \otimes \mathscr{G}_3 \ar[rr]^-{\mathscr{K}_{1, 2} \otimes \id} \ar[d]_-{\id \otimes \mathscr{K}_{2, 3}} & & \mathscr{G}_{1\cdot 2} \otimes \mathscr{G} \ar@{=>}[dll]|*+{\theta} \ar[d]^-{\mathscr{K}_{1\cdot 2,3}} \\
 \mathscr{G}_1 \otimes \mathscr{G}_{2\cdot 3} \ar[rr]_-{\mathscr{K}_{1,2\cdot 3}} & & \mathscr{G}_{1\cdot 2\cdot 3}
 }
 \end{equation*}
 over $G \times G \times G$ that satisfies an evident Pentagon axiom over $G^4$.
\end{Definition}

\begin{Definition}
 Let $(\mathscr{G},\mathscr{K},\theta)$ and $\bigl(\mathscr{G}',\mathscr{K}',\theta'\bigr)$ be multiplicative bundle gerbes over $G$. A~\emph{$1$-morphism}
 $(\mathscr{G},\mathscr{K},\theta) \to \bigl(\mathscr{G}',\mathscr{K}',\theta'\bigr)$ consists of a~1-morphism $\mathscr{B} \colon \mathscr{G} \to \mathscr{G}'$
 over $G$ and a~2-isomor\-phism
 \begin{equation*}
 \xymatrix{
 \mathscr{G}_1 \otimes \mathscr{G}_2 \ar[d]_-{\mathscr{B}_1 \otimes \mathscr{B}_2} \ar[rr]^-{\mathscr{K}} &&\mathscr{G}_{1\cdot 2} \ar[d]^-{\mathscr{B}_{1\cdot 2}} \ar@{=>}[dll]|-*+{\xi} \\
 \mathscr{G}'_1 \otimes \mathscr{G}'_2 \ar[rr]_-{\mathscr{K}'} && \mathscr{G}_{1\cdot 2}'
 }
 \end{equation*}
 over $G \times G$, satisfying an evident compatibility condition with $\theta$, $\theta'$ over $G ^3$.
If $
(\mathscr{B},\xi),\bigl(\hat{\mathscr{B}},\hat{\xi}\bigr) \colon\allowbreak (\mathscr{G},\mathscr{K},\theta) \to \bigl(\mathscr{G}',\mathscr{K}',\theta'\bigr)
$
 are 1-morphisms, then a~\emph{$2$-morphism} $(\mathscr{B},\xi) \Rightarrow \bigl(\hat{\mathscr{B}},\hat{\xi}\bigr)$ consists of a~2-morphism \smash{$\alpha \colon \mathscr{B} \Rightarrow \hat{\mathscr{B}}$} satisfying an evident compatibility
 condition with $\xi$ and \smash{$\hat{\xi}$} over $G \times G$.
\end{Definition}
It is straightforward to show that multiplicative bundle gerbes, together with the 1-morphisms and 2-morphisms form a~bicategory $\MultGrb(G)$.

Multiplicative bundle gerbes over $G$ are classified up to isomorphism by $\h^4(\mathrm{B}G,\mathbb{Z})$
 \cite[Proposition~5.2]{carey4}. The procedure $(\mathscr{G},\mathscr{K},\theta) \mapsto \mathscr{G}$ of forgetting the multiplicative structure realizes in cohomology the ``transgression'' homomorphism \begin{equation}
 \label{transgression-map}
 \h^4(\mathrm{B}G,\mathbb{Z}) \to \h^3(G,\mathbb{Z}).
 \end{equation}
In other words, a~bundle gerbe over $G$ is multiplicative if and only if its ``Dixmier--Douady-class'' in $\h^3(G,\Z)$ lies in the image of that homomorphism.

The transgression homomorphism~\eqref{transgression-map} factors in fact through the $G$-equivariant cohomology of $G$, $\h^3_G(G,\Z)$, where $G$ acts on itself by conjugation.
We describe in the following how this factorization is realized on the level of geometric objects. We consider a~multiplicative bundle gerbe $(\mathscr{G},\mathscr{K},\theta)$ over $G$.
As $G$ acts on itself via conjugation, we may form the action groupoid~$\act GG$.
A bundle gerbe over $\act GG$ (see Definition~\ref{equivariant structure} and Remark~\ref{graded-equivariant-over-a-manifold}) can now be constructed as follows.
The underlying bundle gerbe over $G$ is $\mathscr{G}$. The 1-isomorphism $\mathscr{P}_{\mathscr{K}}$ over~$G \times G$ is defined by the composition
\begin{equation*}
 \alxydim{@C=2em}{
 \mathscr{G}_2 \ar[r] &\! \mathscr{G}_1^* \otimes \mathscr{G}_1 \otimes \mathscr{G}_2 \ar[rr]^-{\id \otimes \mathscr{K}_{1, 2}} &&\!
 \mathscr{G}_1^* \otimes \mathscr{G}_{1\cdot 2} \ar[rrr]^-{\id \otimes \mathscr{K}^{-1}_{1\cdot 2 \cdot 1^{-1},1}} &&&\!\mathscr{G}_1^* \otimes \mathscr{G}_{1\cdot 2\cdot 1^{-1}} \otimes \mathscr{G}_1
 \ar[r] & \!\mathscr{G}_{1\cdot 2 \cdot 1^{-1}},
 }
\end{equation*}
where the unlabel arrows are given by duality/symmetry in the bicategory of bundle gerbes.
The 2-isomorphism $\psi_{\theta}$ over $G^3$ is induced from $\theta$ and diagram~\eqref{equivariant structure: diagram} commutes because of the pentagon axiom
required in Definition~\ref{definition multiplicative bundle gerbe}.
We find the following.
\begin{Proposition} \label{proposition multiplicative gerbes give rise to twistings}
 The assignment $(\mathscr{G},\mathscr{K},\theta) \mapsto (\mathscr{G},\mathscr{P}_{\mathscr{K}},\psi_{\theta})$ extends to a~functor
 \begin{equation*}
 \MultGrb(G) \to \Grb(\act GG) \subset \stwoLineBdl(\act GG)
 \end{equation*}
 between the bicategory of multiplicative bundle gerbes over $G$ and the bicategory of bundle gerbes over the action groupoid $\act GG$.
\end{Proposition}
\begin{proof}
 Let $(\mathscr{B},\xi) \colon (\mathscr{G},\mathscr{K},\theta) \to \bigl(\mathscr{G}',\mathscr{K}',\theta'\bigr)$ be a~1-morphism of multiplicative bundle gerbes. It is tedious but
 straightforward to construct a~2-isomorphism $\eta \colon \mathscr{P}_{\mathscr{K}'} \circ s^*\mathscr{B} \Rightarrow t^*\mathscr{B} \circ \mathscr{P}_{\mathscr{K}}$ from the 2-isomorphism
 $\xi$ and deduce the coherence property of $\eta$ from the one of $\xi$.
 Given a~2-morphism $\alpha \colon (\mathscr{B},\xi) \Rightarrow \bigl(\hat{\mathscr{B}},\hat{\xi}\bigr)$, it is also straightforward to show that $\alpha$ is compatible with
 $\eta$ and~$\hat{\eta}$. Functoriality follows from a~standard pasting diagram argument.
\end{proof}

\begin{Remark}
 If $G$ is finite, and, correspondingly, $\mathscr{G}$ is the trivial bundle gerbe, then Proposition~\ref{proposition multiplicative gerbes give rise to twistings} reproduces a~construction outlined in \cite[equation (4.5)]{freed6}, using that a~1-isomorphism between trivial bundle gerbes
 is precisely a~line bundle.
\end{Remark}

The $(\mathscr{G},\mathscr{P}_{\mathscr{K}},\psi_{\theta})$-twisted vector bundles (in the sense of Definition~\ref{Definition of our twisted vector bundles})
play an important role in Chern--Simons theory, as we review next, following \cite[Section~4]{freed6}.
For the rest of the section, we let $G$ be a~finite group.
Then, $\h^3(G,\mathbb{Z})=0$, so every bundle gerbe over $G$ is trivial.
Let $(\mathscr{K},\theta)$ be a~multiplicative structure on the trivial bundle gerbe over $G$. As an isomorphism between trivial bundle gerbes, we may identify $\mathscr{K}$ with a~line bundle over $G \times G$, and $\theta$ with a~line bundle isomorphism over $G^3$.
Using $(\mathscr{K},\theta)$, the category $\VectBdl(G)$ of complex vector bundles over $G$ can be equipped with a~monoidal structure by twisted convolution: if $\mathcal{W}$, $\mathcal{W}'$ are vector bundles, then
\[
 (\mathcal{W} \otimes \mathcal{W}')_g := \bigoplus_{hh'=g} \mathscr{K}_{h,h'} \otimes \mathcal{W}_h \otimes \mathcal{W}'_{h'},
\]
and the associator is defined using $\theta$ in a~straightforward way.
This yields a~($\C$-linear) tensor category that we denote by $\VectBdl(G)^{\otimes (\mathscr{K},\theta)}$.
The Drinfel'd centre \smash{$\mathscr{Z}\bigl(\VectBdl(G)^{\otimes (\mathscr{K},\theta)}\bigr)$} is a~modular tensor category; via the Reshetikhin--Turaev construction it
 determines an extended topological quantum field theory, with the given modular tensor category as its value on the circle.

 The category
 \smash{$\mathscr{Z}\bigl(\VectBdl(G)^{\otimes (\mathscr{K},\theta)}\bigr)$} can also be described in terms of the $G$-equivariant structure $(\mathscr{P}_{\mathscr{K}},\psi_{\theta})$ obtained from $(\mathscr{K},\theta)$, see \cite[Proposition 4.9]{freed6}. To recall this, we view $\mathscr{P}_{\mathscr{K}}$ as a~line bundle over $G \times G$, and $\psi_{\theta}$ as a~line bundle isomorphism over $G^3$.
 Then, the objects in $\mathscr{Z}\bigl(\VectBdl(G)^{\otimes (\mathscr{K},\theta)}\bigr)$ are vector bundles $\mathcal{W}$ over $G$ together with an iso\-morphism
 $\varphi \colon \mathcal{W}_1 \to \mathcal{W}_{2\cdot 1\cdot 2^{-1}} \otimes (\mathscr{P}_{\mathscr{K}})_{1,2}$ over $G \times G$ satisfying a~compatibility condition with $\psi_{\theta}$, cf.\ \cite[equation~4.9]{freed6}.
 Put differently, $\mathcal{W} \colon \mathscr{I} \to \mathscr{I}$ is a~1-morphism of trivial bundle gerbes and $\varphi$
 is a~2-isomor\-phism~${s^*\mathcal{W} \!\Rightarrow t^*\mathcal{W} \circ \mathscr{P}_{\mathscr{K}}}$ and hence the pair $(\mathcal{W},\varphi)$ is precisely an object in
 $\VectBdl^{(\mathscr{I},\mathscr{P}_{\mathscr{K}},\psi_{\theta})}(\act GG,1)$ in the sense of Definition~\ref{Definition of our twisted vector bundles}. A similar construction applies to the morphisms, and yields the following result.

\begin{Proposition}
There is an equivalence of categories
 \begin{equation*}
 \mathscr{Z}\bigl(\VectBdl(G)^{\otimes (\mathscr{K},\theta)}\bigr) \cong \VectBdl^{(\mathscr{I},\mathscr{P}_{\mathscr{K}},\psi_{\theta})}(\act GG)
 \end{equation*}
 between the Drinfel'd centre $\mathscr{Z}\bigl(\VectBdl(G)^{\otimes (\mathscr{K},\theta)}\bigr)$ and the category of $(\mathscr{I},\mathscr{P}_{\mathscr{K}},\psi_{\theta})$-twisted
 vector bundles over $(\act GG)$.
\end{Proposition}

\subsection{Freed's twisted invertible algebra bundles}
\label{Freed's invertible algebra bundle twistings}

We show in Theorem~\ref{freeds-invertible-algebra-bundles-comparison} below that Freed's twisted invertible algebra bundles~\cite{Freed2012a}, for a~subclass of their involved double covers, are the same as graded-equivariant super 2-line bundles.
The following definition is a~smooth version of \cite[Definition~1.95]{Freed2012a}.

\begin{Definition} \label{Definition: double cover}
 A \emph{double cover} of a~Lie groupoid $\Gamma$ is a~triple $\omega=\bigl(\tilde\Gamma,\tilde\phi,\pi\bigr)$ consisting of a~graded Lie groupoid \smash{$\bigl(\tilde \Gamma,\tilde \phi\bigr)$} and a~weak equivalence \smash{$\pi\colon \tilde \Gamma \to \Gamma$}. \end{Definition}

Next we recall the definitions of Freed's twisted invertible algebra bundles following {\cite[Definitions~1.59 and~1.99]{Freed2012a}}.

\begin{Definition}\label{Invertible algebra bundle objects}
 Let $\omega=\bigl(\tilde\Gamma,\tilde\phi,\pi\bigr)$ be a~double cover of a~Lie groupoid $\Gamma$. An \emph{$\omega$-twisted invertible algebra bundle over $\Gamma$} is a~triple $\mathscr{A}=(\mathcal{A},\mathcal{B},\lambda)$
 consisting of the following data:
 \begin{itemize}\itemsep=0pt
\item
 a~central simple super algebra bundle $\mathcal{A}$ over $\tilde \Gamma_0$,

\item
 an invertible super bimodule bundle $\mathcal{B}$ over $\tilde \Gamma_1$
 whose fibre over $\gamma \in \tilde \Gamma_1$ is an \smash{$\mathcal{A}_{t(\gamma)}^{\tilde\phi(\gamma)}$}\text{-}$\mathcal{A}_{s(\gamma)}$-bimodule, and

\item
 an invertible even intertwiner $\lambda$ over $\tilde \Gamma_2$ that is fibre-wise over $(\gamma_1,\gamma_2) \in \tilde \Gamma_2$ an intertwiner of \smash{$\mathcal{A}_{t(\gamma_1)}^{\tilde\phi(\gamma_1 \circ \gamma_2)}$}\text{-}$\mathcal{A}_{s(\gamma_2)}$ bimodules
 \[
 \lambda_{\gamma_1,\gamma_2} \colon\ \mathcal{B}_{\gamma_1}^{\tilde\phi(\gamma_2)} \otimes_{\mathcal{A}_{t(\gamma_2)}^{\tilde\phi(\gamma_2)}} \mathcal{B}_{\gamma_2} \to \mathcal{B}_{\gamma_1 \circ \gamma_2},
 \]
\end{itemize}
 such that the following diagram of intertwiners
is commutative for all $(\gamma_1,\gamma_2,\gamma_3) \in \tilde\Gamma_3$
 \begin{equation*}
\begin{gathered}
 \begin{xy}
 \xymatrix@C=-1.5em{
 \bigl(\mathcal{B}_{\gamma_1}^{\tilde\phi(\gamma_2)} \otimes_{\mathcal{A}_{t(\gamma_2)}^{\tilde\phi(\gamma_2)}} \mathcal{B}_{\gamma_2}\bigr)^{\tilde\phi(\gamma_3)} \otimes_{\mathcal{A}_{t(\gamma_3)}^{\tilde\phi(\gamma_3)}} \mathcal{B}_{\gamma_3} \ar[rr] \ar[d]_-{\lambda_{\gamma_1,\gamma_2}^{\tilde\phi(\gamma_3)} \otimes \id} && \mathcal{B}_{\gamma_1}^{\tilde\phi(\gamma_2) \cdot \tilde\phi(\gamma_3)} \otimes_{\mathcal{A}_{t(\gamma_2 \circ \gamma_3)}^{\tilde\phi(\gamma_2 \circ \gamma_3)}} \bigl(\mathcal{B}_{\gamma_2}^{\tilde\phi(\gamma_3)} \otimes_{\mathcal{A}_{t(\gamma_3)}^{\tilde\phi(\gamma_3)}} \mathcal{B}_{\gamma_3}\bigr) \ar[d]^-{\id \otimes \lambda_{\gamma_2,\gamma_3}}\\
 \mathcal{B}_{\gamma_1\circ \gamma_2}^{\tilde\phi(\gamma_3)} \otimes_{\mathcal{A}_{t(\gamma_3)}^{\tilde\phi(\gamma_3)}} \mathcal{B}_{\gamma_3} \ar[dr]_{\lambda_{\gamma_1\circ \gamma_2,\gamma_3}} && \mathcal{B}_{\gamma_1}^{\tilde\phi(\gamma_2 \circ \gamma_3) } \otimes_{\mathcal{A}_{t(\gamma_2 \circ \gamma_3)}^{\tilde\phi(\gamma_2 \circ \gamma_3)}} \mathcal{B}_{\gamma_2\circ \gamma_3} \ar[dl]^{\lambda_{\gamma_1,\gamma_2\circ \gamma_3}} \\
 & \mathcal{B}_{\gamma_1\circ \gamma_2\circ \gamma_3}.
 }
 \end{xy}
\end{gathered}
 \end{equation*}
 \end{Definition}

Here, we have employed our notation~\eqref{complex conjugation according to sign} saying that $(..)^{\tilde\phi}$ means complex conjugation when~${\tilde\phi=-1}$ and no operation otherwise.
For the sake of readability, we will omit writing out the algebra in relative tensor products in the following. We remark that -- despite of the terminology -- a~twisted invertible algebra bundle over $\Gamma$ makes no contact with $\Gamma$ or the functor~${\pi\colon \tilde\Gamma \to \Gamma}$. Next we recall the morphisms, following \cite[Definition~1.62\,(i)]{Freed2012a}.

\begin{Definition}\label{Invertible algebra bundle 1-morphisms}
 Let $\omega=\bigl(\tilde\Gamma,\tilde\phi,\pi\bigr)$ be a~double cover of $\Gamma$, and let $\mathscr{A}=(\mathcal{A},\mathcal{B},\lambda)$ and $\mathscr{A}'=(\mathcal{A}',\mathcal{B}',\lambda')$ be $\omega$-twisted invertible algebra bundles over $\Gamma$.
 A \emph{$1$-morphism} $\mathscr{A} \to \mathscr{A}'$ is a~pair $(\mathcal{C},\mu)$ consisting of an invertible super bimodule bundle $\mathcal{C}$ over $\tilde\Gamma_0$, whose fibre over $x \in \tilde \Gamma_0$
 is an $\mathcal{A}'_x$\text{-}$\mathcal{A}_x$-bimodule, and of an invertible even intertwiner $\mu$ that is fibrewise over $\gamma \in \tilde\Gamma_1$ an intertwiner of \smash{${\mathcal{A}'}_{t(\gamma)}^{\tilde\phi(\gamma)}$}\text{-}$\mathcal{A}_{s(\gamma)}$-bimodules
\smash{$
 \mu_\gamma\colon \mathcal{C}_{t(\gamma)}^{\tilde\phi(\gamma)} \otimes \mathcal{B}_{\gamma} \to \mathcal{B}'_{\gamma} \otimes \mathcal{C}_{s(\gamma)}
$}
 such that the following diagram of intertwiners is commutative for all $(\gamma_1,\gamma_2) \in \tilde\Gamma_2$
 \begin{equation}
 \begin{gathered}
 \begin{xy}
 \xymatrix{
 \mathcal{C}_{t(\gamma_1 \circ \gamma_2)}^{\tilde\phi(\gamma_1 \circ \gamma_2)} \otimes \mathcal{B}_{\gamma_1}^{\tilde\phi(\gamma_2)} \otimes \mathcal{B}_{\gamma_2} \ar[r]^-{\mu_{\gamma_1}^{\tilde\phi(\gamma_2)} \otimes \id} \ar[d]_{\id \otimes \lambda_{\gamma_1,\gamma_2}} &
 {\mathcal{B}'}_{\gamma_1}^{\tilde\phi(\gamma_2)} \otimes \mathcal{C}_{t(\gamma_2)}^{\tilde\phi(\gamma_2)} \otimes \mathcal{B}_{\gamma_2} \ar[r]^-{\id \otimes \mu_{\gamma_2}} &
 {\mathcal{B}'}_{\gamma_1}^{\tilde\phi(\gamma_2)} \otimes \mathcal{B}'_{\gamma_2} \otimes \mathcal{C}_{s(\gamma_2)} \ar[d]^{\lambda'_{\gamma_1,\gamma_2} \otimes \id}\\
 \mathcal{C}_{t(\gamma_1 \circ \gamma_2)}^{\tilde\phi(\gamma_1 \circ \gamma_2)} \otimes \mathcal{B}_{\gamma_1 \circ \gamma_2}\ar[rr]_-{\mu_{\gamma_1 \circ \gamma_2}} && \mathcal{B}'_{\gamma_1 \circ \gamma_2} \otimes \mathcal{C}_{s(\gamma_1 \circ \gamma_2)}.
 }
 \end{xy}
 \end{gathered}
 \end{equation}
 \end{Definition}

We remark that morphisms are only defined between twisted invertible algebra bundles with the same double cover $\omega$. Finally, we recall the definition of 2-morphisms from \cite[Definition~1.62\,(ii)]{Freed2012a}.

\begin{Definition}\label{Invertible algebra bundle 2-morphisms}
 Let $\bigl(\mathcal{C}^1,\mu^1\bigr)$ and $\bigl(\mathcal{C}^2,\mu^2\bigr)$ be 1-morphisms between $\omega$-twisted invertible algebras~$\mathscr{A} $ and~$ \mathscr{A}'$. A \emph{$2$-morphism}
 $\bigl(\mathcal{C}^1,\mu^1\bigr) \Rightarrow \bigl(\mathcal{C}^2,\mu^2\bigr)$ is an invertible even intertwiner $\nu$ of $\mathcal{A}'_x$-$\mathcal{A}_x$-bimodules,
$
 \nu_x \colon \mathcal{C}^1_x \to \mathcal{C}^2_x
$
 such that the following diagram of intertwiners is commutative for all $\gamma \in \tilde \Gamma_1$
 \[
 \begin{gathered}
 \begin{xy}
 \xymatrix{
 {\mathcal{C}^1}_{t(\gamma)}^{\tilde\phi(\gamma)} \otimes \mathcal{B}_{\gamma} \ar[r]^-{\mu^1_\gamma} \ar[d]_{\nu_{t(\gamma)}^{\tilde\phi(\gamma)} \otimes \id} &
 \mathcal{B}'_{\gamma} \otimes {\mathcal{C}^1}_{s(\gamma)} \ar[d]^{\id \otimes \nu_{s(\gamma)}}\\
 {\mathcal{C}^2}_{t(\gamma)}^{\tilde\phi(\gamma)} \otimes \mathcal{B}_{\gamma} \ar[r]_-{\mu^2_\gamma} & \mathcal{B}'_{\gamma} \otimes {\mathcal{C}^2}_{s(\gamma)}.
 }
 \end{xy}
 \end{gathered}
 \]
\end{Definition}

Freed remarks that $\omega$-twisted invertible algebra bundles over $\Gamma$ form a~bigroupoid, which we denote by $\mathrm{Fr}(\Gamma)^{\omega}$ in order to reflect the terminology appropriately. Let us first note the following result, which translates one-to-one between Freed's and our setting.

 \begin{Lemma}
 \label{Main theorem of the twisted invertible algebra bundle section}
 Let $\omega=\bigl(\tilde\Gamma,\tilde\phi,\pi\bigr)$ be a~double cover of a~Lie groupoid $\Gamma$. There is a~canonical equivalence of bigroupoids
\smash{$
 \mathrm{Fr}(\Gamma)^{\omega} \cong \cssAlgBdlbi\bigl(\tilde \Gamma,\tilde \phi\bigr)^\times
$}
 between Freed's $\omega$-twisted invertible algebra bundles and the maximal sub-bigroupoid of the bicategory of central simple super algebra bundles over $\bigl(\tilde\Gamma,\tilde\phi\bigr)$.
\end{Lemma}

\begin{proof}
 The equivalence is a~straightforward comparison of the Definitions~\ref{Invertible algebra bundle objects},
 \ref{Invertible algebra bundle 1-morphisms}, and~\ref{Invertible algebra bundle 2-morphisms} with Definitions~\ref{equivariant structure}, \ref{equivariant 1-morphisms}, and~\ref{equivariant 2-morphisms}. We only note that, on the level of 1-morphisms,
 the comparison inverts the map $\mu$ in Definition~\ref{Invertible algebra bundle 1-morphisms} to get the map $\eta$
 from Definition~\ref{equivariant 1-morphisms}, and vice versa.
\end{proof}

The equivalence of Lemma~\ref{Main theorem of the twisted invertible algebra bundle section} does not properly reflect the meaning of twisted invertible algebra bundles, as they are meant to represent twistings on $\Gamma$ and not twistings on $\tilde\Gamma$.
In order to resolve this, we will restrict to double covers $\omega=\bigl(\tilde\Gamma,\tilde\phi,\pi\bigr)$, where $\tilde\phi=\phi \circ \pi$, for a~grading $\phi$ on $\Gamma$.
We will call double covers of this kind \emph{$\phi$-induced double covers}.

We assume that $(\Gamma,\phi)$ is a~graded Lie groupoid and $\omega=\bigl(\tilde\Gamma,\tilde\phi,\pi\bigr)$ is a~$\phi$-induced double cover of $\Gamma$. We are then in position to construct a~2-functor
\begin{equation}
\label{descent-for-twisted-algebra-bundles}
\mathrm{Fr}(\Gamma)^{\omega} \cong \cssAlgBdlbi\bigl(\tilde \Gamma,\tilde \phi\bigr)^\times \incl \stwoLineBdl\bigl(\tilde \Gamma,\tilde \phi\bigr) \cong \stwoLineBdl(\Gamma,\phi),
\end{equation}
where the first equivalence is Lemma~\ref{Main theorem of the twisted invertible algebra bundle section}, the arrow is the inclusion of central simple super algebra bundles into super 2-line bundles, and the second equivalence is Theorem~\ref{descent-for-graded-equivariant-2-line-bundles}.
Since the right-hand side in~\eqref{descent-for-twisted-algebra-bundles} is independent of the double cover, the 2-functor~\eqref{descent-for-twisted-algebra-bundles} extends to the disjoint union over all $\phi$-induced double covers of $(\Gamma,\phi)$.

\begin{Theorem}
\label{freeds-invertible-algebra-bundles-comparison}
Let $(\Gamma,\phi)$ be a~graded Lie groupoid. Then, the $2$-functor
\begin{equation*}
\coprod_{\omega}\mathrm{Fr}(\Gamma)^{\omega}\to \stwoLineBdl(\Gamma,\phi),
\end{equation*}
where the disjoint union goes over all $\phi$-induced double covers of $(\Gamma,\phi)$,
is essentially surjective.
\end{Theorem}

\begin{proof}
Let $(\mathscr{L},\mathscr{P},\psi)$ be a~super 2-line bundle over $(\Gamma,\phi)$. Then, $\mathscr{L}$ has some surjective submersion $\pi\colon Y \to \Gamma_0$. The corresponding covering groupoid $\Gamma^{\pi}$ (see Definition~\ref{Example : covering groupoid and weak equivalences}) has a~natural grading $\tilde\phi := \phi \circ \pi$, where $\pi\colon \Gamma^{\pi} \to \Gamma$, and $\pi\colon \Gamma^{\pi} \to \Gamma$ is a~weak equivalence. Thus, $\omega := \bigl(\Gamma^{\pi},\pi,\tilde\phi\bigr)$ is a~$\phi$-induced double cover. The pullback $\pi^{*}\mathscr{L}$ is canonically isomorphic to a~central simple super algebra bundle $\mathcal{A}$. Transferring the equivariant structure along this isomorphism produces a~central simple algebra bundle $(\mathcal{A},\mathscr{P}',\psi')$ over $\big(\Gamma^{\pi},\tilde\phi\big)$; which, via Lemma~\ref{Main theorem of the twisted invertible algebra bundle section}, is a~$\omega$-twisted invertible algebra bundle. Since the equivalence of Theorem~\ref{descent-for-graded-equivariant-2-line-bundles}, which is used in the definition of~\eqref{descent-for-twisted-algebra-bundles}, is established by pullback along the weak equivalence $\pi$, and $\pi^{*}\mathscr{L}\cong \mathcal{A}$, it follows that~${(\mathcal{A},\mathscr{P}',\psi')}$ is indeed an essential preimage of $(\mathscr{L},\mathscr{P},\psi)$.
\end{proof}

We expect that, if one enlarges the definitions above by allowing for non-invertible morphisms and morphisms between twisted invertible algebra bundles for different double covers, then
Theorem~\ref{freeds-invertible-algebra-bundles-comparison} lifts to an equivalence of bicategories.

\subsection{Distler--Freed--Moore twistings} \label{Distler Freed Moore twistings}

In this section, we want to relate the twistings introduced by Distler, Freed, and Moore to
graded-equivariant super 2-line bundles~\cite{Distler2011,distler2010spin}; we will use a~version described in the slides~\cite{Freed2009Slides}. The following is a~smooth version
of that definition.

\begin{Definition} \label{definition: DFM-twisting}
 Let $\Gamma$ be a~Lie groupoid and $\omega=\bigl(\tilde\Gamma,\tilde\phi,\pi\bigr)$ a~double cover of $\Gamma$ in the sense of Definition~\ref{Definition: double cover}.
 A \emph{DFM $\omega$-twisting} of $\Gamma$ is a~triple $\tau = (d,\mathcal{L},\lambda)$ consisting of a~locally constant function
 $d \colon \tilde\Gamma_0 \to \mathbb{Z}$, a~super line bundle $\mathcal{L}$ over $\tilde\Gamma_1$,
 and an isomorphism of super line bundles over $\tilde\Gamma_2$, which fibre-wise over $(g,f) \in \tilde\Gamma_2$ reads
\smash{$
 \lambda_{g,f} \colon \mathcal{L}_g^{\tilde\phi(f)} \otimes \mathcal{L}_f \to \mathcal{L}_{g\circ f}$}.
 The map $d$ is subject to the condition that $d \circ s=d\circ t$, i.e.,
 $d(s(f)) = d(t(f))$ for all $f \in \tilde\Gamma_1$, and the isomorphism $\lambda$ is required
 to be associative over $\tilde\Gamma_3$.
\end{Definition}

We add the following notion of morphisms between DFM $\omega$-twistings.
 A \emph{morphism} of DFM $\omega$-twistings
 $\tau \to \tau'$ is a~morphism of $\tilde \phi$-twisted super extensions of $\tilde\Gamma$ in the sense of Definition~\ref{The category of twisted super extensions}, and it is required that the maps $d$ and $d'$ coincide.
We denote the category of DFM $\omega$-twistings of $\Gamma$ by $\mathrm{DFM}^\omega(\Gamma)$.
By comparison of definitions, we find the following.

\begin{Lemma} \label{lemma: DFM twistings via FM category}
 Let $\Gamma$ be a~Lie groupoid and $\omega=\bigl(\tilde\Gamma,\tilde\phi,\pi\bigr)$ a~double cover of $\Gamma$.
 There is an equivalence of categories
\smash{$
 \mathrm{DFM}(\Gamma)^\omega \cong C^{\infty}\bigl(\tilde \Gamma,\Z\bigr)_{\rm dis} \times \mathrm{FM}\bigl(\tilde\Gamma,\tilde\phi\bigr)$},
where \smash{$C^{\infty}\bigl(\tilde \Gamma,\Z\bigr)$} denotes the set of smooth maps \smash{$d \colon \tilde\Gamma_0 \to \Z$} such that $d \circ s=d\circ t$.
\end{Lemma}

Next we exhibit the relation to super 2-line bundles, and consider first the functor
\begin{equation*}
C^{\infty}\bigl(\tilde \Gamma,\Z\bigr)_{\rm dis} \to\cssAlgBdlbi\bigl(\tilde\Gamma,\tilde\phi\bigr)^{\times}
\end{equation*}
defined in the following way. The function $d$ will select a~Morita equivalence class of central simple super algebras, i.e., an element in the Brauer--Wall group $\mathrm{BW}_{\C}\cong \Z_2=\{0,1\}$. We shall fix representatives: $A_0:=\C$ for the neutral element, and $A_{1}:=\C \oplus \C u$ (where $u$ is odd and~${u^2=1}$) for the non-trivial element. Then, we assign to $d$ the central simple super algebra bundle $\mathcal{A}^{d} := \tilde\Gamma_0 \times A_{d \text{ mod } 2}$.
In order to equip $\mathcal{A}^{d}$ with a~$\bigl(\tilde\Gamma,\tilde\phi\bigr)$-equivariant structure $(\mathscr{P},\psi)$ in the sense of Definition~\ref{equivariant structure}, we note that there are canonical super algebra isomorphisms $A_0 \cong \overline{A_0}$ (just complex conjugation) and $A_1 \cong \overline{A_1}$, given by $x\oplus yu \mapsto \overline{x}\oplus\overline{y}u$. Thus, we have a~bundle isomorphism $\rho\colon \mathcal{A}^{d} \to \overline{\mathcal{A}^{d}}$. The required 1-isomorphism $\mathscr{P}\colon s^*\mathcal{A}^{d} \to \bigl(t^*\mathcal{A}^{d}\bigr)^{\phi}$ is defined to be the identity on connected components with $\phi=1$, and to be induced by the isomorphism $\rho$ over connected components with $\phi=-1$, under the assignment of invertible bimodules to algebra isomorphisms via the framing. The 2-isomorphism $\psi$ is the identity (using that $\rho^2=\id$).

Second, we consider the functor
\smash{$
\mathrm{FM}\bigl(\tilde\Gamma,\tilde\phi\bigr) \to \stwoLineBdltriv^{\text{ref}}\bigl(\tilde\Gamma,\tilde\phi\bigr)
\subset\cssAlgBdlbi\bigl(\tilde\Gamma,\tilde\phi\bigr)^{\times}
$}
of
Theorem~\ref{The main theorem of the Freed-Moore chapter}. Together with the symmetric monoidal structure on super 2-line bundles, see also Remark~\ref{Remark: the symmetric monoidal structure extends to the graded-equivariant case}, we obtain a~functor
\begin{align}
 C^{\infty}\bigl(\tilde \Gamma,\Z\bigr)_{\rm dis} \times \mathrm{FM}\bigl(\tilde\Gamma,\tilde\phi\bigr) &\to \cssAlgBdlbi\bigl(\tilde\Gamma,\tilde\phi\bigr)^{\times}\times \cssAlgBdlbi\bigl(\tilde\Gamma,\tilde\phi\bigr)^{\times}\nonumber \\
 &\stackrel\otimes\to \cssAlgBdlbi\bigl(\tilde\Gamma,\tilde\phi\bigr)^{\times} .\label{DFM-twistings-to-algebra-bundles}
\end{align}

We assume now that $(\Gamma,\phi)$ is a~graded Lie groupoid and
$\omega=\bigl(\tilde\Gamma,\tilde\phi,\pi\bigr)$ is a~$\phi$-induced double cover of $\Gamma$, see Section~\ref{Freed's invertible algebra bundle twistings}. Combining Lemma~\ref{lemma: DFM twistings via FM category}, Theorem~\ref{descent-for-graded-equivariant-2-line-bundles}, and
the functor~\eqref{DFM-twistings-to-algebra-bundles}, we obtain a~functor
$\mathrm{DFM}(\Gamma)^\omega
\to
\stwoLineBdl(\Gamma,\phi)$.

\begin{Theorem}
\label{freed-distler-moore-comparison}
The functor
\begin{equation*}
\coprod_{\omega}\mathrm{DFM}(\Gamma)^\omega\to \stwoLineBdl(\Gamma,\phi),
\end{equation*}
where the disjoint union goes over all $\phi$-induced double covers of $(\Gamma,\phi)$, is essentially surjective.
\end{Theorem}

\begin{proof}
Let $(\mathscr{L},\mathscr{P},\psi)$ be a~super 2-line bundle over $(\Gamma,\phi)$. We first extract the map $d$. On each connected component $c$ of $\Gamma_0$, $\mathscr{L}$ has some Morita class $d_c \in \mathrm{BW}_{\C}\cong \Z_2$ in the Brauer--Wall group of $\C$, see \cite[Section~3.1]{Kristel2020}.
We regard this as a~smooth map $d \colon \Gamma_0 \to \Z_2 \subset \Z$. The isomorphism~${\mathscr{P}\colon s^{*}\mathscr{L} \to t^{*}\mathscr{L}^{\phi}}$ implies that $d(s(\gamma))= d(t(\gamma))$. 
We let $\mathcal{A}^{-d}$ be the central simple algebra bundle over $(\Gamma,\phi)$ that corresponds to $-d$ via the construction above.
Then, $\mathscr{L}':=\mathscr{L} \otimes \mathcal{A}^{-d}$ is a~super 2-line bundle over $(\Gamma,\phi)$ with trivial Morita class, i.e., a~super bundle gerbe \cite[Proposition~3.2.4]{Kristel2020}. As in the proof of Theorem~\ref{freeds-invertible-algebra-bundles-comparison}, we then construct a~$\phi$-induced double cover $\omega=\bigl(\Gamma^{\pi},\pi,\tilde\phi\bigr)$ from the surjective submersion of $\mathscr{L}'$, so that the pullback $\pi^{*}\mathscr{L}'$ is trivializable. By Theorem~\ref{The main theorem of the Freed-Moore chapter}, it can then be identified with a~\smash{$\tilde\phi$}-twisted super extension of~$\Gamma^{\pi}$.
\end{proof}

\begin{Remark}
Our functor~\eqref{DFM-twistings-to-algebra-bundles} above may be used, together with Lemmas~\ref{lemma: DFM twistings via FM category} and~\ref{Main theorem of the twisted invertible algebra bundle section}, to obtain a~functor
$
\mathrm{DFM}(\Gamma)^\omega \to \mathrm{Fr}(\Gamma)^{\omega}
$
that relates Distler--Freed--Moore's twistings directly to Freed's invertible algebra bundles.
\end{Remark}

\appendix

\section{Appendix: Lie groupoids and weak equivalences}
\label{appendix}
\label{bicategory-of-lie-groupoids}

Lie groupoids, smooth functors, and smooth natural transformations form a~strict 2-category $\liegrpd^{\rm fun}$. The 1-isomorphisms in this 2-category are called \emph{strong equivalences}. There is also a~bicategory $\liegrpd$ whose objects are Lie groupoids, whose 1-morphisms are smooth anafunctors, and whose 2-morphisms are smooth anafunctor transformations, see, for example, \cite{DuLi2014,Nikolaus,pries2}.
This bicategory $\liegrpd$ is equivalent to the bicategory of differentiable stacks~\cite{pronk}. Smooth anafunctors are also known as bibundles, Hilsum--Skandalis maps, or generalized morphisms; the precise notion is not of relevance for this article. Invertible smooth anafunctors establish Morita equivalences between Lie groupoids.

Every smooth functor induces a~smooth anafunctor, this way making up a~functor
$
\liegrpd^{\rm fun}\allowbreak\to \liegrpd
$
that is the identity on the level of objects, and fully faithful on 2-morphism level.
A smooth functor is called a~\emph{weak equivalence} if and only if its induced smooth anafunctor is invertible, see~\cite[Definition~58/Proposition~60]{metzler} or \cite[Lemma~3.34]{lerman1}, and hence a~Morita equivalence.
The following definition characterizes them precisely.
\begin{Definition}\label{Definition: weak equivalence of Lie groupoids}
 A smooth functor $F\colon \Gamma \to \Lambda$ between Lie groupoids is called a~\emph{weak equivalence} if
 both of the following conditions are satisfied:
{\samepage \begin{itemize}\itemsep=0pt
 \item $F$ is smoothly essentially surjective: the map
$
 s \circ \pr_{\Lambda_1} \colon \Gamma_0 \ttimes{F_0}{t} \Lambda_1 \to \Lambda_0
$
 is a~surjective submersion.
 \item $F$ is smoothly fully faithful: the commutative square
 \[
 \alxydim{}{\Gamma_1 \ar[rr]^{F_1} \ar[d]_{(s,t)} & & \Lambda_1 \ar[d]^{(s,t)} \\
 \Gamma_0 \times \Gamma_0 \ar[rr]_{F_0 \times F_0} & & \Lambda_0 \times \Lambda_0}
 \]
 is a~pullback diagram in the category of smooth manifolds.
 \end{itemize}}
\end{Definition}

\begin{Remark}
\label{Lemma: facts about weak equivalences}
The following facts follow immediately from the fact that weak equivalences become 1-isomorphisms in a~bicategory:
 \begin{enumerate}[(i)]\itemsep=0pt
 \item If $T,T'\colon \Lambda \to \Gamma$ are smooth functors, and there is a~smooth natural transformation $T \Rightarrow T'$, then
 $T$ is a~weak equivalence if and only if $T'$ is a~weak equivalence.
 \item\label{2-out-of-3}
 Weak equivalences satisfy the 2-out-of-3 property: if $T \colon \Omega \to \Lambda$ and $F\colon \Lambda \to \Gamma$ are smooth functors,
 and any two of the three smooth functors $T$, $F$, and $F \circ T$ are weak equivalences, then so is the third.
 \item \label{Weak equivalence filling} If $F,F' \colon \Omega \to \Lambda$ are smooth functors, $T \colon \Lambda \to \Gamma$ is a~weak equivalence, then there is a~smooth natural transformation
 $T \circ F \Rightarrow T \circ F'$ if and only if there is a~smooth natural transformation $F \Rightarrow F'$.
 \end{enumerate}
\end{Remark}

\begin{Proposition}[{\cite[Theorem 2.16]{nikolaus2}}] \label{Nikolaus-Schweigert weak equivalences}
 Suppose $F\colon \Gamma \to \Lambda$ is a~weak equivalence between Lie groupoids.
 Let $\mathfrak{X}$ be a~$2$-stack on the site of smooth manifolds. Then, the functor
$
F^{*}\colon \mathfrak{X}(\Lambda) \to \mathfrak{X}(\Gamma)
$
 is an equivalence of bicategories.
\end{Proposition}

Next we discuss the notion of a~covering groupoid and provide some results we use in the main text.

\begin{Definition}
\label{Example : covering groupoid and weak equivalences}
 Let $\Gamma$ be a~Lie groupoid and $\pi\colon Y \to \Gamma_0$ be a~surjective submersion. The
 \emph{covering groupoid} $\Gamma^\pi$ has objects $\Gamma^\pi_0 := Y$ and morphisms $\Gamma^\pi_1 := Y \ttimes\pi t \Gamma_1 \ttimes s\pi Y$, and its structural maps are defined by
\begin{equation*}
s(y_2,\gamma,y_1):= y_1
,\qquad
t(y_2,\gamma,y_1) := y_2
,\qquad
(y_3,\gamma_{23},y_2) \circ (y_2,\gamma_{12},y_1) := (y_3,\gamma_{23} \circ \gamma_{12},y_1).
\end{equation*}
The \emph{projection functor} $P^{\pi}\colon \Gamma^{\pi} \to \Gamma$ is given by $P_0 :=\pi$ and $P_1(y_2,\gamma,y_1) := \gamma$.
\end{Definition}

\begin{Example}\label{Example : Cech groupoid and weak equivalences}
 Let $M$ be a~smooth manifold and $\pi\colon Y \to M$ a~surjective submersion. Then, $\text{\v C}(\pi):=M_{\rm dis}^{\pi}$ is the usual \v Cech groupoid.
\end{Example}

\begin{Lemma}
Let $\Gamma$ be a~Lie groupoid. If $\pi\colon Y \to \Gamma_0$ is a~surjective submersion, then the projection functor $P^{\pi}\colon \Gamma^{\pi} \to \Gamma$ is a~weak equivalence. Moreover, if $\zeta\colon Z \to Y$ is a~further surjective submersion, then the evident functor $\Gamma^{\pi \circ \zeta} \to \Gamma^{\pi}$ is a~weak equivalence.
\end{Lemma}

\begin{proof}
The first statement can be checked directly. The second statement follows from Remark~\ref{Lemma: facts about weak equivalences}\,\ref{2-out-of-3}.
\end{proof}

\begin{Lemma}
\label{pullback-of-coverings}
Let $\Lambda\colon \Gamma \to \Omega$ be smooth functor, and let $\pi\colon Y \to \Omega_0$ be a~surjective submersion. Then, there exists a~surjective submersion $\zeta\colon Z \to \Gamma_0$ and a~functor $\tilde\Lambda\colon \Gamma^{\zeta} \to \Omega^{\pi}$ such that the diagram
\begin{equation*}
\alxydim{}{\Gamma^{\zeta} \ar[d]_{P^{\zeta}} \ar[r]^{\tilde\Lambda} & \Omega^{\pi} \ar[d]^{P^{\pi}} \\ \Gamma \ar[r]_{\Lambda} & \Omega}
\end{equation*}
is strictly commutative. Moreover, if $\Lambda$ is a~weak equivalence, then $\tilde\Lambda$ is a~weak equivalence.
\end{Lemma}
\begin{proof}
The surjective submersion $\zeta$ is the pullback of $\pi$ along $\Lambda_0$, i.e., we have a~pullback diagram
\begin{equation*}
\alxydim{}{Z \ar[d]_{\zeta} \ar[r]^{\tilde\Lambda_0} & Y \ar[d]^{\pi} \\ \Gamma_0\ar[r]_{\Lambda_0} & \Omega_0.}
\end{equation*}
This defines $\tilde\Lambda$ on the level of objects, and it is straightforward to complete its definition.
\end{proof}

Coverings groupoids are important because weak equivalences can be inverted on covering groupoids; see, e.g., \cite[Proposition 5.7]{nikolaus2}.
\begin{Proposition} \label{Lemma: Extracting covering groupoids from weak equivalences}
 A smooth functor $T \colon \Lambda \to \Gamma$ is a~weak equivalence if and only if there exists a~surjective submersion $\pi \colon Y \to \Gamma_0$ and a~strong equivalence
 $S \colon \Lambda \to \Gamma^\pi$ such that $T \cong P^{\pi} \circ S$.
\end{Proposition}

Finally, we discuss fibre products of Lie groupoids.

\begin{Definition}
\label{definition: fibre product of Lie groupoids}
 Let $T\colon \Lambda \to \Gamma$ and $T' \colon \Lambda' \to \Gamma$ be smooth functors of Lie groupoids for which either
 \smash{$t \circ \pr_{\Gamma_1} \colon \Lambda_0 \ttimes{T_0}{s} \Gamma_1 \to \Gamma_0$} or \smash{$s \circ \pr_{\Gamma_1} \colon \Lambda_0' \ttimes{T_0'}{t} \Gamma_1 \to \Gamma_0$}
 is a~surjective submersion. The \emph{fibre product} $\Lambda \times_{\Gamma} \Lambda'$
 is the Lie groupoid whose objects are tuples $(p,\phi,p')$, where $p \in \Lambda$, $p' \in \Lambda'$ are objects and $\phi \colon T(p) \to T'(p')$ is an isomorphism.
 Morphisms $(p,\phi,p') \to (q,\psi,q')$ are pairs~${(\alpha,\beta)}$ where $\alpha \colon p \to q$ and $\beta\colon p' \to q'$ are morphisms such that
 the following diagram commutes:
 \begin{equation*}
 \xymatrix{
 T(p) \ar[r]^-{T(\alpha)} \ar[d]_-{\phi} & T(q) \ar[d]^-{\psi}\\
 T'(p') \ar[r]_-{T'(\beta)} & T'(q').
 }
 \end{equation*}
 Composition is given componentwise. Put differently, the smooth manifold of objects is $\Lambda_0 \ttimes{T}{s} \Gamma_1 \ttimes{t}{T'} \Lambda_0'$ and
 the smooth manifold of morphisms is $\Lambda_1 \ttimes{s \circ T}{s} \Gamma_1 \ttimes{t}{s \circ T'} \Lambda_1'$.

\end{Definition}
\begin{Remark}\label{Remark: the fibre product is a~2-pullback}
 The fibre product of Lie groupoids is a~2-pullback in the bicategory $\liegrpd^{\rm fun}$, see \cite[Section~5.3]{Moerdijk2003foliations} and
 \cite[Remark 1.18]{DuLi2014} for further references.
\end{Remark}

\begin{Remark} \label{Remark: facts about fibre products}
 The following elementary consequence of these definitions is \cite[Proposition~1.26]{DuLi2014}.
 For a~weak equivalence $T \colon \Lambda \to \Gamma$ and a~smooth functor $F \colon \Lambda' \to \Gamma$, their fibre product~${\Lambda \times_{\Gamma} \Lambda'}$ exists.
 Moreover, the projection $\Lambda \times_{\Gamma} \Lambda' \to \Lambda'$ is a~weak equivalence.
\end{Remark}

\begin{Lemma}\label{Lemma: weak equivalences induce weak equivalences of fibre products}
 Let $T \colon \Lambda \to \Gamma$, $T' \colon \Lambda' \to \Gamma$, $S \colon \Omega \to \Lambda$ and $S' \colon \Omega' \to \Lambda'$ weak equivalences of Lie groupoids.
 Then, there is a~weak equivalence
$
 S' \times_{\Gamma} S \colon \Omega' \times_{\Gamma} \Omega \to \Lambda' \times_{\Gamma} \Lambda$.
\end{Lemma}
\begin{proof}
 Combine Remarks~\ref{Remark: the fibre product is a~2-pullback}, \ref{Remark: facts about fibre products}, and~\ref{Lemma: facts about weak equivalences}.
\end{proof}

Fibre products of surjective submersions and fibre products of Lie groupoids are compatible in the following way.

\begin{Lemma}
 \label{Common refinements of common refinements induce weak equivalences NEU}
 Let $\Gamma$ be a~Lie groupoid and let
 $\pi \colon Y \to \Gamma_0$ and $\pi' \colon Y' \to \Gamma_0$ be surjective submersions. Let $Z:=Y \times_{\Gamma_0} Y'$ be the fibre product, and let $\zeta\colon Z \to \Gamma_0$ be the projection to the base. Then, there is a~canonical weak equivalence $\Gamma^{\zeta} \to \Gamma^{\pi} \times_{\Gamma} \Gamma^{\pi'}$ such that the diagram
 \begin{equation*}
 \alxydim{@C=2em}{\Gamma^{\zeta} \ar[dr]_{\zeta} \ar[rr] && \Gamma^{\pi} \times_{\Gamma} \Gamma^{\pi'} \ar[dl] \\ & \Gamma }
 \end{equation*}
 is strictly commutative.
 \end{Lemma}
 \begin{proof}
 We construct a~functor
 $F\colon \Gamma^{\zeta} \to \Gamma^{\pi} \times_{\Gamma} \Gamma^{\pi'}$. To an object $(y,y')\in Y \times_{\Gamma_0}Y'$ we assign the object $F(y,y') := (y,y',\id_{\pi(y)})$ of $\Gamma^{\pi} \times_{\Gamma} \Gamma^{\pi'}$. To a~morphism $(y_1,y_1',\gamma,y_2,y_2')$ in
 $\Gamma^{\zeta}$ we assign the morphism
$ ((y_1,\gamma,y_2),(y_1',\gamma,y_2'))$. This is well-defined, smooth, functorial, and the diagram is obviously commutative. Since the vertical arrows are weak equivalences, the 2-out-of-3 property of weak equivalences (Remark~\ref{Lemma: facts about weak equivalences}\,\ref{2-out-of-3}) shows that $F$ is a~weak equivalence.
 \end{proof}

\subsection*{Acknowledgements}
We would like to thank the anonymous referees for carefully reading the paper and providing valuable feedback.
Tim L\"uders is supported by a~scholarship of the Studienstiftung des deutschen Volkes e.V.

\pdfbookmark[1]{References}{ref}
\LastPageEnding

\end{document}